\documentclass[11pt]{amsart}
\usepackage{mathtools}
\usepackage{tabularx}
\usepackage{enumerate}
\usepackage{bbm, dsfont}
\usepackage{appendix}
\usepackage{amssymb,amsmath,amsfonts,eurosym,geometry,ulem,caption,color,setspace,comment,footmisc,pdflscape,subfigure,array,hyperref}

\usepackage{float}
\usepackage{mdframed}
\usepackage{lipsum}
\allowdisplaybreaks
\newtheorem{theorem}{Theorem}
\newtheorem{lemma}{Lemma}
\def\indic{{\rm {\large 1}\hspace{-2.3pt}{\large
l}}}
\def\N{\mathbb{N}}
\def\tr{{\rm tr}}
\def\rank{{\rm rank}}
\def\spn{{\rm span}}
\newtheorem{corollary}{Corollary}
\newtheorem{proposition}{Proposition}
\newtheorem{assumption}{Assumption}
\newtheorem{remark}{Remark}

\newcommand{\mat}[1]{\text{mat}}
\newcommand{\F}[1]{\text{F}}

\newcolumntype{L}[1]{>{\raggedright\let\newline\\arraybackslash\hspace{0pt}}m{#1}}
\newcolumntype{C}[1]{>{\centering\let\newline\\arraybackslash\hspace{0pt}}m{#1}}
\newcolumntype{R}[1]{>{\raggedleft\let\newline\\arraybackslash\hspace{0pt}}m{#1}}

\newmdtheoremenv{theo}{Theorem}
\usepackage{amsmath, amssymb, amsfonts}
\usepackage{amsmath}
\usepackage{listings}

\usepackage{setspace}

\usepackage{lmodern}
\usepackage{listings}
\usepackage{indentfirst}

\DeclarePairedDelimiter\norm{\lvert }{\rvert }
\usepackage{float}

\newcommand\op{\mathrm{op}}

\newcommand{\argmin}[1]{\underset{#1}{\text{argmin}}\;}

\usepackage[linesnumbered,ruled,vlined]{algorithm2e}
\RestyleAlgo{boxruled}
\LinesNumbered
\graphicspath{ {images/} }\graphicspath{ {images/} }

\usepackage{geometry}
\title[Panel Data Models with Approximately Low-rank Unobserved Heterogeneity]{Square-root Nuclear Norm Penalized Estimator for Panel Data Models with Approximately Low-rank Unobserved Heterogeneity}

\author{Jad Beyhum}
\author{Eric Gautier}
\address{Toulouse School of Economics, Universit\'e Toulouse Capitole, 21 all\'ee de Brienne, 31000 Toulouse, France}
\email{jad.beyhum@gmail.com, eric.gautier@tse-fr.eu}

\thanks{The authors acknowledge financial support from the grant ERC POEMH 337665. They are grateful to Thierry Magnac for helpful discussions}

\begin{document}
\maketitle

\begin{abstract}
This paper considers a nuclear norm penalized estimator for panel data models with interactive effects. The low-rank interactive effects can be an approximate model and  the rank of the best approximation  unknown and grow with sample size. The estimator is solution of a well-structured convex optimization problem and can be solved in polynomial-time. We derive rates of convergence, study the low-rank properties of the estimator, estimation of the rank and of annihilator matrices when the number of time periods grows with the sample size. We propose and analyze a two-stage estimator and prove its asymptotic normality. We can also use the baseline estimator as an initialization for any sequential algorithm. None of the procedures require knowledge of the variance of the errors.
\end{abstract}

\section{Introduction}
Panel data allow to estimate models with flexible unobserved heterogeneity using the fact that each individual is observed repeatedly. The high-dimensional statistics literature enables estimation in the presence of a high-dimensional parameter, provided that it has a low-dimensional structure. This paper studies a model that borrows from the two aforementioned strands of literature. We consider a linear panel data model with interactive effects 
\begin{equation}\label{model 1}
 Y_{it} = \sum_{k=1}^K\beta_kX_{kit} + \lambda_i^{\top}f_t+\Gamma^d_{it}+E_{it},\  \mathbb{E}[ E_{it}]=0,
\end{equation}
where $i\in\{1,...,N\}$ indices the individuals and $t\in\{1,...,T\}$ the time periods, $Y_{it}$ is the outcome, $X_{kit}$ is the $k^{\text{th}}$ regressor, $\beta\in \mathbb{R}^K$ is a vector of parameters, $\lambda_i$ and $f_t$ are vectors in $\mathbb{R}^r$ of factor loadings and factors, $\Gamma^d_{it}$ is a remainder which accounts for the fact that the usual interactive effects specification (when  $\Gamma^d_{it}=0$) can be an approximation, and $E_{it}$ is an error. Only $\beta$ is nonrandom.  Precise assumptions on the joint distribution of the right-hand side random elements is given later. Only the regressors and outcomes are available to the researcher. The regressors correspond to observed heterogeneity and the remaining right-hand side elements to unobserved heterogeneity. 
The interactive effects or statistical factor structure generalizes the usual individual plus time effects where 
$\lambda_i^{\top}f_t=c_i+d_t$. %It allows, for example, for group time effects of the form $d_{gt}$ for individuals in group $g$.  
One can think that $\lambda_i^{\top}f_t+\Gamma^d_{it}+E_{it}$ accounts for the contribution of regressors which are not available to the researcher but have an effect on the outcome if we believe these have a statistical factor structure plus remainder plus error term. In such a case, the error $E_{it}$ is a composite error which accounts for a linear combination of those coming from the missing regressors and the usual error from the long regression model which includes both observed and unobserved regressors. When the regressors and $\lambda_i^{\top}f_t+\Gamma^d_{it}$ are correlated, the least-squares estimator is inconsistent. This is a situation where we say that the regressors are endogenous or that there is an omitted variable bias. 
The specification is very flexible to model unobserved heterogeneity and can be broadly applied (see, \emph{e.g.}, \cite{gobillonmagnac} in the context of public policy evaluation). In matrix form, \eqref{model 1} becomes
\begin{equation}
\label{model 2}
                Y=\sum_{k=1}^{K} \beta_k X_k+ \Gamma^l+\Gamma^d+ E,
\end{equation}
where $Y,X_1,...,X_K, \Gamma^l,\Gamma^d$ and $E$ are random $N\times T$ matrices, $\Gamma_{it}^l=\lambda_i^{\top}f_t$, $\mathrm{rank}\left(\Gamma^l\right)=r$, and $\Gamma^d$ has small nuclear norm. The nuclear norm is the sum of  the singular values. We denote by $\Gamma=\Gamma^l+\Gamma^d$. %In this paper, $\beta$ is most of the time the parameter of interest and $\Gamma$ a nuisance. 
Many variations on model \eqref{model 1} have been considered and we name only a few. In \cite{chudik2015common,pesaran2006estimation} the regressors have a factor structure and $\beta$ can vary across individuals. In \cite{hansen2016factor,lu2016shrinkage} the number of regressors grows with the sample size. \cite{chudik2015common,moon2017dynamic} allow for lags of the outcome in \eqref{model 1}. In the setup where $\Gamma^d=0$ and $r$ is fixed and known, \cite{bai2009panel} analyses the least-squares estimator 
\begin{equation}
\label{Baiest}
                \left(\widehat{\beta}^B,\widehat{\Lambda}^B,\widehat{F}^B\right)\in\argmin{\substack{\beta\in\mathbb{R}^p\\ \Lambda^{\top}\Lambda\in\mathcal{D}_{rr},\ F^{\top}F=T I_r}}\left|Y-\sum_{k=1}^{K} \beta_k X_k-\Lambda F^{\top}\right|_2^2,
\end{equation}
where $|\cdot|_2$ is the $\ell_2$ norm of the vectorized matrix, $\Lambda$ (resp. $F$) is a $N\times r$ (resp. $T\times r$) matrix, $\mathcal{D}_{rr}$ the set of diagonal $r\times r$ matrices, and  $I_r$ the identity of size $r$. It is shown to be $\sqrt{NT}$-consistent and asymptotically normal when, among other things, the factors are strong. This means that 
$\Lambda^\top\Lambda/N$ converges in probability to a nonsingular matrix, hence the ratio of any singular value of $\Gamma^l$ and $\sqrt{NT}$ has a positive and finite limit in probability as $N$ goes to infinity and $T$ increases with $N$. \cite{moon2015linear} shows that using the same estimator with an upper bound on the number of factors leads to the same asymptotic properties. However, \eqref{Baiest} is a nonconvex optimization problem. For this reason, an iterative algorithm is used starting from an initial estimator which could be the least-squares estimator
 \begin{equation}
\label{LS}
                \widehat{\beta}^{LS}\in\argmin{\beta\in\mathbb{R}^p}\left|Y-\sum_{k=1}^{K} \beta_k X_k\right|_2^2
\end{equation}
or based on grid values. %When $\Gamma_{it}$ are interpreted as coming from missing regressors, it is natural that the regressors and $\Gamma_{it}$ are correlated so the regressors are endogenous. In this case, the initial estimator 
Using $\widehat{\beta}^{LS}$ can be problematic because it can be inconsistent and thus far away from a (the?) global minimum of \eqref{Baiest}. \cite{Hsiao2} analyses the asymptotic properties of $m^{\rm th}$ iterates of one of \cite{bai2009panel}'s iterative algorithm treating them as an estimator. It is found that it can be consistent when adding several additional assumptions among which the consistency of $\widehat{\beta}^{LS}$ if used as an initialization. It is argued in Remark 4.2 in \cite{Hsiao} that iterative algorithms are consistent if the initialisation is by a consistent estimator. %Because least-squares can fail, it is important to provide alternatives. 
Corollary 1 in \cite{moon2018nuclear} give a condition on the rate of convergence of a preliminary estimator so that an iterative algorithm in the spirit of those proposed by \cite{bai2009panel} is asymptotically equivalent to \cite{bai2009panel}'s theoretical estimator.

The tools in this paper are related to those used in matrix completion. There, the problem consists in estimating the unobserved entries of a low-rank matrix from an observed subset of its entries, sometimes with additive noise (see, \emph{e.g.}, \cite{candes2011tight, candes2010power,klopp2014noisy,klopp2011high,koltchinskii2011nuclear,recht2011simpler,recht2010guaranteed,rohde2011estimation}). 
%\blue{(\cite{klopp2011high} pourrait apparaitre autrement que dans cette liste non ordonnee d'autant qu'ils etudient la variance inconnue comme ici)} \red{Ok, j'écrirai quelque chose} 
The usual $\ell_0$ and $\ell_1$-norms are replaced by the rank and nuclear norm, soft and hard thresholding are carried on the singular values.  These methods have recently been used in econometrics (see in particular \cite{atheyimbens,BaiNgHD,galichon}).  The problem in this paper differs in that we observe all the entries of the matrices $Y$ and $X_1,\dots,X_K$ but none of $\Gamma+E$ and both $\Gamma$ and $E$ are  random. 

The iterative procedures in \cite{bai2009panel} could yield a local minimum while the theoretical properties are for the global minimum. In contrast, the estimators in \cite{moon2018nuclear} and in this paper involve convex programs for which convergence to a global minimum is achieved in polynomial time. The additional novelties of this paper are as follows. This paper considers a square-root nuclear norm penalized estimator (see \cite{belloni2011square} for the Lasso), where the sum of squared residuals is replaced by its square-root. It can be viewed as the estimator in \cite{moon2018nuclear} using a data-driven penalty level so it is directly implementable by the researcher and does not require an additional diverging multiplicative factor which can result in over-penalization. % and is useful in finite samples. 
We provide a straightforward iterative algorithm to compute the estimator. Our results do not rely on conditioning on realizations of $\Gamma$ and we state the conditions on the joint distribution of $\Gamma$ and the regressors. We allow the interactive effect to be an approximate model and hence many non-strong factors (see \cite{pesaran2015time}) via the additional term $\Gamma^d$. The rank of $\Gamma^l$ is treated as random and can grow with the sample size and be unknown. We obtain low-rank oracle type inequalities for various loss functions and results on the rank of our estimator of $\Gamma$, introduce a thresholded estimator which can be used to estimate the rank of $\Gamma^l$ as well as projectors on the vector spaces spanned by the factors and factor loadings which we analyze theoretically. We also obtain rates of convergence for the estimation of $\beta$. These results do not rely on a strong-factor assumption.  %(see \cite{bai2009panel} and \cite{moon2015linear}) 
Finally, we propose a two-stage estimator and show its asymptotic normality. Based on our procedure and result on the estimation of the rank of $\Gamma^l$, we can proceed as analyzed in \cite{moon2018nuclear}
and use an iterative algorithm as a second stage. %In all procedures, the variance of the errors if homoskedastic is always assumed unknown.

\section{Preliminaries}
%\subsection{Notations} 
$\N$ denotes the positive integers, $\N_0$ denotes $\N\cup\{0\}$. For $a\in \mathbb{R}$, we set $a_+=\max(a,0)$ and, for $a>0$, $a/0=\infty$. $\left\{\mu_{N}\right\}$ denotes a numerical sequence of generic term $\mu_N$. $\mathcal{M}_{NT}$ is the set of matrices with real coefficients of size $N\times T$. 
%\blue{(si on garde) $L$ denotes the lag operator}. 
The transpose of $A\in \mathcal{M}_{NT}$ is $A^{\top}$, its trace is $\mathrm{tr}(A)$, and its rank is $\rank(A)$. For $A\in \mathcal{M}_{NT}$ and $v\in \mathbb{R}^{NT}$, $\text{vec}(A)$ is obtained by stacking the columns of $A$ and $\text{mat}(v)$ is the unique matrix in $\mathcal{M}_{NT}$ such that $v=\text{vec}\left(\text{mat}(v)\right)$. Matrices are denoted using capital letters and their vectorization using lowercase letters. 
The $k^{\text{th}}$ singular value of  $A\in \mathcal{M}_{NT}$ (arranged in decreasing order and repeated according to multiplicty) is $\sigma_k(A)$. % and $\mathrm{rank}(A)$ is its rank. 
$A=\sum_{k=1}^{\mathrm{rank}(A)}\sigma_k(A)u_k(A)v_k(A)^\top$ is the singular value decomposition of $A$, where $\left\{u_k\left(A\right)\right\}_{k=1}^{\mathrm{rank}\left(A\right)}$ is a family of orthonormal vectors of $\mathbb{R}^N$ and $\left\{v_k\left(A\right)\right\}_{k=1}^{\mathrm{rank}\left(A\right)}$ of $\mathbb{R}^T$. The scalar product is $\left\langle A,B\right\rangle=\mathrm{tr}\left(A^{\top}B\right)$. The $\ell_2$-norm (or Frobenius norm) is $\norm{A}_2^2=\left\langle A,A\right\rangle=\sum_{k=1}^{\mathrm{rank}(A)}\sigma_k(A)^2$, the nuclear norm is $\norm{A}_*=\sum_{k=1}^{\mathrm{rank}(A)}\sigma_k(A)$, and the operator norm is $|A|_{\op}=\sigma_1(A)$. %=\max\limits_{h\in\mathbb{R}^T\ \text{s.t.}\ \norm{h}_2=1} \norm{Ah}_2$. 
$P_{u(A)}$ and $P_{v(A)}$ are the orthogonal projectors onto $\spn\{u_1(A),\dots,u_{\rank(A)}(A)\}$ and $\spn\{v_1(A),\dots,v_{\rank(A)}(A)\}$ and $M_{u(A)}$ and $M_{v(A)}$ %the projectors 
onto the orthogonal complements. For $A,\Delta\in\mathcal{M}_{NT}$, we define $\mathcal{P}_{A}(\Delta)=\Delta-M_{u(A)}\Delta M_{v(A)}$ and 
$\mathcal{P}_{A}^{\perp}(\Delta)=  M_{u(A)} \Delta  M_{v(A)}$. 
\noindent Recall that, if $\widetilde{\Delta}\in\mathcal{M}_{NT}$ (see lemma 2.3 and 3.4 in \cite{recht2010guaranteed} for \eqref{orth}-\eqref{add}), 
\begin{align}
\mathcal{P}_{A}(\Delta)&=
M_{u(A)}\Delta P_{v(A)}+P_{u(A)}\Delta\label{eP},\\
\rank\left(\mathcal{P}_{A}( \Delta)\right)&\le 2\min\left(\rank\left( \Delta\right),\rank(A)\right),\label{eq:rank}\\
\left\langle\mathcal{P}_{A}( \Delta),\mathcal{P}_{A}\left(\widetilde{\Delta}\right)\right\rangle&=\left\langle\mathcal{P}_{A}( \Delta),\widetilde{\Delta}\right\rangle,\\
\left\langle\mathcal{P}_{A}( \Delta),\mathcal{P}_{A}^{\perp}( \Delta)\right\rangle&=0,\label{orth} \\
\left|A+\mathcal{P}_{A}^{\perp}( \Delta)\right|_*&=
\left|A\right|_*+\left|\mathcal{P}_{A}^{\perp}( \Delta) \right|_*\label{add}.
\end{align}
The cone %(see \cite{koltchinskii2011nuclear}) 
$C_{A,c}=\left\{\Delta\in\mathcal{M}_{NT}:\ 
\left|\mathcal{P}_{A}^{\perp}\left(\Delta\right)\right|_*
\le c\left|\mathcal{P}_{A}\left(\Delta\right)\right|_{*}\right\}$, defined for $A\in\mathcal{M}_{NT}$ and $c>0$, 
is important for our analysis. We use the notations $A^l,A^d\in\mathcal{M}_{NT}$ for two components such that $A=A^l+A^d$. The role of and assumptions on $A^l$ and $A^d$ will be clear from the text. $A^l$ stands for a ``low-rank" (the rank can diverge with sample size) component with a large operator norm while $A^d$ is a small 
remainder term.

We denote by $P_X$ (resp. $M_X$) the orthogonal projector on the vector space spanned by $\left\{X_k\right\}_{k=1}^K$ (resp. on its orthogonal) and 
 $X=\left(x_1,\dots, x_K\right)$. %The number of regressors $K$ does not increase with the sample size. 
 We consider an asymptotic where $N$ goes to infinity and $T$ is a function of $N$ that goes to infinity when $N$ goes to infinity. The probabilistic framework consists of a sequence of data generating processes (henceforth DGPs) that depend on $N$.  We write that an event occurs w.p.a. $1$ ("with probaility approaching $1$") when its probability converges to $1$ as $N$ goes to infinity. All limits are when $N$ goes to infinity. We denote convergence in probability and in distribution by respectively $\xrightarrow{\mathbb{P}}$ and $\xrightarrow{d}$. 
We allow the researcher to apply annihilator matrices $M_u$ (to the left) and $M_v$ (to the right) on both sides of \eqref{model 2} and still denote by $Y,X_1,\dots,X_K,\Gamma^l,\Gamma^d,E$ the transformed matrices. She can apply a within-group (or first-difference or Helmert) transform on the left to annihilate individual effects and a similar on the right to annihilate time effects. %, two matrices are required to annihilate group specific time effects. 
This is important if the researcher thinks there are individual and time effects and there could be additional interactive effects and wants to avoid relying on penalization to figure out that there are classical individual and time effects. %These annihilators can be estimate. %d and used to annihilate low-rank components of the regressors with large operator norm and obtain transformed regressors with smaller norm. 
The regressors can be transformations of the baseline regressors as in Section \ref{s:transform} to ensure their operator norm is not too large. %, a feature sometimes useful in the analysis. 
We do not write these transformations explicitly to simplify the exposition.

\section{First-stage estimator}
\label{sec:2}
%\subsection{Definition}
The estimator is defined, for $\lambda >0$, as 
   \begin{equation} \label{estimator2}
 \left(\widehat{\beta},\widehat{\Gamma}\right) \in \argmin{\beta \in \mathbb{R}^K,\  \Gamma \in \mathcal{M}_{NT}
} \frac{1}{\sqrt{NT}}\left\vert  Y-\sum_{k=1}^K\beta_k X_k-\Gamma  \right \vert_2
+\frac{\lambda}{NT} \norm{\Gamma}_*.
 \end{equation}
The nuclear norm plays the role of the $\ell_1$-norm in the Lasso estimator. It yields low-rank solutions, that is  a sparse vectors of singular value of $\widehat{\Gamma}$ (see Proposition \ref{l1} for a formal result). The nuclear norm penalization is the convex relaxation of a penalization involving $\rank\left(\Gamma\right)$. 
This estimator can be viewed as a type of square-root Lasso estimator of \cite{belloni2011square} for parameters which are matrices. As for the square-root Lasso, the $\ell_2$-norm is not squared in \eqref{estimator2} which implies that we do not need to know the variance of $E_{it}$ when these are iid to choose $\lambda$. For large classes of DGP, this choice will be canonical. 
 \begin{proposition}\label{equiv}  
A solution $\left(\widehat{\beta},\widehat{\Gamma}\right)$ of \eqref{estimator2} is such that
 \begin{equation*}
\widehat{\Gamma}
\in 
\arg \min_{\Gamma\in\mathcal{M}_{NT}} \frac{1}{\sqrt{NT}}\left\vert M_X\left( Y-\Gamma\right)  \right \vert_2
+\frac{\lambda}{NT} \norm{\Gamma}_*.
\end{equation*}
 \end{proposition}
Let us provide another interpretation for this estimator. In the least-squares problem  $\min\left\vert  Y-\sum_{k=1}^K\beta_k X_k-\Gamma  \right \vert_2^2$, even if $\Gamma$ were given, estimation of $\beta$ is an inverse problem. 
For this reason, properties of the design matrix matter for estimation. While, if $\beta$ were given, estimation of $\Gamma$ would not be an inverse problem, it is when $\beta$ is unknown. Indeed, applying $M_X$ to \eqref{model 2}, we obtain $M_X(Y)=M_X(\Gamma)+M_X(E)$. Because $\Gamma$ appears via $M_X(\Gamma)$, estimation of $\Gamma$ is an inverse problem with correlated errors which can be correlated with $M_X(\Gamma)$ via $M_X$. The trace of the covariance operator of the error term is $\mathbb{E}\left[\left|M_X(E)\right|_2^2\right]$ and diverges with $N$.  We will see that $\left|M_X(E)\right|_2/\sqrt{NT}$ can converge to the standard error of the entries of $E$ when these are iid mean zero with finite variance. The nonstandard framework here is that $X$, hence the operator, and the ``parameter'' $\Gamma$ are random. Also $M_X$ is not invertible. But estimation of $\Gamma$ becomes feasible under shape restrictions. This paper considers a generalization of the restriction that $\Gamma$ has low rank by allowing for approximately low-rank matrices.
 
For $u\ge0$, $u=\min_{\sigma>0}\left\{\frac{\sigma}{2}+\frac{1}{2\sigma}u^2\right\}$ and the minimum is attained at $\sigma=u$ if $u>0$ or using minimizing sequences going to 0 if $u=0$. Thus, any solution $\left(\widehat{\beta},\widehat{\Gamma}\right)$ of \eqref{estimator2} is solution of 
%\begin{equation} \notag \left(\widehat{\beta},\widehat{\Gamma}\right)\in \argmin{\beta\in\mathbb{R}^K,\Gamma\in \mathcal{M}_{NT}} \frac{1}{\sqrt{NT}}\left\vert Y-\sum_{k=1}^K\beta_kX_k-\Gamma \right\vert_2+\frac{\lambda}{NT}\left|\Gamma\right|_*
%\end{equation}
\begin{equation}\label{optiglob}
\left(\widehat{\beta},\widehat{\Gamma},\widehat{\sigma}\right)\in \argmin{\beta\in\mathbb{R}^k,\Gamma\in\mathcal{M}_{NT},\sigma>0}\sigma+ \frac{1}{\sigma NT}
\left\vert Y-\sum_{k=1}^K\beta_kX_k-\Gamma \right\vert_2^2+\frac{2\lambda}{NT}\left|\Gamma\right|_*
\end{equation}
and 
\begin{equation}\label{sigmahat}
\widehat{\sigma}=\frac{1}{\sqrt{NT}}\left\vert Y-\sum_{k=1}^K\widehat{\beta}_kX_k-\widehat{\Gamma}\right\vert_2.
\end{equation}
This objective function in \eqref{optiglob} has the advantage that the new objective function only has one nonsmooth convex function in $(\beta,\Gamma)$: the nuclear norm.  Because $f(x,y)=x^2/y$ is convex on the domain $\{(x,y)\in\mathbb{R}^2|y>0\}$, the objective function is convex in $(\beta,\Gamma,\sigma)$. This formulation is analogous to the concomitant Lasso or scaled-Lasso for linear regression (see \cite{owen,sun2012scaled}). 

\subsection{First-order conditions and consequences}
%In this section we use similar arguments as those recalled page 112 in \cite{giraud2014introduction}.
 The formulation is used in Section \ref{sec:ca} for implementation of our estimator. It is also useful to obtain by subdifferential calculus the first order-conditions of program \eqref{estimator2}. Indeed, the differential with respect to $\beta_k$ at $(\beta,\Gamma,\sigma)$ on the domain (hence $\sigma>0$) is, for $k=1,\dots,K$,  
   %  We introduce 
     %$y=\text{vec}(Y)$, , $\gamma=\text{vec}(\Gamma)$, $e=\text{vec}(E)$ and $\widehat{\gamma}=\text{vec}\left(\widehat{\Gamma}\right)$.
%\subsection{Computational aspects}
\begin{equation}\label{eq:grad}
-\frac{2}{\sigma NT}\left\langle X_k,Y-\sum_{k=1}^K\beta_kX_k-\Gamma\right\rangle
\end{equation}
and the subdifferential with respect to $\Gamma$ at $(\beta,\Gamma,\sigma)$ (see (2.1) in \cite{koltchinskii2011nuclear}) is 
\begin{equation}\label{eq:sub}
\left\{-\frac{2}{\sigma NT}\left(Y-\sum_{k=1}^K\beta_kX_k-\Gamma\right)+\frac{2\lambda}{NT}Z,\ Z=\sum_{k=1}^{\rank(\Gamma)}u_k(\Gamma)v_k(\Gamma)^{\top}+M_{u(\Gamma)}W M_{v\left(\Gamma\right)},\ |W|_{\op}\le1\right\},
\end{equation}
in particular 
$\left|Z\right|_{\op}\le 1$ and $\left\langle\Gamma,Z\right\rangle=\left|\Gamma\right|_*$. 
Due to \eqref{sigmahat}, if $\widehat{\sigma}=0$ then clearly $\widehat{\beta}$ is the least-squares estimator which minimizes $\left|Y-\sum_{k=1}^K\beta_kX_k-\widehat{\Gamma}\right|_2^2$. Else, setting \eqref{eq:grad} to 0 at $\left(\widehat{\beta},\widehat{\Gamma},\widehat{\sigma}\right)$ yields the same conclusion. 
Hence, if $X^\top X$ is positive definite, we have %the usual formula
\begin{equation}\label{OLS}
\widehat{\beta}=\left(X^\top X\right)^{-1}X^\top (y-\widehat{\gamma}).
\end{equation}
Because, if $\widehat{\sigma}>0$, 
$0$ belongs to the set defined in \eqref{eq:sub} at $\left(\widehat{\beta},\widehat{\Gamma},\widehat{\sigma}\right)$, 
there exists $\widehat{W}\in\mathcal{M}_{NT}$ and   
$\widehat{Z}=\sum_{k=1}^{\rank(\widehat{\Gamma})}u_k\left(\widehat{\Gamma}\right)v_k\left(\widehat{\Gamma}\right)^{\top}+M_{u\left(\widehat{\Gamma}\right)}\widehat{W} M_{v\left(\widehat{\Gamma}\right)}$ such that $\left|\widehat{W}\right|_{\op}\le1$ and 
%\begin{equation*}
$Y-\sum_{k=1}^K\widehat{\beta}_kX_k-\widehat{\Gamma}=\lambda\widehat{\sigma}\widehat{Z}$, 
%\end{equation*}
hence, for all $k=1,\dots,K$, $\left\langle X_k,\widehat{Z}\right\rangle=0$, thus $M_{X}\left(\widehat{Z}\right)=\widehat{Z}$ and
\begin{equation}\label{eq:apend}
Y-\sum_{k=1}^K\widehat{\beta}_kX_k-\widehat{\Gamma}=M_{X}\left(Y-\widehat{\Gamma}\right)=\lambda\widehat{\sigma}\widehat{Z}.
\end{equation}
Again, due to \eqref{sigmahat}, if $\widehat{\sigma}=0$ then \eqref{eq:apend} also holds.
As a consequence, we have 
\begin{align*}
&\widehat{\sigma}=\frac{1}{\sqrt{NT}}\left\vert M_X\left( Y-\widehat{\Gamma}\right)  \right \vert_2
%,\\&\widehat{\Gamma}\in \argmin{\Gamma\in\mathcal{M}_{NT}} \frac{1}{\sqrt{NT}}\left\vert M_X\left( Y-\Gamma\right)  \right \vert_2
%+\frac{\lambda}{NT} \norm{\Gamma}_*.
\end{align*}
and any solution $\left(\widehat{\beta},\widehat{\Gamma}\right)$ of \eqref{estimator2} is also solution of 
\begin{equation} \label{estimator3}
 \left(\widehat{\beta},\widehat{\Gamma}\right) \in \argmin{\beta \in \mathbb{R}^K,\  \Gamma \in \mathcal{M}_{NT}
} \frac{1}{NT}\left\vert  Y-\sum_{k=1}^K\beta_k X_k-\Gamma  \right \vert_2^2
+\frac{2\lambda\widehat{\sigma}}{NT} \norm{\Gamma}_*.
 \end{equation}
 So $ \left(\widehat{\beta},\widehat{\Gamma}\right)$ given by \eqref{estimator2} is a solution to a type of matrix Lasso estimator %considered in \cite{moon2018nuclear} 
 with data-driven penalty
 $\lambda\widehat{\sigma}\norm{\Gamma}_*/NT$. 
The estimator in \cite{moon2018nuclear} corresponds to \eqref{estimator3} without the data-driven  $\widehat{\sigma}$. 
\begin{remark} %The penalty involving the nuclear norm is used to enforce a low-rank estimator of $\Gamma$.
Due to the nuclear norm, \eqref{eq:apend} %, hence  we do not know $\rank(\Gamma)$ and rely on a
 and the expression of $\widehat{Z}$ yield 
$$\left(Y-\sum_{k=1}^K\widehat{\beta}_kX_k\right)M_{v\left(\widehat{\Gamma}\right)}=\lambda\widehat{\sigma}M_{u\left(\widehat{\Gamma}\right)}\widehat{W} M_{v\left(\widehat{\Gamma}\right)}$$ which, unlike \cite{bai2009panel}, is not zero. Applying the annihilator $M_{u\left(\widehat{\Gamma}\right)}$ does not change this. %This is because, unlike \cite{bai2009panel}, 
\end{remark}

\subsection{Computational aspect}\label{sec:ca}
Based on \eqref{optiglob}, where the objective function is convex,  we can minimize  iteratively over $\beta$, $\Gamma$, and $\sigma$: start from $\left(\beta^{(0)},\Gamma^{(0)},\sigma^{(0)}\right)$ and repeat %, for $t\in\N_0$ 
until convergence,   
\begin{enumerate}[\textup{(}1\textup{)}] 
\item\label{s1} $\beta^{(t+1)}$ is obtained by least-squares minimizing $\left\vert Y-\sum_{k=1}^K\beta_kX_k-\Gamma^{(t)}\right\vert_2^2$,
\item\label{s2} Setting $Z^{(t+1)}=Y-\sum_{k=1}^K\beta_k^{(t+1)}X_k$, $\Gamma^{(t+1)}$ is obtained by solving the matrix Lasso
$$\min_{\Gamma}  \left\vert Z^{(t+1)}-\Gamma \right\vert_2^2+2\lambda\sigma^{(t)}\left|\Gamma\right|_*,$$
\emph{i.e.} applying soft-thresholding to the singular value decomposition (henceforth SVD)
$$\Gamma^{(t+1)}=\sum_{k=1}^{\min(N,T)}\left(\sigma_k\left(Z^{(t+1)}\right)-\lambda\sigma^{(t)}\right)_+u_k\left(Z^{(t+1)}\right)v_k\left(Z^{(t+1)}\right)^{\top},$$
\item\label{s3} $\sigma^{(t+1)}= \left\vert Z^{(t+1)}-\Gamma^{(t+1)}\right\vert_2/\sqrt{NT}$.
\end{enumerate}
\begin{remark} 
The estimator in \cite{moon2018nuclear} can be obtained by repeating (1) and (2) for a fixed value of $\sigma^{(t)}$. $\lambda_N\sigma^{(t)}$ corresponds to $\sqrt{NT}\Psi_{NT}$ in their notations and they assume $1/\left(\Psi_{NT}\sqrt{\min(N,T)}\right)+\Psi_{NT}\to0$ to circumvent the unavailability of an upper bound on the variance of the errors.   
\end{remark}
\begin{remark} Various numerical approximations of the theoretical estimator \eqref{Baiest} are discussed in \cite{bai2009panel}. 
The method page 1237 in \cite{bai2009panel} considers iterates step (1) and a modified step (2) where 
 $\lambda=0$ and under the restriction that $\rank(\Gamma)=r$, from which we extract the factor and factor loadings. The second step corresponds to hard-thresholding the SVD of $Z^{(t+1)}$ to keep only the part corresponding to the $r$ largest singular values. It is argued that the approach from iterating (a) augmented least-squares given factors
 \begin{equation*}
\label{LSA}
                \left(\beta^{(t+1)},\Lambda^{(t+1)}\right)\in\argmin{\beta\in\mathbb{R}^p,\ \Lambda\in\mathcal{M}_{Nr}}\left|Y-\sum_{k=1}^{K} \beta_k X_k-\Lambda \left(F^{(t)}\right)^\top\right|_2^2
\end{equation*}
or, by partialling out and the fact that $M_{v\left(\Gamma^{(t)}\right)}$ is also the projector onto the orthogonal of the columns of $F^{(t)}$, 
 \begin{equation*}
%\label{LSA1}
\beta^{(t+1)}\in\argmin{\beta\in\mathbb{R}^p}\left|\left(Y-\sum_{k=1}^K\beta_kX_k\right)M_{v\left(\Gamma^{(t)}\right)}\right|_2^2
\end{equation*}
 and (b) PCA to obtain $F^{(t+1)}$ %the $r$ common factors 
 is less numerically robust. 
  \end{remark}

\section{Results}
\subsection{Error bound on the estimation of $\Gamma$ and $\left|M_X(E)\right|_2/\sqrt{NT}$}\label{sec:3}
The result of this section are upper bounds on the error made by estimating $\Gamma$ and $\left|M_X(E)\right|_2/\sqrt{NT}$ by $\left(\widehat{\Gamma},\widehat{\sigma}\right)$. They hold without assumption. A key quantity is the compatibility constant (see \cite{buhlmann2011statistics} for the high-dimensional linear regression). It is defined, for all realizations of $X$ and $A\in\mathcal{M}_{NT}$, by
$$\kappa_{A,c}  =\inf_{ \Delta \in C_{A,c}:\ \Delta\ne 0} \frac{\sqrt{2\mathrm{rank}\left(A\right)}\norm{M_X(\Delta)}_{2}}{\left|\mathcal{P}_{A}\left(\Delta\right)\right|_*}.$$
\begin{remark}
A few remarks are in order. First, if $X=0$, we have $M_X(\Delta)=\Delta$. Second, the denominator of the ratio cannot be 0 because, for $\Delta\in C_{A,c}$,
$\left|\Delta\right|_*\le (1+c) \left|\mathcal{P}_{A}\left(\Delta\right)\right|_*$, hence the function of $\Delta$ in the infimum is continuous. Third, because the ratio  involves two linear operators, the infimum is the same if we restrict $\Delta$ to have norm 1 and the intersection with the cone is compact. Hence, the infimum is a minimum. 
Fourth, for all $A\in\mathcal{M}_{NT}$ and $c>0$, the minimum is the limit of minima over finite sets so it is a measurable function of $X$. Fifth, we work with $\kappa_{\widetilde{\Gamma},c}$ for a random $\widetilde{\Gamma}$ which depends on the random $\Gamma$ and $X$ via $\kappa_{\widetilde{\Gamma},c}$ itself and we allow $\Gamma$ and $X$ to be dependent. We  make a slight abuse of notations and consider that $\kappa_{\widetilde{\Gamma},c}$ is a measurable function of $\widetilde{\Gamma}$ and $X$. In practice, it is a measurable lower bound on it for every fixed $\widetilde{\Gamma}\in\mathcal{M}_{NT}$ and $X$ in the support of the corresponding random matrix.
\end{remark}
\begin{remark}
When $X=0$ one has, for all $A\in\mathcal{M}_{NT}$ and $c>0$, $\kappa_{A,c}  \ge1$.
 \end{remark}
\begin{proposition}\label{pkappa} The following lower bounds hold
\begin{equation}\label{comp}\kappa_{A,c} \ge\min_{\Delta \in C_{A,c}:\ \Delta\ne 0} \frac{\norm{M_X(\Delta)}_{2}}{\left|\mathcal{P}_{A}\left(\Delta\right)\right|_2}\ge\min_{\Delta \in C_{A,c}:\ \Delta\ne 0}\frac{\norm{M_X(\Delta)}_{2}}{\left|\Delta\right|_2}.
\end{equation}
\end{proposition}
The quantity in the middle is the restricted eigenvalue (see \cite{koltchinskii2011nuclear}). The smaller one is used in \cite{moon2018nuclear}. These constants are essential elements in the upper bounds of Theorem \ref{orac} and Proposition \ref{l2} and further discussed below. Throughout the rest of the paper, $\rho\in(0,1)$. Define
\begin{align*}
&c\left(\rho,\widetilde{\rho}\right)=\frac{1+\rho+\widetilde{\rho}}{1-\rho},\ d\left(\rho,\widetilde{\rho}\right)=\max\left(1+\widetilde{\rho},\rho \left(1+c\left(\rho,\widetilde{\rho}\right)\right)\right),\ e\left(\rho,\widetilde{\rho}\right)=d\left(\rho,\widetilde{\rho}\right)+\rho \left(1+c\left(\rho,\widetilde{\rho}\right)\right),\\
&\theta_\infty\left(\widetilde{\Gamma},\rho,\widetilde{\rho}\right)=2\left(1-   \left(\frac{d\left(\rho,\widetilde{\rho}\right)\sqrt{2\mathrm{rank}\left(\widetilde{\Gamma}\right)}\lambda}{\sqrt{NT}\kappa_{\widetilde{\Gamma},c\left(\rho,\widetilde{\rho}\right)}}\right)^2\right)_+^{-1}e\left(\rho,\widetilde{\rho}\right),\\
&\theta\left(\rho,\widetilde{\rho}\right)=\inf_{\widetilde{\Gamma}\in\mathcal{M}_{NT}}
\max\left(\theta_\infty\left(\widetilde{\Gamma},\rho,\widetilde{\rho}\right)\frac{\lambda\mathrm{rank}\left(\widetilde{\Gamma}\right)\norm{M_X(E)}_2}{\sqrt{NT}\kappa_{\widetilde{\Gamma},c\left(\rho,\widetilde{\rho}\right)}^2},\frac{1}{\widetilde{\rho}}\left|\Gamma-\widetilde{\Gamma}\right|_*\right),\\
&\theta_*(\rho)=\inf_{\widetilde{\rho}>0}\left(1+c\left(\rho,\widetilde{\rho}\right)\right)\theta\left(\rho,\widetilde{\rho}\right),\ 
\theta_\sigma(\rho)=\inf_{\widetilde{\rho}>0}d\left(\rho,\widetilde{\rho}\right)\theta\left(\rho,\widetilde{\rho}\right).
\end{align*}
\begin{theorem}\label{orac} 
On the event $\mathcal{E}=\left\{\rho\lambda\left|M_X(E)\right|_2/\sqrt{NT}\ge\norm{M_X(E)}_{\op}\right\}$, we have
\begin{align}
&\left|\widehat{\Gamma}-\Gamma\right|_* \le 2\theta_*(\rho),\label{OracleBNF} \\
&\left \vert\widehat{\sigma}-\frac{1}{\sqrt{NT}}\left\vert M_X\left(E\right)\right \vert_2\right \vert \le \frac{2\lambda}{NT}
\theta_\sigma(\rho).\label{Oraclevar} 
\end{align}
\end{theorem}
Note that $\theta_*(\rho)\le\theta_\sigma(\rho)/\rho$. For example, we can take $\widetilde{\rho}=1$ and $\rho=2/5$, in which case $c\left(\rho,\widetilde{\rho}\right)=4$, $d\left(\rho,\widetilde{\rho}\right)=2$, $e\left(\rho,\widetilde{\rho}\right)=4$, $\theta_*(\rho)=5\theta\left(\rho,\widetilde{\rho}\right)$, and $ \theta_\sigma(\rho)=2\theta\left(\rho,\widetilde{\rho}\right)$. 
We state a more general result to allow the theory to handle the case where $\rho$ is close to 1 which allows a smaller $\lambda$ (what matters is the product $\rho\lambda$ in the definition of $\mathcal{E}$ and $\rho\psi_N$ in Assumption \ref{Pen}) which we find works well in small samples. If we write $\Gamma=\Gamma^l+\Gamma^d$ and take $\widetilde{\Gamma}=\Gamma^l$ in the maximum in the expression of $\theta(\rho,\widetilde{\rho})$, we obtain
\begin{equation}\label{ubtheta}
\theta(\rho,\widetilde{\rho})\le\max\left(\theta_\infty\left(\Gamma^l,\rho,\widetilde{\rho}\right)\frac{\lambda\mathrm{rank}\left(\Gamma^l\right)\norm{M_X(E)}_2}{\sqrt{NT}\kappa_{\Gamma^l,c\left(\rho,\widetilde{\rho}\right)}^2},\frac{1}{\widetilde{\rho}}\left|\Gamma^d\right|_*\right).
\end{equation}
Under the premises of Proposition \ref{sigmacons} \eqref{ub1}, the quantity $\norm{M_X(E)}_2/\sqrt{NT}$ is consistent and we can obtain an upper bound like \eqref{ubtheta} where it is replaced by a constant $\sigma$.  
We obtain a tight bound if $\Gamma^l$ and $\Gamma^d$ in the decomposition $\Gamma=\Gamma^l+\Gamma^d$ are such that $\theta_\infty\left(\Gamma^l,\rho,\widetilde{\rho}\right)\mathrm{rank}\left(\Gamma^l\right)/\kappa_{\Gamma^l,c\left(\rho,\widetilde{\rho}\right)}^2$ and in particular $\rank\left(\Gamma^l\right)$ is small and $\left|\Gamma^d\right|_*$ is small. But $\Gamma^d$ could have high-rank. The decomposition  $\left(\Gamma^l,\Gamma^d\right)$ which realizes the tradeoff depends on $\Gamma$ (via their sum) and $M_X$ which appears in the definition of $\kappa_{\Gamma^l,c\left(\rho,\widetilde{\rho}\right)}$. Theorem 1 and Proposition \ref{l2} show that $\widehat{\Gamma}$ performs, up to a multiplicative constant, as well as an oracle who would know $\Gamma$ and choose the best misspecified model to keep the number of incidental parameters moderate while incurring a bias which is not too large. The term  involving $(\cdot)_+^{-1}$ in the definition of $\theta_\infty\left(\widetilde{\Gamma},\rho,\widetilde{\rho}\right)$ could be $\infty$ if $\kappa_{\widetilde{\Gamma},c\left(\rho,\widetilde{\rho}\right)}$ is too small. It appears because we do not know the variance of the errors or use a sequence of penalties that diverge faster than necessary. Finally, a smaller constant $c\left(\rho,\widetilde{\rho}\right)$ implies a smaller cone and a larger $\kappa_{\widetilde{\Gamma},c\left(\rho,\widetilde{\rho}\right)}$.

Let us comment the first term in the maximum in the right-hand side of \eqref{ubtheta}. $\rank\left(\Gamma^l\right)$ plays the same role as the number of nonzeros for estimation of sparse vectors of coefficients in linear regression models. 
The other key ingredient is $\kappa_{\Gamma^l,c\left(\rho,\widetilde{\rho}\right)}^2$. Because it appears in the denominator of an upper bound, it is desirable to have it as large as possible. The compatibility constant is the sharpest of the three quantities in \eqref{comp}. To gain insight, we make an analogy with linear regression. Denoting by $X$ the design matrix and $N$ the sample size, the rate of estimation of the vector of coefficients depends on the largest eigenvalue of $(X^{\top}X/N)^{-1}$ in a numerator, or the smallest eigenvalue of $X^{\top}X/N$ in a denominator (the square of the smallest singular value of $X/\sqrt{N}$). Sharper constants can be used when the  vector is sparse. Because of the $\ell_1-$norm in the Lasso, the difference between the estimator and the sparse vector belongs  to a cone with probability close to 1. 
As a result, the smallest singular value which, by the Courant-Fisher theorem, is solution to a minimization problem, can be replaced by a minimization on the cone rather than on the whole space. This is important in high-dimensions because the minimum singular value is zero when the dimension is larger than the sample size. In his paper, the smaller quantity in \eqref{comp} restricts the whole space to $C_{A,c}$ in the minimization problem defining the smallest singular value of the operator $M_X$.  Recall that, without restriction, the minimum is 0 because $M_X$ is not invertible. The relevant cone is $C_{\Gamma^l,c}$. It contains $\Gamma^l$ but also, if $\Gamma^d=0$, we show that it contains, with probability close to 1, $\Delta$ defined as $\widehat{\Gamma}-\Gamma$, for all minimizer $\widehat{\Gamma}$. The definition of the compatibility constant yields 
$$\left|\Delta\right|_*\le (1+c)\frac{\sqrt{2\rank\left(\Gamma^l\right)}}{\kappa_{\Gamma^l,c}}\left|M_X(\Delta)\right|_2.$$  
This allows to relate the error in terms of a loss involving the nuclear norm to the loss derived from the least-squares criterion in the optimization program of Proposition \ref{equiv}. The restricted eigenvalue replaces $\Delta$ in the denominator by a type of projection $\mathcal{P}_{A}\left(\Delta\right)$ of $\Delta$ onto a subspace spanned by few columns and few rows. The additional gain from using the compatibility constant is obtained because we use $\sqrt{2\mathrm{rank}\left(A\right)}/\left|\mathcal{P}_{A}\left(\Delta\right)\right|_*$ instead of 
$1/\left|\mathcal{P}_{A}\left(\Delta\right)\right|_2$ and hence avoid a type of Cauchy-Schwartz inequality $\left|\mathcal{P}_{A}\left(\Delta\right)\right|_*\le\sqrt{2\mathrm{rank}\left(A\right)}\left|\mathcal{P}_{A}\left(\Delta\right)\right|_2$.

\subsection{Restriction on the joint distribution of $X$ and $E$}\label{sec:XE}
%We maintain the following baseline assumption on the DGP.  
\begin{assumption}\label{ER} The following hold:
\begin{enumerate}[\textup{(}i\textup{)}]  
\item\label{ERi} There exists $\sigma>0$ such that $\left|E\right|_2^2/(NT)\xrightarrow{\mathbb{P}}\sigma^2$,
\item\label{ERii} There exists $\Sigma\in\mathcal{M}_{KK}$ positive definite such that $X^{\top}X/(NT)\xrightarrow{\mathbb{P}}\Sigma$,
\item\label{ERiii}   $X^{\top}e=O_{P}\left(\sqrt{NT}\right)$,
\item\label{ERiv}  There exists $\left\{\mu_{N}\right\}$ such that $\sum_{k=1}^K\left|X_k\right|_{\op}^2=O_P\left(\mu_N^2\right)$.
   \end{enumerate}
\end{assumption}
Condition \eqref{ERiii} is satisfied if, for all $k$, $\left\langle X_{k},E\right\rangle=\sum_{t=1}^T\sum_{i=1}^NX_{kit}E_{it}=O_P\left(\sqrt{NT}\right)$. This can allow for so-called predetermined regressors. This can be satisfied if, for some family of filtrations $\left(\mathcal{F}_{Nt}\right)_{t=1,\dots,T}$, for all $t=1,\dots,T$ and $i=1,\dots,N$, $X_{kit}$ is $\mathcal{F}_{Nt-1}$-measurable and $E_{it}$ is $\mathcal{F}_{Nt}$ measurable and, for example, under cross sectional independence. The role of \eqref{ERiv} is to introduce the notation $\{\mu_N\}$, this is not a restriction. Due to \eqref{ERii}, $\mu_N=O\left(\sqrt{NT}\right)$. $\{\mu_N\}$ sometimes appears in upper bounds in the results below and 
 slowly diverging sequences provide sharper results than when $\mu_N$ is of the order as large as $\sqrt{NT}$. 
 $\{\mu_N\}$ can diverge as slowly as $\sqrt{\max\left(N,T\right)}$ if the regressors satisfy the same assumptions as $E$ in Proposition \ref{proposal} or more generally those that can be found in \cite{onatski2015asymptotic,vershynin2010introduction} (see also Appendix A.1 in \cite{moon2015linear}). This will usually not hold if the regressors have a nonzero mean and more generally under the setup of Section \ref{s:transform}. In these cases, there exists $C>0$ such that, for all $N\in\N$, $\mu_N\ge C \sqrt{NT}$. Section \ref{s:transform} presents how to work with transformed regressors to obtain sharper results and complements the solution presented in the paragraph after Proposition \ref{proposal}.
 \begin{proposition}\label{sigmacons}
Under Assumption \ref{ER} with $\mu_N=\sqrt{NT}$, we have 
\begin{align}
&\left|\frac{\left|M_X(E)\right|_2}{\sqrt{NT}}-\sigma\right|=O_P\left(\frac{1}{\sqrt{NT}}\right)\label{ub1}\\
&\left|\norm{M_X(E)}_{\op}-\norm{E}_{\op}\right|=O_P\left(\frac{\mu_N}{\sqrt{NT}}\right)\label{ub2}.
\end{align}
\end{proposition}
The next assumption is % on the choice of $\left\{\lambda_{N}\right\}$. It is 
a sufficient condition for the event $\mathcal{E}$ in Theorem \ref{orac} to have a probability which converges to 1. 
\begin{assumption}\label{Pen}  
Maintain Assumption \ref{ER} and, if an upper bound $\mu_N$ for Assumption \ref{ER} \eqref{ERiv} is available else take $\mu_N=\sqrt{NT}$, take $\left\{\lambda_{N}\right\}$ of the form 
\begin{equation}\label{choiceL}
\lambda_N=\left(1-\frac{\phi_{1N}}{\sqrt{NT}}\right)^{-1}\left(\psi_{N}+\phi_{2N}\frac{\mu_N}{\sqrt{NT}}\right),
\end{equation}
where $\left\{\phi_{1N}\right\}$ and $\left\{\phi_{2N}\right\}$ are arbitrary sequences going to infinity, as slowly as we want but no faster than $\sqrt{NT}$ for $\left\{\phi_{1N}\right\}$, and \begin{enumerate}[\textup{(}i\textup{)}]  
  \item\label{Peni} $\psi_N=O\left(\sqrt{NT}\right)$, 
  \item\label{Penii} $\lim_{N\to \infty}\mathbb{P}\left(\rho\psi_N\sigma\ge \norm{E}_{\op}
\right)=1$.
\end{enumerate}
\end{assumption}
We can take $\phi_1=\phi_2$ in which case we write $\phi=\phi_1=\phi_2$. %\eqref{choiceL} holds whether $\mu_N=\sqrt{NT}$ or we have a sharper bound on it. 
Under the premises of Section \ref{s:transform}, we can take $\mu_N=\lambda_N$ and
\begin{equation}\label{choiceLb}
\lambda_N=\left(1-\frac{\phi_N}{\sqrt{NT}}\right)^{-1}\psi_N.
\end{equation}
We have
$$\mathcal{E}=\left\{\rho\psi_N\sigma+\rho\frac{\phi_{2N}\mu_N
}{\sqrt{NT}}\sigma+\rho\frac{\phi_{1N}\lambda_N
}{\sqrt{NT}}\sigma\ge \norm{E}_{\op}+\left(\norm{M_X(E)}_{\op}
-\norm{E}_{\op}\right)+\rho\lambda_N\left(\sigma-\frac{\left|M_X(E)\right|_2}{\sqrt{NT}}\right)\right\},$$
hence
\begin{align*}
\mathbb{P}\left(\mathcal{E}\right)&\ge 
\mathbb{P}\left(\left\{\rho\psi_N\sigma\ge \norm{E}_{\op}\right\}\bigcap
\left\{\rho\frac{\phi_{2N}\mu_N
}{\sqrt{NT}}\sigma\ge \norm{M_X(E)}_{\op}
-\norm{E}_{\op}\right\}\bigcap\left\{\frac{\phi_{1N}
}{\sqrt{NT}}\sigma\ge\sigma-\frac{\left|M_X(E)\right|_2}{\sqrt{NT}}
\right\}\right)
%&=\mathbb{P}\left(\psi_N\sigma+\phi_N\sigma\ge 2\norm{E}_{\op}+2\left(\norm{M_X(E)}_{\op}
%-\norm{E}_{\op}\right)+\lambda_N\left(\sigma-\frac{\left|M_X(E)\right|_2}{\sqrt{NT}}\right)
%\right),\\
%&\le\mathbb{P}\left(\psi_N\sigma\ge 2\norm{E}_{\op}\right)
%+\mathbb{P}\left(\phi_N\sigma>2\left(\norm{M_X(E)}_{\op}
%-\norm{E}_{\op}\right)+\lambda_N\left(\sigma-\frac{\left|M_X(E)\right|_2}{\sqrt{NT}}\right)
%\right),
\end{align*}
and the 3 events have probability going to 1 by \eqref{Penii} and Proposition \ref{sigmacons} so
%hence Assumption \ref{Pen} (ii) guarantees that the first probability in the right-hand side goes to 1, 
%(i) allows to treat the second term in the second probability, and $\left\{\phi_N\right\}$ is introduced to deal with \eqref{ub1} and \eqref{ub2} which in the worst case scenario give rise to terms which are $O_P(1)$. With such a choice of $\left\{\lambda_N\right\}$, the event appearing in Theorem \ref{orac} satisfies
$\lim_{N\to \infty}\mathbb{P}\left(\mathcal{E}\right)=1$.
% If $\mu_N=o\left(\sqrt{NT}\right)$, we can replace (i) by $\lambda_N=o\left(\sqrt{NT}\right)$ and take a constant $\left\{\phi_N\right\}$.

%\begin{proposition}\label{choicenuc}
%Under Assumption \ref{ER}, we have 
%$$2\sqrt{NT}\frac{\norm{M_X(E)}_{\op}}{\norm{M_X(E)}_2}=2\frac{\norm{E}_{\op}}{\sigma}+o_P(\norm{E}_{\op})+O_P(1).$$
%If $\left\{\psi_N\right\}$ is such that $\lim \limits_{N\to \infty}\mathbb{P}\left(\psi_N\ge 2\frac{\norm{E}_{\op}}{\sigma}\right)=1$ and $\varphi_N\to \infty$, then for any $c>1$, $\lambda_N=c\psi_N+\varphi_N$ satisfies Assumption \ref{Pen}.
%\end{proposition}
%The choice of the penalty given in Assumption \ref{Pen} is independent of the distribution of the regressors. 
We can handle large classes of joint distributions of $X$ and $E$, including ones where the errors have heavy tails. 
It is usual, but not necessary, to work with classes of distributions such that $\norm{E}_{\op}=O_P\left(\sqrt{\max(N,T)}\right)$.  For such distributions, it is enough to take $\psi_N=C\sqrt{\max(N,T)}$ for large enough $C$ for Assumption \ref{Pen} to hold. %This does not tell which constant to pick. 
An easy way to circumvent the problem that $C$ is unknown is to take $\psi_N=\phi_{2N}\sqrt{\max(N,T)}$ but this results in over penalization. 
%Regarding the choice of the penalty level, we maintain the following assumption:
% \begin{assumption}[Penalty] \label{Pen} 
%\end{assumption}
%At this level, the choice of the penalty level is not very practical but it allows to 
%It plays the same role as condition (15) in \cite{moon2018nuclear}. 
At the cost of additional assumptions on the distribution, one can obtain the following more precise proposal based on Corollary 5.35 and Theorem 5.31 in \cite{vershynin2010introduction}. 
\begin{proposition} \label{proposal} If $E=M_u\eta M_v$, where $M_u$ and $M_v$ are, possibly random, matrices such that  $\left|M_u\right|_{\op}\le1$ and $\left|M_v\right|_{\op}\le 1$ and 
either of the following holds
\begin{enumerate}[\textup{(}i\textup{)}]  
\item\label{proposali} $\{\eta_{it}\}_{i,t}$ are i.i.d. centered Gaussian random variables,
\item\label{proposalii} $\{\eta_{it}\}_{i,t}$ are i.i.d. centered random variables with finite fourth moments and $T/N$ converges to a constant in $[0,1]$, 
\end{enumerate}
then the sequence defined by $\psi_N=\left(\sqrt{N}+\sqrt{T}\right)/\rho+\varphi_N$, where $\varphi_N\to \infty$ arbitrarily slowly in case \eqref{proposali} and 
$\left\{\varphi_N/\sqrt{T}\right\}$ is bounded away from 0 %\to \infty$ 
in case \eqref{proposalii},  satisfies Assumption \ref{Pen} \eqref{Penii}.
\end{proposition}
The matrices $M_u$ and $M_v$ can be known or estimated (see, \emph{e.g.}, Section \ref{s:transform}) and have been applied to the data. Applying such matrices cannot increase $\rank\left(M_u\Gamma^l M_v\right)$, $\left|M_u \Gamma^d M_v\right|_{\op}$, or $\left|M_u \eta M_v\right|_\op$ but can reduce the operator norm of the regressors and give rise to  a smaller sequence $\left\{\mu_N\right\}$. These matrices can be unknown and the baseline error $E$ can have temporal and cross-sectional dependence. Because the operator norm of a matrix is equal to the operator norm of its transpose, the role of $N$ and $T$ can be exchanged in \eqref{proposalii}.  The proposed choice of the penalty level is almost completely explicit and does not depend on the variance of the errors. The remaining sequences are arbitrary. In contrast to \eqref{choiceL} where $\left(1-\phi_{1N}/\sqrt{NT}\right)^{-1}$ converges to 1, %though for the case Lasso,  
\cite{moon2018nuclear} employs a factor converging to infinity. Hence, %the method of this paper provides a data-driven thresholding rule in Step \eqref{s2} in Section \ref{sec:ca} of order smaller than that implicitly used in \cite{moon2018nuclear}. In other words, 
the data-driven method of this paper provides less shrinkage, less bias, and a better bias/variance tradeoff. 

\subsection{Restriction on the joint distribution of $X$ and $\Gamma$}
%We now discuss restrictions so that the bounds in Theorem \ref{orac} are small. 
\begin{assumption}
\label{ALR} $\rho$ and $\widetilde{\rho}$ are given and the random matrix $\Gamma$ can be decomposed as $\Gamma=\Gamma^l+\Gamma^d$, where, for $\left\{r_N\right\}$, 
\begin{enumerate}[\textup{(}i\textup{)}]  
 \item\label{ALRi} $\mathrm{rank}\left(\Gamma^l\right)=O_P(r_N)$,
 \item\label{ALRii} $\left|\Gamma^d\right|_*=O_P\left(\lambda_N r_{N}\right)$,
 \item\label{ALRiii} There exists $ \kappa>0$ independent of $N$ such that $\kappa_{\Gamma^l,c\left(\rho,\widetilde{\rho}\right)}\ge\kappa$ w.p.a.  1.
 \end{enumerate}
 \end{assumption}
We maintain  Assumption \ref{ALR} to translate the result of Theorem \ref{orac} into rates of convergence.  \eqref{ALRi} allows for ranks which are random and can vary with the sample size which is more general and realistic than the usual assumption that the rank is fixed. Condition \eqref{ALRiii} is a condition on the second random element of the first term in the right-hand side of \eqref{ubtheta}. Condition \eqref{ALRiii} is introduced because $\kappa_{\Gamma^l,c\left(\rho,\widetilde{\rho}\right)}$ is random. Such an assumption would not be required if $X$ and $\Gamma$ were fixed.  Unlike other papers on the topic, this paper allows for $\Gamma^d\ne0$. By the oracle type inequalities of Theorem \ref{orac}, the estimator performs as well as the best infeasible trade-off. Assumptions \eqref{ALRi} and \eqref{ALRii} are not restrictive because $r_N$ can be anything. The idea though to obtain tight results is to have $r_N$ small and realize a trade-off. \eqref{ALRi} is the reason why the component $\Gamma^l$ is called low-rank. The component $\Gamma^d$ can be viewed as a remainder which can have an arbitrary rank. Before the statement of Assumption \ref{ALR}, $\Gamma^l$ and $\Gamma^d$ are not precisely defined. Their sum is $\Gamma$ so parts of  $\Gamma^l$ can be transferred to  $\Gamma^d$ and vice versa. Assumption \ref{ALR} makes it more precise which component is which.
\begin{proposition}\label{propLB2}
Assumption \ref{ALR} \eqref{ALRiii} for a cone with constant $c$ holds with the lower bound $\kappa$ if, w.p.a. 1, 
$\kappa^2+2\rank\left(\Gamma^l\right)Q(b,b_\perp)
\le 1$, where $b,b_\perp\in\mathbb{R}^K$ are defined, for $k=1,\dots,K$, as $b_k=a\min\left(\left|\mathcal{P}_{\Gamma^l}\left(X_k\right)\right|_{\op},\left|X_k\right|_{\op}\right)$, $b_{\perp k}=a\left|\mathcal{P}_{\Gamma^l}^{\perp}\left(X_k\right)\right|_{\op}$, $a=\left|X^{\top}X/(NT)\right|_{\op}^{-1}
|X|_2/(NT)$, 
\begin{align*}
Q(b,b_\perp)=&
|b|_2^2\indic\left\{p_N\left|b_\perp\right|_2^2\ge1\right\}+
\left(\left|b+b_{\perp}c\right|_2^2-\frac{c^2}{p_N}\right)\indic\left\{
1-\frac{p_N\langle b_{\perp},b\rangle}{c}\le
p_N\left|b_\perp\right|_2^2<1\right\}\\
&+\left(\left|b+b_{\perp }
\frac{p_N\langle b_{\perp},b\rangle}{1-p_N\left|b_{\perp}\right|_2^2}\right|_2^2-\frac{p_N\langle b_{\perp},b\rangle^2}{\left(1-p_N\left|b_{\perp}\right|_2^2\right)^2} \right)\indic\left\{
p_N\left|b_{\perp}\right|_2^2<1-\frac{p_N\langle b_{\perp},b\rangle}{c}\right\},
\end{align*}
and $p_N=\min\left(N-\rank\left(\Gamma^l\right),T-\rank\left(\Gamma^l\right)\right)$.
\end{proposition}
Note that $Q(b,b_\perp)<\left|b+b_{\perp}c\right|_2^2$ and, if $K=1$, 
$a=1/\left|X_1\right|_2$ and $$\left|b+b_{\perp}c\right|_2^2=\frac{1}{\left|X_1\right|_2^2}
\left(\min\left(\left|\mathcal{P}_{\Gamma^l}\left(X_1\right)\right|_{\op},\left|X_1\right|_{\op}\right)+\left|\mathcal{P}_{\Gamma^l}^{\perp}\left(X_1\right)\right|_{\op} c\right)^2.$$ 
The quantity $\left|\mathcal{P}_{\Gamma^l}^{\perp}\left(X_k\right)\right|_{\op}= \left|M_{u\left(\Gamma^l\right)}X_kM_{v\left(\Gamma^l\right)}\right|_{\op}$ in the definition of $b_{\perp k}$ can be not too large because the projectors can reduce the operator norm if $X_k$ has a component with a factor structure and shares some factors in common with $\Gamma^l$ which are annihilated by $M_{v\left(\Gamma^l\right)}$ (see Remark \ref{ra} for further discussion of this aspect). 
Due to Assumption \ref{ER} \eqref{ERii}, $a=O_P\left(1/\sqrt{NT}\right)$.  In the worst case, by the crude bound $\left|X_k\right|_{\op}\le\left|X_k\right|_2$, $b$ and $b_\perp$, hence $Q(b,b_\perp)$ are bounded. If $\mu_{N}=o\left(\sqrt{NT}\right)$, the condition in Proposition \label{propLB} holds for arbitrary constants $\kappa<1$ for $N$ large enough, but this is not necessary. Section \ref{s:transform} presents solutions to work with regressors with smaller operator norm.  Lemma A.7 in \cite{moon2018nuclear} provides an alternative sufficient condition for Assumption \ref{ALR} \eqref{ALRii}. Lemma A.8 is another sufficient condition when $K=1$. In our framework $r_{1N}$ can grow, $c$ can be different from $3$, and we do not work contionnal on $\Gamma^l$,  condition (iii) has to be modified with a denominator of $\sqrt{NTr_{N}}$ and the probabilities are with respect to the distribution of $(\Gamma,X_1)$. It is claimed in Remark (a) in \cite{moon2018nuclear} that the condition in Lemma A.8 holds when $X_1=\Pi_1^l+U_1$, $\Pi_1^l$ has a fixed rank, and $U_1$ has iid mean zero normal entries. 

\subsection{Rates of convergence}\label{s:estbeta}
Theorem \ref{orac} and the assumptions on the DGP yield the following. % theorem.
\begin{theorem} \label{const1}
Under assumptions \ref{Pen} and \ref{ALR}, 
\begin{align}
&\left|\widehat{\Gamma}-\Gamma\right|_*=O_P\left(\lambda _Nr_{N} \right),\label{rategamma}\\
&\widehat{\sigma}-\sigma
=O_P\left(\frac{\lambda _N^2r_{N}}{NT} \right),\label{Oraclevar} \\
&\widehat{\beta}-\beta=O_P\left(\frac{\lambda_N r_{N}\mu_N}{NT}\right).\label{ratebeta}
\end{align}
\end{theorem}
In \eqref{ratebeta}, we have implicitly assumed that $\sqrt{NT}=O\left(\lambda_N r_{N}\mu_N\right)$ but this always occurs when $X\ne0$ and the problem is to have $\lambda_N r_{N}\mu_N$ as close as possible in rate to $\sqrt{NT}$. Under usual assumptions where we can take $\lambda_N$ proportional to $\sqrt{\max(N,T)}$, $r_{N}$ fixed, and make no restriction on $\left\{\mu_N\right\}$ so that  $\mu_N=O(\sqrt{NT})$, we obtain the rate  convergence of $1/\sqrt{\min(N,T)}$ which is the one in \cite{moon2018nuclear}. Theorem \ref{const1} shows  that $\widehat{\beta}$ remains consistent if $r_N=o\left(\sqrt{\min(N,T)}\right)$. Obviously $r_N$ can be larger if $\mu_N$ is smaller. The most favorable situation, when $\mu_N=O\left(\sqrt{\max\left(N,T\right)}\right)$ and $\lambda_N$ is proportional to $\sqrt{\max(N,T)}$, yields
$\widehat{\beta}-\beta=O_P\left(\max(N,T)  r_{N}/(NT)\right)$, 
hence, when $N/T$ has a positive limit, this becomes $O_P\left(r_{N}/\sqrt{NT}\right)$. Achieving  $\mu_N=o\left(\sqrt{NT}\right)$ and in some cases $\mu_N=O\left(\sqrt{\max\left(N,T\right)}\right)$ using transformed regressors is sometimes possible under the premises of Section \ref{s:transform} and this paper allows to obtain such an estimator and transformed regressors in a data-driven way.  Section \ref{sec:modif} proposes an alternative approach where we can obtain the $1/\sqrt{NT}$ rate  and have asymptotic normality.

\subsection{Additional results using the relation to the matrix Lasso}\label{sec:add}
Because any solution $\left(\widehat{\beta},\widehat{\Gamma}\right)$ of \eqref{estimator2} is solution of 
\eqref{estimator3}, we prove the following additional results. They would also apply to \eqref{estimator3} even if, rather than $\widehat{\sigma}$, we used an upper bound on the standard error of the errors. The results that we state involve $\widehat{\sigma}$ but, under the assumptions of Theorem \ref{const1}, $\widehat{\sigma}$ is a consistent estimator of $\sigma$. In order to guarantee $\mathbb{P}\left(\rho\lambda_N\min\left(\widehat{\sigma},\sigma\right)\ge\norm{M_X(E)}_{\op}
\right)\to 1$ we need the following assumption.
\begin{assumption}\label{Pen2} 
Assumption \ref{Pen} holds and $\left\{\phi_{1N}\right\}$ satisfies the additional restriction that, for $N$ large enough,  
$$\left(1-\frac{\phi_{1N}}{\sqrt{NT}}\right)^2\phi_{1N}\ge\phi_{2N}\frac{r_{N}}{\sqrt{NT}}\left(\psi_N+\phi_{2N}\frac{\mu_N}{\sqrt{NT}}\right)^2.$$
\end{assumption}
Indeed, we can replace $\left\{\phi_{1N}\sigma/\sqrt{NT}\ge\sigma-\left|M_X(E)\right|_2/\sqrt{NT}
\right\}$ by $\left\{\phi_{2N}\sigma\lambda_N^2r_N/NT\ge\sigma-\widehat{\sigma}
\right\}$ in the previous analysis which converges to 1 due to \eqref{Oraclevar} because, due to Assumption \ref{Pen2}, 
$\phi_{1N}\ge\phi_{2N}\lambda _N^2r_{N}/\sqrt{NT} $, hence 
$$
\left(1-\phi_{2N}\frac{\lambda _N^2r_{N}}{NT} \right)\lambda_N\ge\left(1-\frac{\phi_{1N}}{\sqrt{NT}} \right)\lambda_N=\psi_N+\phi_{2N}\frac{\mu_N}{\sqrt{NT}}.
$$
A  conservative choice is $\phi_{1N}=c_1\sqrt{NT}$ for a small $c_1\in(0,1)$. Now on, we use cones with constant $c=c\left(\rho\right)=(1+\rho)/(1-\rho)$. First, with a proof similar to the computations in \cite{koltchinskii2011nuclear}, we obtain a result which is an oracle inequality with constant 1 if $X$ and $\Gamma$ are not random.

\begin{proposition}\label{l2}
If $\rho\lambda\min\left(\widehat{\sigma},\sigma\right)\ge \norm{M_X(E)}_{\op}$, we have
%For $\delta\in(0,1)$, on the event $\mathcal{E}_{\lambda\delta}$, we have
 \begin{align*}
\frac{1}{NT}\left|M_{ X}\left(\Gamma-\widehat{\Gamma}\right)\right|_2^2&\le \inf_{\widetilde{\Gamma}}\left\{
\frac{1}{NT} \left|M_{ X}\left(\Gamma-\widetilde{\Gamma}\right)\right|_2^2+
\frac{2(\lambda(1+\rho)\min\left(\widehat{\sigma},\sigma\right))^2}{NT}\frac{\rank\left(\widetilde{\Gamma}\right)}{\kappa_{\widetilde{\Gamma},c(\rho)}^2}\right\}.
%\\
%\frac{1}{NT}\left|\Pi-\widehat{\Pi}\right|_2^2&\le \inf_{\widetilde{\Pi}}\left\{
%\frac{1}{NT} \left|\Pi-\widetilde{\Pi}\right|_2^2+\left(1+\sqrt{2}\delta\right)^2\lambda\widehat{\sigma}^2\rank(\Pi)\right\}.
\end{align*}
\end{proposition}
This inequality yields a slightly different notion of approximately sparse solution because the first term in the maximum involves $\left|M_{ X}\left(\Gamma-\widetilde{\Gamma}\right)\right|_2^2/(NT)$ rather than $\left|\Gamma-\widetilde{\Gamma}\right|_*$. The next result provides a bound on $\rank\left(\widehat{\Gamma}\right)$ as a function of the previous bound. 
\begin{proposition}\label{l1}
If $\rho\lambda\widehat{\sigma}\ge \norm{M_X(E)}_{\op}$ then 
%and, for $\eta\in[0,1)$, $\Gamma=\Gamma^l+\Gamma^d$, where 
%$\left|M_{ X}\left(\Gamma^d+ E\right)\right|_{\op}\le 2(1-\eta)\left|M_{ X}\left( E\right)\right|_{\op}$,
we have
$$
 \left(\lambda(1-\rho)\widehat{\sigma}-\left|\Gamma^d\right|_{\op}\right)_+^2\rank\left(\widehat{\Gamma}\right)
\le 
\left|P_{u\left(\widehat{\Gamma}\right)}M_{ X}\left(\Gamma^l-\widehat{\Gamma}\right)P_{v\left(\widehat{\Gamma}\right)}\right|_2^2\le\left|M_{ X}\left(\Gamma^l-\widehat{\Gamma}\right)\right|_2^2.$$
\end{proposition}
% $$\left|P_{u\left(\widehat{\Gamma}\right)}M_{ X}\left(\Gamma^l-\widehat{\Gamma}\right)P_{v\left(\widehat{\Gamma}\right)}\right|_2\ge \left(\lambda(1-\rho)\widehat{\sigma}-\left|\Gamma^d\right|_{\op}\right)\sqrt{\rank\left(\widehat{\Gamma}\right)}.$$
As a result, under the above conditions and Assumtion \ref{ALR} (ii), % if $\Gamma=\Gamma^l$ and $\eta=1/2$,  
$$\rank\left(\widehat{\Gamma}\right)\le 2\left((1+\rho)/((1-\rho)\kappa_{\Gamma^l,c(\rho)})\right)^2\rank\left(\Gamma^l\right).$$ 
%But the result is more general because it allows for the presence of a nonzero $\Gamma^{d}$ such that $\left|M_{ X}\left(\Gamma^{d}\right)\right|_2^2/NT$ is not too large. 
We can combine propositions \ref{l2} and \ref{l1} with Proposition \ref{c2} in the appendix to obtain results for other loss functions, in particular the Frobenius norm.

Our estimator has desirable low-rank properties but it can fail to obtain $\rank(\Gamma)$, $\rank\left(\Gamma^l\right)$, or annihilator matrices. Thus, we introduce the hard-thresholded estimator 
$$\widehat{\Gamma}^{t}=\sum_{k=1}^{\rank\left(\widehat{\Gamma}\right)}
\sigma_k\left(\widehat{\Gamma}\right)\indic\left\{\sigma_k\left(\widehat{\Gamma}\right)\ge t\right\}
u_k\left(\widehat{\Gamma}\right)v_k\left(\widehat{\Gamma}\right)^{\top}.$$
\begin{proposition}\label{p:op}
Under the assumptions of Theorem \ref{const1} and Assumption \ref{Pen2}, if %$\lambda _N^3r_{N}=O(NT)$ and 
$\left|\Gamma^d\right|_{\op}=o_P\left(\lambda_N\sigma\right)$, we have
\begin{equation}\label{emax}
\max\left(\left|\Gamma-\widehat{\Gamma}\right|_{\op},\left|\Gamma^l-\widehat{\Gamma}\right|_{\op}\right)\le(\rho+1)\lambda_N\left(\sigma+O_P\left(\frac{ r_{N}\mu_N^2}{NT}\right)\right).
\end{equation}
\end{proposition}
\begin{assumption}\label{ass:SF} Let $h>1$. The following conditions hold
\begin{enumerate}[\textup{(}i\textup{)}]  
\item\label{ass:SFi} $r_{N}\mu_N^2=o(NT)$,
\item\label{ass:SFii} $\mathbb{P}\left(\sigma_{\rank\left(\Gamma^l\right)}\left(\Gamma^l\right)\ge (\rho+1)\lambda_Nh^2(h+1)\sigma\right)\rightarrow 1$.
\end{enumerate}
\end{assumption}
Condition \eqref{ass:SFi} %in Assumption \ref{ass:SF} 
guarantees %, under the assumptions of Proposition \ref{p:op},
 the $O_P$ in \eqref{emax} is $o_P(1)$. It allows the pivotal thresholding methods below but imposes a slight restriction on the operator norms of the regressors. Section \ref{s:transform} allows to come back to a case where \eqref{ass:SFi} holds for a large class of regressors. Without \eqref{ass:SFi} 
$$\max\left(\left|\Gamma-\widehat{\Gamma}\right|_{\op},\left|\Gamma^l-\widehat{\Gamma}\right|_{\op}\right)=O_P\left(\lambda_N\right)$$
and can adapt the results which follow at the expense of a theoretical but unfeasible thresholding level or using conservative levels $\lambda_N/t=o(1)$. Condition \eqref{ass:SFii} %in Assumption \ref{ass:SF} 
is weaker than a strong-factor assumption on $\Gamma^l$. %, in the benchmark case where $\lambda_N$ is proportional to $\sqrt{\max(N,T)}$. %\ll\sqrt{NT}$. 
%Based on Proposition \ref{p:op}, 
We now show that we can recover $\rank(\Gamma)$ with a data-driven threshold.
\begin{proposition}\label{p:th}
Under the assumptions of Proposition \ref{p:op} and Assumption \ref{ass:SF}, then setting $t=(\rho+1)\lambda_Nh^2\widehat{\sigma}$ yields
$$\mathbb{P}\left(\rank\left(\widehat{\Gamma}^t\right)=\rank\left(\Gamma^l\right)\right)\rightarrow1.$$
Moreover, if we remove \eqref{ass:SFii}, then we have 
$$\mathbb{P}\left(\rank\left(\widehat{\Gamma}^t\right)\le\rank\left(\Gamma^l\right)\right)\rightarrow1,$$
if we replace \eqref{ass:SFii} by the weaker assumption $\mathbb{P}\left(\sigma_{\rank\left(\Gamma^l\right)}\left(\Gamma^l\right)\ge (\rho+1)\lambda_Nh^3\sigma\right)\rightarrow 1$, we have
$$\mathbb{P}\left(\rank\left(\widehat{\Gamma}^t\right)\ge\rank\left(\Gamma^l\right)\right)\rightarrow1,$$
and 
\begin{equation}\label{rateopt}
\max\left(\left|\Gamma-\widehat{\Gamma}^t\right|_{\op},\left|\Gamma^l-\widehat{\Gamma}^t\right|_{\op}\right)
\le(\rho+1)\lambda_N(h^2+1)\left(\sigma+o_P\left(1\right)\right).
\end{equation}
\end{proposition}
%\blue{La section 5 de Moon et Weidner propose une solution differente. Peut etre que celle-ci est plus simple? Mieux/moins bien.} \red{Je crois qu'il y a une erreur dans le Lemme 3 de Moon et Weidner, on ne peut pas avoir à la fois $\psi_{NT}\to 0$ and $\psi_{NT}/\sqrt{\min(N,T)}\to \infty$}
%Note as well that under the assumptions of Proposition \ref{p:op}, Assumption \ref{ass:SF} (i), and setting $t=3\lambda_Nc^2\widehat{\sigma}/2$, we have the following result for the thresholded estimator
We strengthen Assumption \ref{ass:SF} \eqref{ass:SFii} as follows. When $v_N$ increases like $\sqrt{NT}$, it is a strong-factor assumption. 
\begin{assumption}\label{ass:SF2} 
Let $\left\{v_N\right\}$ be such that $v_N\ge(\rho+1)\lambda_Nh^2(h+1)\sigma$. Assume that 
$$\mathbb{P}\left(\sigma_{\rank\left(\Gamma^l\right)}\left(\Gamma^l\right)\ge v_N \right)\rightarrow 1.$$
\end{assumption}
\begin{proposition}\label{p:proj} Under the assumptions of Proposition \ref{p:th} and Assumption \ref{ass:SF2}, we have
\begin{align*}
\left|P_{v\left(\widehat{\Gamma}^t\right)}-P_{v\left(\Gamma^l\right)}\right|_2&=\left|M_{v\left(\widehat{\Gamma}_v^t\right)}-M_{v\left(\Gamma^l\right)}\right|_2
\le(\rho+1)\frac{\sqrt{2r_N}\lambda_N}{v_N}\left((h^2+1)\sigma+o_P\left(1\right)\right)\\
\left|P_{u\left(\widehat{\Gamma}^t\right)}-P_{u\left(\Gamma^l\right)}\right|_2&=\left|M_{u\left(\widehat{\Gamma}^t\right)}-M_{u\left(\Gamma^l\right)}\right|_2\le(\rho+1)\frac{\sqrt{2r_N}\lambda_N}{v_N}\left((h^2+1)\sigma+o_P\left(1\right)\right).
\end{align*}
\end{proposition}
%\red{(En soit, nous n'avons pas besoin d'une strong-factor assumption pour obtenir ce résultat, non ?} 
Under a strong-factor assumption, when $\lambda_N$ is proportional to $\sqrt{\max(N,T)}$ %\ll\sqrt{NT}$, 
and $r_N$ is fixed, we obtain the same rate of convergence as using PCA and as in Lemma A.7 in \cite{bai2009panel}.  % if we observed $\Gamma^l+E$. 
Here we obtain an upper bound with known constant. The rates that we obtain are also more general because we do not need to maintain the strong-factor assumption or that $r_N$ is fixed, $\left\{\lambda_N\right\}$ could also allow for errors with larger tails of the operator norm. 
% \blue{(quelle reference mettrais tu?)} \red{On pourrait mettre notre papier sur le "toutPCA" justement, il y a un gros travail sur la PCA qui a été fait}  .

\subsection{Working with transformed regressors}\label{s:transform}
In the previous sections, $\left\{\mu_N\right\}$ sometimes plays an important role and we might want it to be not too large.  However, this can be as large as $O(\sqrt{NT})$ if the next assumption holds. %So we devote this section to the analysis of this difficult situation. 
%Assumption \ref{ER} (iv), Proposition \ref{propLB}, and our analysis of the convergence rate of $\widehat{\beta}$ in Section \ref{s:estbeta} depend on the operator norms of the regressors via $\left\{\mu_N\right\}$ and this should not to be too large. The second condition in Assumption \ref{ER} (iv) does not hold 
\begin{assumption}\label{assX}
For at least one $k\in\{1,\dots,K\}$, 
\begin{equation}\label{efactorX}
X_k=\Pi_k^l+\Pi_k^d+U_k,
\end{equation}
and $\Pi_k^d$, $U_k$, $\sigma_k$, $r_{kN}$, $\lambda_{kN}$, and $v_{kN}$  play the role of $\Gamma^d$, $E$, $\sigma$, $r_N$, $\lambda_N$ and $v_N$ and satisfy the assumptions of Proposition \ref{proposal}, Assumption \ref{ALR} \eqref{ALRi} and \eqref{ALRii}, and Assumption \ref{ass:SF} \eqref{ass:SFii}, %$\Pi_k^l\ne0$ %and $\left|\Pi_k^l\right|_{\op}=O_P\left(\sqrt{NT}\right)$
%and 
$\left|\Pi_k^l\right|_{\op}^{-1}=O_P\left(1/\sqrt{NT}\right)$.
\end{assumption} 
The problem is difficult due to $\left|\Pi_k^l\right|_{\op}^{-1}=O_P\left(1/\sqrt{NT}\right)$. This occurs under a strong-factor assumption when the ratio of any singular value of $\Pi_k^l$ and $\sqrt{NT}$ has a positive and finite limit in probability.  %deterministic limit as $N$ goes to infinity. 
%No transformation is required if $\Pi_k^l=0$ or if $\left|\Pi_k^l\right|_{\op}=o_P\left(\sqrt{NT}\right)$.  
The problem would be even harder if $\Pi_k^l$ does not have a small rank (\emph{i.e.}, with ``many'' strong factors) and there is obviously a problem related to identification when $X_k=\Pi_k^l$ and $\Pi_k^l$ has small rank. Under the aforementioned assumptions,  we can take $\lambda_{kN}=\lambda_N$. The matrix $\Pi_k^l$, $\sigma_k$, and the annihilators $M_{u(\Pi_k^l)}$ and $M_{v(\Pi_k^l)}$ can be estimated like in the previous sections and one can replace $X_k$ by $\widetilde{X}_k$, where $X_k-\widetilde{X}_k$ has low rank, and $\Gamma^l$ by $\widetilde{\Gamma}^l=\Gamma^l+\sum_{k=1}^K\beta_k\left(X_k-\widetilde{X}_k\right)$. For simplicity of exposition, we apply a transformation to all regressors. When $X=0$, \eqref{estimator2} can be computed as an iterated soft-thresholding estimator. % (see Section \ref{sec:ca}). 

One can work with an estimator $\widetilde{\Pi}_k$ of $\Pi_k$ of the form 
$\widetilde{\Pi}_k=\widehat{\Pi}_k$ or $\widetilde{\Pi}_k=\widehat{\Pi}_k^t$ obtained as described in the previous sections, with (1) $\widetilde{X}_k=X_k-\widetilde{\Pi}_k$, (2) $\widetilde{X}_k=M_{u\left(\widetilde{\Pi}_k\right)}X_k$, (3) $\widetilde{X}_k=X_kM_{v\left(\widetilde{\Pi}_k\right)}$, (4) $\widetilde{X}_k=\mathcal{P}_{\widetilde{\Pi}_k}^\perp\left(X_k\right)$, (5) $\widetilde{X}_k=X_k-X_k^{(l_k)}$ where $X_k^{(l_k)}$ is obtained from $X_k$ by keeping the low rank component from a SVD corresponding to the $l_k=\rank\left(\widetilde{\Pi}_k\right)$ largest singular values. An alternative %for which we do not give details in this paper 
is to rely on Principal Component Analysis (henceforth PCA) using the eigenvalue-ratio (see \cite{ahn2013eigenvalue}) to select the number of factors. 
By the previous results, using such transformed regressors gives rise to additional terms in $\widetilde{\Gamma}$ of rank each at most $18r_{kN}+o_P(1)$ if  $\Pi_k=\Pi_k^l$ or of same rank as $\widetilde{X}_k^l$ w.p.a. 1 if we use hard-thresholding as well. % and allow for $\Pi_k^d\ne0$. 
Assuming we transform all regressors, the rank of $\widetilde{\Gamma}$ is at most $\widetilde{r}_N+o_P(1)$, where $\widetilde{r}_N=r_N+2((1+\rho)/(1-\rho))^2\sum_{k=1}^Kr_{kN}$ if $\widetilde{\Pi}_k=\widehat{\Pi}_k$
and $l_k=\rank\left(\widehat{\Pi}_k\right)$ and else $\widetilde{r}_N=r_N+\sum_{k=1}^Kr_{kN}$.
Using $\widetilde{\Pi}_k=\widehat{\Pi}_k^t$ has the advantage that if $\Pi_k^d\ne0$ we have guarantees on the low rank of $\widetilde{\Gamma}$. %Else we might have to also modify $\Gamma^d$.

\begin{remark}
In Assumption \ref{assX} we have assumed that we maintain the assumption of Proposition \ref{proposal} and Assumption \ref{ass:SF} \eqref{ass:SFii} for simplicity of exposition. But we can also handle heavy tailed errors $U_k$ by choosing a penalty level $\lambda_{kN}$ large enough as disscussed before Proposition \ref{proposal}  . 
We maintain Assumption \ref{ass:SF} \eqref{ass:SFii} to allow for a simple thresholding rule but it is enough to use a thresholding at any level of smaller order than $\sqrt{NT}$ to obtain $\mu_N=o\left(\sqrt{N}\right)$.  
\end{remark}

\subsection{Second-stage estimator of $\beta$}\label{sec:modif}
As seen at the end of Section \ref{s:estbeta}, the estimator $\widehat{\beta}$ could sometimes achieve the $1/\sqrt{NT}$ rate. But under weaker conditions we obtain a slower rate of convergence. This section presents three different two-stage approaches which deliver an asymptotically normal estimator of $\beta$.

\subsubsection{Approach 1: Annihilation of low-rank components of $\Gamma$ and the regressors}
This section analyzes another approach under Assumption \ref{assX} where, for simplicity of exposition, the last statement holds for all regressors, 
%\eqref{efactorX} and the conditions under it 
and we use the transformed regressors with transformation (1) or (2). We obtain estimators of
$\Pi_u^l=\left(\Gamma^l,\Pi_1^l,\dots,\Pi_K^l\right)$ and
$\Pi_v^l=\left(\left(\Gamma^l\right)^\top,\left(\Pi_1^l\right)^\top,\dots,\left(\Pi_K^l\right)^\top\right)^\top$ by plug-in using $\widetilde{\Pi}_k=\widehat{\Pi}_k$ or 
 $\widetilde{\Pi}_k=\widehat{\Pi}_k^t$ (preferably) for $k=1,\dots,K$ and
\begin{equation}\label{estGamma}
\widehat{\Gamma}=\widehat{\widetilde{\Gamma}}-\sum_{k=1}^K\widehat{\beta}_k\widetilde{\Pi}_k.
\end{equation} We denote by $\widehat{\Pi}_u$ and $\widehat{\Pi}_v$ the estimators, by
$\overline{\sigma}^2=\sigma^2+\sum_{k=1}^K\sigma_k^2$ and $\widehat{\overline{\sigma}}^2=\widehat{\sigma}^2+\sum_{k=1}^K\widehat{\sigma}_k^2$, by $\widetilde{\sigma}=\overline{\sigma}$ and $\widehat{\widetilde{\sigma}}=\widehat{\overline{\sigma}}$ if $\widetilde{\Pi}_k=\widehat{\Pi}_k$, and by  
 $\widetilde{\sigma}=(h^2+1)\overline{\sigma}$ and $\widehat{\widetilde{\sigma}}=(h^2+1)\widehat{\overline{\sigma}}$ if $\widetilde{\Pi}_k=\widehat{\Pi}_k^t$.
Because, for $K\in\N$ and $A_1,\dots,A_K$ with same number of rows, $\left|(A_1,\dots,A_K)\right|_{\op}^2\le\sum_{k=1}^K|A_k|_{\op}^2$, and
$$\widehat{\Gamma}-\Gamma^l=\widehat{\widetilde{\Gamma}}-\widetilde{\Gamma}^l+\sum_{k=1}^K\left(\beta_k-\widehat{\beta}_k\right)\left(\widetilde{\Pi}_k-\Pi_k\right)+\sum_{k=1}^K\left(\beta_k-\widehat{\beta}_k\right)\Pi_k,$$
we obtain the following corollary of Proposition \ref{p:op} and \eqref{rateopt}.
\begin{corollary}\label{p:opc}
Under the assumptions \ref{ER}, % and either of those of Proposition \ref{proposal}, 
\ref{ALR}, where in \eqref{ALRiii} we have $\widetilde{\Gamma}^l$ instead of $\Gamma^l$, \ref{Pen2},  \ref{assX}, 
%\blue{(une partie etant estimee ca demande quelques verifications)}, 
$\lambda _N^2\widetilde{r}_{N}=o(NT)$, %$\lambda_N\widetilde{r}_N=o\left(\sqrt{NT}\right)$, 
and $\left|\Gamma^d\right|_{\op}=o_P\left(\lambda_N\sigma\right)$, we have
\begin{align*}
&\left|\Gamma^l-\widehat{\Gamma}\right|_{\op}\le(\rho+1)\lambda_N\left(\sigma+o_P\left(1\right)\right)\\
&\max\left(
\left|
%\left(\left(\Gamma^l\right)^\top,\left(\Pi_1^l\right)^\top,\dots,\left(\Pi_K^l\right)^\top\right)
\Pi_u^l-\widehat{\Pi}_u^l
%\left(\widehat{\Gamma}^\top,\widehat{\Pi}_1^\top,\dots,\widehat{\Pi}_K\right)
\right|_{\op},
\left|\Pi_v^l-\widehat{\Pi}_v^l
%\left(\widehat{\Gamma}^\top,\widehat{\Pi}_1^\top,\dots,\widehat{\Pi}_K\right)
\right|_{\op}\right)
\le(\rho+1)\lambda_N\left(\widetilde{\sigma}+o_P\left(1\right)\right).
%\\
%&\left|\left(\Gamma^l,\Pi_1^l,\dots,\Pi_K^l\right)-\left(\widehat{\Gamma},\widehat{\Pi}_1,\dots,\widehat{\Pi}_K\right)\right|_{\op}\le(\rho+1)\lambda_N\left(\overline{\sigma}+o_P\left(1\right)\right),
\end{align*}
\end{corollary}
Based on this corollary, we can rely on hard-thresholding of these estimators that we denote by $\widehat{\Gamma}^t$, $\widehat{\Pi}_u^t$ and $\widehat{\Pi}_v^t$ and estimate the rank of $\Gamma^l$ and the annihilator matrices $M_{u\left(\Gamma^l\right)}$,
$M_{v\left(\Gamma^l\right)}$, $M_{u\left(\Pi_u^l\right)}$, and 
$M_{v\left(\Pi_v^l\right)}$ by $M_{u\left(\widehat{\Gamma}^t\right)}$,
$M_{v\left(\widehat{\Gamma}^t\right)}$, $M_{u\left(\widehat{\Pi}_u^t\right)}$, and 
$M_{v\left(\widehat{\Pi}_v^l\right)}$. Again, the first two annihilators are estimated at the same rate as in Lemma A.7 in \cite{bai2009panel} if $\Gamma^l$ satisfies a strong-factor assumption. 
Proposition \ref{p:th}  and Proposition \ref{p:proj} hold with the annihilator matrices of this section replacing $\sigma$ by $\widetilde{\sigma}$ and $\widehat{\sigma}$ by $\widehat{\widetilde{\sigma}}$ and Assumption \ref{ass:SF} \eqref{ass:SFii} by 
$$\mathbb{P}\left(
\min\left(\sigma_{\rank\left(\Pi_u^l\right)}\left(\Pi_u^l\right),\sigma_{\rank\left(\Pi_v^l\right)}\left(\Pi_v^l\right)\right)
%\sigma_{\rank\left(\Pi_k^l\right)}\left(\Pi_k^l\right)
\ge (\rho+1)\lambda_Nh^2(h+1)\widetilde{\sigma}\right)\rightarrow 1$$ 
and Assumption \ref{ass:SF2} by $\lambda _N^2\widetilde{r}_{N}=o(NT)$ maintained in Corollary \ref{p:opc} and the next assumption.
\begin{assumption}\label{ass:SF2b} 
Let $\left\{\overline{v}_{N}\right\}$ be such that $\overline{v}_{N}\ge  (\rho+1)\lambda_Nh^2(h+1)\widetilde{\sigma}$, we have 
$$\mathbb{P}\left(\min\left(\sigma_{\rank\left(\Pi_u^l\right)}\left(\Pi_u^l\right),\sigma_{\rank\left(\Pi_v^l\right)}\left(\Pi_v^l\right)\right)\ge \overline{v}_N \right)\rightarrow 1$$
and, for a sequence $\left\{\overline{r}_N\right\}$, 
%There exists a random $N\times T$ matrix $\Gamma^l$ \blue{(fixe ou aleatoire?)} \red{(aleatoire si on prend $\Gamma$ aléatoire)} such that
$\max\left(\mathrm{rank}\left(\Pi_u^l\right),\mathrm{rank}\left(\Pi_v^l\right)\right)=O_P(\overline{r}_N).$
\end{assumption}
We denote by $\mathcal{P}_{\widehat{\Pi}^t}^\perp$ (resp. $\mathcal{P}_{\Pi}^\perp$) the operator which applied to $A\in\mathcal{M}_{NT}$ is $\mathcal{P}_{\widehat{\Pi}}^\perp (A)=M_{u\left(\widehat{\Pi}_u^t\right)}AM_{v\left(\widehat{\Pi}_v^t\right)}$ (resp. $\mathcal{P}_{\Pi}^\perp (A)=M_{u\left(\Pi_u\right)}AM_{v\left(\Pi_v\right)}$) 
 and define the estimator \begin{equation}\label{ewan}
\widetilde{\beta}^{(1)}\in\argmin{\beta\in\mathbb{R}^K}\left|\mathcal{P}_{\widehat{\Pi}^t}^\perp\left(Y-\sum_{k=1}^K\beta_kX_k\right)\right|_2^2.
%\left|M_{u\left(\widehat{\Pi}_u^t\right)}\left(Y-\sum_{k=1}^K\beta_kX_k\right)M_{v\left(\widehat{\Pi}_v^t\right)}\right|_2^2.
\end{equation}
%The researcher can work with transformed regressors like in Section \ref{s:transform} to obtain the annihilator matrices but the matrices of regressors in \eqref{ewan} could be the original or other transformed ones. 
Also $\mathcal{P}_{\widehat{\Pi}^t}^\perp(X)$ (resp. $\mathcal{P}_{\widehat{\Pi}^t}^\perp(U)$, $\mathcal{P}_{\Pi}^\perp(X)$, and $\mathcal{P}_{\Pi}^\perp(U)$) is the matrix 
formed like $X$, replacing the matrices $X_k$ by $\mathcal{P}_{\widehat{\Pi}^t}^\perp(X_k)$ (resp. $\mathcal{P}_{\Pi}^\perp(X_k)$, $\mathcal{P}_{\Pi}^\perp(U_k)$, and $\mathcal{P}_{\Pi}^\perp(U_k)$) for $k=1,\dots,K$. %This estimator shares some similarity with the one in 
%Unlike \cite{bai2009panel}, %. Two annihilators rather than one on the left are used. Also, 
%the algorithm is not iterative and we can be certain that we have computed $\widetilde{\beta}^{(1)}$ for which the theory below applies rather than a local minimum. 
%The annihilators can be estimated as fast as in \cite{bai2009panel} under the premises of Section \ref{sec:add} \blue{(il faudrait etre plus precis)} even without knowing the ranks of the matrices. Using two of them allows to obtain the following result. 
%Note that the objective function is $\left|\mathcal{P}_{\widehat{\Gamma}^t}^\perp\left(Y-\sum_{k=1}^K\beta_kX_k\right)\right|_2^2$. 
%We use the notation, for $\Delta\in\mathcal{M}_{NT}$, $\dot{\Delta}=M_{u\left(\Pi_u^l\right)}\Delta M_{v\left(\Pi_v^l\right)}$ and construct the stacked matrix $\dot{X}$ similarly to $X$ and use dots for vectorized elements as well. Similarly we use the notations $\ddot{\Delta}$ for all matrices and vectors when we replace
%$M_{u\left(\Pi_u^l\right)}$ by $M_{u\left(\widehat{\Pi}_u^t\right)}$ and $M_{v\left(\Pi_v^l\right)}$ by $M_{v\left(\widehat{\Pi}_v^t\right)}$.
\begin{assumption}\label{ass:asnorm} Maintain the assumptions of Corollary \ref{p:opc} and Assumption \ref{ass:SF2b} and 
\begin{enumerate}[\textup{(}i\textup{)}]  
%\item [(i)] $r_N^{3/2}\lambda_N^2/\overline{v}_N=o\left(NT\right)$, 
\item\label{ass:asnormi} $\overline{r}_N\lambda_N^2\left(
\lambda_N+\sqrt{\overline{r}_N}\mu_N^2/\overline{v}_N\right)/\overline{v}_N=o\left(NT\right)$,
\item\label{ass:asnormii} $\overline{r}_N\lambda_N^3/\overline{v}_N=o\left(\sqrt{NT}\right)$, 
%$\overline{r}_N\lambda_N^2\max_{l=1,\dots,K}\left(\max\left(\left|X_lM_{v\left(\Pi_v^l\right)}\right|_{\op},\left|M_{u\left(\Pi_u^l\right)}X_l\right|_{\op}\right)+\sqrt{\overline{r}_N}\lambda_N\mu_N/v_N\right)/v_N=o_P\left(\sqrt{NT}\right)$,
\item\label{ass:asnormiii} $\overline{r}_N^{3/2}\lambda_N^3\left(|\Gamma|_{\op}+\lambda_N\right)/\overline{v}_N^2=o_P\left(\sqrt{NT}\right)$,
\item\label{ass:asnormiv} 
%%$\max_k\left|\left\langle \mathcal{P}_{\Pi^l}^\perp\left(U_k\right),\Gamma^d\right\rangle\right|=o_P\left(\sqrt{NT}\right)$,
$\left|\mathcal{P}_{\Pi^l}^\perp(\Pi^d)\right|_2^2=o_P(NT)$,
%
%\left|\Gamma^d\right|_{*}+\sum_{k=1}^K\left|\Pi_k^d\right|_*=o_P\left(\sqrt{NT}/\lambda_N\right)$,
\item\label{ass:asnormv} There exists $\Sigma_{\perp}\in\mathcal{M}_{KK}$ positive definite such that $\mathcal{P}_{\Pi^l}^\perp(U)^{\top}\mathcal{P}_{\Pi^l}^\perp(U)/(NT)\xrightarrow{\mathbb{P}}\Sigma_{\perp}$, 
\item\label{ass:asnormvi} $\mathcal{P}_{\Pi^l}^\perp(U)^\top e/\sqrt{NT}\xrightarrow{d}\mathcal{N}\left(0,\sigma^2\Sigma_{\perp}\right)$. 
\end{enumerate}
\end{assumption}
%A sufficient condition for Assumption \ref{ass:asnorm} (i) is 
%$r_N\lambda_N\mu_{N}^2/v_N=o\left(NT\right)$. We write a more complicated expression because $\mu_{uN}$ and $\mu_{vN}$ can be smaller than $\mu_N$ when the regressors and $\Gamma$ share factors or factor loadings 
%\blue{(je suis d'accord avec toi que c'est peu realiste, je l'ais remis car tu as recommence a parler de cela dans tes derniers commentaires...)} %\red{(Il faut que tous les facteurs forts soient tous communs dans ce cas, je pense )} 
%or if we work with transformed regressors as in Section \ref{s:transform}. A simple sufficient condition is
Regarding Assumption \ref{ass:asnorm} \eqref{ass:asnormiii}, $|\Gamma|_{\op}$ is usually $O_P\left(\sqrt{NT}\right)$ if it has a nontrivial low-rank component. \eqref{ass:asnormi}-\eqref{ass:asnormiii} can be satisfied under weaker assumptions than a strong-factor assumption ($\overline{v}_N$ is of the order of $\sqrt{NT}$) and when $\overline{r}_N$ goes to infinity. 
\eqref{ass:asnormv} is satisfied, for example, if $(\Pi_u^l,\Pi_v^l)$ and $U$ are independent and \eqref{ass:asnormvi} when $(X,\Gamma^l)$ and $E$ are independent.
\begin{theorem}\label{twostep} Let Assumption \ref{ass:asnorm} holds. We have 
\begin{align*}
&\frac{\sqrt{NT}}{\widehat{\sigma}}\left(\widetilde{\beta}^{(1)}-\beta\right)\xrightarrow{d}\mathcal{N}\left(0,\Sigma_{\perp}^{-1}\right),\\
&\mathcal{P}_{\widehat{\Pi}^t}^\perp(X)^{\top}\mathcal{P}_{\widehat{\Pi}^t}^\perp(X)/(NT)\xrightarrow{\mathbb{P}}\Sigma_{\perp}.
\end{align*} Also, if %$\Gamma,X$ and $E$ are independent and 
$\left|\mathcal{P}_{\Pi^l}(U)\right|_2^2=o_P(\left|U\right|_2^2)$ 
%\blue{(je suis pour enlever et fusionner la preuve avec celle de l'ennonce suivant)} 
then $\Sigma_\perp=\mathbb{E}[U^\top U]$. This occurs if $\mathbb{E}\left[\max\left(\rank\left(\Pi_u^l\right),\rank\left(\Pi_v^l\right)\right)\right]=o\left(\sqrt{\min(N,T)}\right)$ and $U$ and $(\Pi_u^l,\Pi_v^l)$ are independent.
% then $\Sigma_{\perp}=\mathbb{E}[U^\top U]$.
%\blue{(Une idee qui ne va peut etre pas marcher: Moreover, if the original regressors satisfy \eqref{efactorX} \blue{and ...} and we have applied transformation \eqref{eProj}, then $\Sigma_\perp$ is the same whether we use the transformed regressors or the orignal ones, finir la preuve si c'est possible)}.
\end{theorem}
%By the Pythagorean theorem, $\left|\mathcal{P}_{\Gamma^l}(X)\right|_2^2=\left|M_{u\left(\Pi_u^l\right)}XP_{v\left(\Gamma^l\right)}\right|_2^2+\left|P_{u\left(\Gamma^l\right)}X\right|_2^2$ so the last assumption is meaningful when $r_N/\min(N,T)\to0$. 
%The result that that $\Sigma_\perp=\Sigma$ means that the asymptotic covariance is the same as least-squares in a model where $\Gamma^l$ is known. When this holds, shows there is no loss in efficiency by using transformed regressors. %\blue{(La condition est-elle realiste?)}

\subsubsection{Approach 2: Using \cite{bai2009panel}'s estimator as a second stage}
An alternative approach discussed in \cite{Hsiao} is to rely on a preliminary consistent estimator to initialize \cite{bai2009panel}'s non convex estimator. \cite{moon2018nuclear} put forward this approach and the possibility to rely on a preliminary estimator like their matrix Lasso as a first-step. 
%to initialize \cite{bai2009panel}'s non convex estimator. % as justified in Section 5 in \cite{moon2018nuclear}. 
Among other conditions, using such a two-stage approach %\cite{moon2018nuclear} 
requires that the rate of convergence of the first-step estimator of $\beta$ is at least $(NT)^{1/6}$, a consistent estimator of $\rank(\Gamma)$, which is assumed constant, a strong-factor assumption on $\Gamma$, and $\Gamma^d=0$. % or of an upper bound on it \red{(Je reregarderai quand on aura fini le reste)}. 
This methodology can be applied using as a first-stage the thresholded or nonthresholded square-root estimator of this paper. We denote this estimator by $\left(\widetilde{\beta}^{(2)},\widetilde{\Gamma}^{(2)}\right)$. This paper provides a consistent estimator of $\rank\left(\Gamma^l\right)$ via hard-thresholding of \eqref{estGamma} or an upper bound on it without thresholding. Lemma 3 in \cite{moon2018nuclear} proposes an other consistent estimator but probably has a typo due to contradictory assumptions. The advantage of the estimator of this paper is that the level of thresholding is less conservative and makes use of the consistent estimator of the variance of errors. Recall  that if $\Gamma^d=0$ and $\Pi_1^l=\dots,\Pi_K^l$, from the discussion after Proposition \ref{l1} and \eqref{estGamma}, 
$$\rank\left(\widehat{\Gamma}\right)\le 
2\left(\frac{1+\rho}{1-\rho}\right)^2\left(\frac{\widetilde{r}_N}{\kappa_{\widetilde{\Gamma}^l}^2}+
\sum_{k=1}^Kr_{kN}\right)+o_P(1).$$ 
%The main advantage of this procedure is  that the second stage does not involve annihilation of low-rank components of the regressors and can yield a smaller asymptotic covariance matrix than the one in Theorem
%\ref{twostep} 
%However, it relies on  
An estimator of the asymptotic covariance matrix of the second-stage estimator, given a consistent estimator of $\widehat{r}=\rank\left(\Gamma^l\right)$, is given by % from \cite{bai2009panel}
(see page 1552 of \cite{moon2015linear})
$\widehat{\sigma}_B\widehat{\Sigma}_B$, where
\begin{align*}
\widehat{\sigma}_B&=\frac{1}{\sqrt{(N-\widehat{r})(T-\widehat{r})-K}}\left|Y-\sum_{k=1}^K\widetilde{\beta}_k^{(2)}X_k-\widetilde{\Gamma}^{(2)}\right|_2\\
\left(\widehat{\Sigma}_B\right)_{kl}&=\frac{1}{NT}\left\langle M_{u\left(\widetilde{\Gamma}^{(2)}\right)}X_kM_{v\left(\widetilde{\Gamma}^{(2)}\right)},X_l\right\rangle\quad \forall k,l\in\{1,\dots,K\}^2.
\end{align*}

\section{Simulations}
We take the same data generating process as in \cite{moon2018nuclear} with a single regressor and two factors:
\begin{align*}
Y_{it}&= X_{1it}+\sum_{l=1}^2\left(1+\lambda_{0,il}\right)f_{0,tl}+E_{it},\\
X_{1it}&=1+ \sum_{l=1}^2(2+\lambda_{0,il}+\lambda_{1,il})(f_{0,tl}+f_{0,t-1\ r})+U_{it},
\end{align*}
where $f_{0,tl}$, $\lambda_{0,il}$, $\lambda_{1,il}$, $U_{it}$, and $E_{it}$ for all indices are mutually independent and i.i.d. standard normal. The matrix $X_1$ has a statistical factor structure with a low-rank component of rank 3 due to the constant 1. Recall that $\widehat{\beta}^{LS}$ is the least-squares estimator which ignores the presence of $\Gamma$ is inconsistent because $X_{it}$ and $\Gamma_{it}$ are correlated. By the analysis of the paper, the square-root estimator coincides with the estimator in \cite{moon2018nuclear} with a smaller penalization. % if we were to take the theoretical one in \cite{moon2018nuclear}. 
The results in \cite{moon2018nuclear} are obtained with a penalty much smaller than allowed by the theory. %and we do not consider such values. %In Table~\ref{fig:Cov}, we present the bias and the standard deviation of various estimators from this paper computed
We compare the performance of the least-squares estimator $\widehat{\beta}^{LS}$, the square-root estimator $\widehat{\beta}$ obtained with the baseline regressors, the square-root estimator $\widehat{\beta}_{pt}$ obtained with the transformed regressors, where we apply (2) from Section \ref{s:transform} with $\widetilde{\Pi}_1=\widehat{\Pi}_1^t$, and the two-stage estimators from Section \ref{sec:modif}. We use $\widehat{\beta}^{LS}$
%the least-square estimator of $\beta$ 
to initialize the iterative estimators. The number of iterations is 100 to obtain the estimator of $\rank\left(\Gamma\right)$, as explained after Corollary \ref{p:opc}, useful to compute $\widetilde{\beta}^{(2)}$. We use the same number of iterations to obtain $\widehat{\beta}_{pt}$. We consider an additional 100 iterations for $\widehat{\beta}$, $\widehat{\beta}_{pt}$, and $\widetilde{\beta}^{(2)}$. As a result, $\widetilde{\beta}^{(1)}$ and $\widetilde{\beta}^{(2)}$ have been computed with the same number of iterations. We consider two sample sizes: (a) $N=T=50$ and (b) $N=T=150$. We use $7300$ Monte-Carlo replications to allow for an accuracy of $\pm0.005$ with 95\% for the coverage probabilities of 95\% confidence intervals. We choose $\lambda_N=1.01\left(\sqrt{N}+\sqrt{T}\right)$ and the hard-thresholding levels are $2\lambda_N$ times an estimator of the standard error from the first-stages. 

A first approach is to not apply any matrix to the data as described after Proposition \ref{proposal}. The results in tables \ref{fig:Cov} and \ref{fig:Cov1} compare the performance of the estimators in terms of MSE, bias, and standard error (henceforth std). In case (a), $\rank\left(\widehat{\Pi}_1^t\right)$ has been found to be always equal to 2 while $\rank\left(\widehat{\Pi}_1\right)$ to 3 (the true rank), $\rank\left(\widehat{\Gamma}^t\right)$ has been found to be always equal to 2 (the true rank) in 89\% of the cases and else to 1. % in 11\% of the cases. 
We used $\rank\left(\widehat{\Pi}_1^t\right)$ for $\widehat{\beta}_{pt}$ and subsequently  $\rank\left(\widehat{\Gamma}^t\right)$, $\widetilde{\beta}^{(1)}$ and $\widetilde{\beta}^{(2)}$, even though it did not perform well for such small sample size.   In case (b), $\rank\left(\widehat{\Pi}_1^t\right)$ has been found to be always equal to 3 while $\rank\left(\widehat{\Pi}_1\right)$ and $\rank\left(\widehat{\Gamma}^t\right)$ have been found to be always equal to 2 (the true rank). 
%We used $\rank\left(\widehat{\Pi}_1^t\right)$ for $\widehat{\beta}_{pt}$ and subsequently  $\rank\left(\widehat{\Gamma}^t\right)$, $\widetilde{\beta}^{(1)}$ and $\widetilde{\beta}^{(2)}$, even though it did not performed well for such a small sample size.  
\begin{table}[!ht]
\begin{minipage}{0.48\textwidth}
  \centering
  \caption{$N=T=50$}
\label{fig:Cov} 
{\small
             \begin{tabular}{|c|c|c|c|c|c|}
  \hline
 & $\widehat{\beta}^{LS}$  & $\widehat{\beta}$& $\widehat{\beta}_{pt}$& $\widetilde{\beta}^{(1)}$ & $\widetilde{\beta}^{(2)}$ \\
  \hline
   MSE  & 0.053 & 0.020 & 5 $10^{-4}$&  0.002 & 9 $10^{-4}$ \\
   bias  & 0.230 & 0.142 & -$10^{-4}$&  0.019& 0.009\\
   std  &   0.017 & 0.015 & 0.023 &   0.035&  0.029 \\
  \hline
      \end{tabular}  }
%\end{table}
%\begin{table}[!ht]
\end{minipage}
\begin{minipage}{0.48\textwidth}
  \centering
  \caption{$N=T=150$}
    \label{fig:Cov1}
    {\small
             \begin{tabular}{|c|c|c|c|c|c|}
  \hline
 & $\widehat{\beta}^{LS}$ & $\widehat{\beta}$& $\widehat{\beta}_{pt}$& $\widetilde{\beta}^{(1)}$ & $\widetilde{\beta}^{(2)}$ \\
  \hline
   MSE  & 0.053 & 0.011 & 4 $10^{-5}$&  4 $10^{-5}$ & 1 $10^{-5}$ \\
   bias  & 0.231 & 0.103 & 4 $10^{-4}$&  2 $10^{-5}$ & -8 $10^{-5}$\\
   std  &   0.009 & 0.008 & 0.006 &   0.006&  0.003 \\
  \hline
        \end{tabular}  
%               \begin{tabular}{|c|c|c|c|c|c|}
%  \hline
% & $\widehat{\beta}^{LS}$ & $\widehat{\beta}$& $\widehat{\beta}_{pt}$& $\widetilde{\beta}^{(1)}$ & $\widetilde{\beta}^{(2)}$ \\
%  \hline
%   MSE  & 0.053 & 0.009 & 3 $10^{-5}$&  2 $10^{-5}$ & 6 $10^{-6}$ \\
%   bias  & 0.231 & 0.095 & 7 $10^{-5}$&  2 $10^{-4}$ & 7 $10^{-5}$\\
%   std  &   0.009 & 0.007 & 0.005 &   0.005&  0.003 \\
%  \hline
%        \end{tabular}        
        }
\end{minipage}
\end{table}

A second approach is to apply Within transforms $M_u=I_N-J_N/N$ and $M_v=I_T-J_T/T$ to the left and right of $Y$ and $X_1$, where  $J_N\in\mathcal{M}_{NN}$ (resp. $J_T\in\mathcal{M}_{TT}$) has all entries equal to 1. These allow to get rid of the mean 1 of $X_1$ but more generally of any individual and time effects in both $\Pi^l$ and $\Gamma^l$. %It can be expected to be useful because $\rank(\Pi_1)$ has not been well estimated in case (a) above. 
The results are in tables \ref{fig:Covw} and \ref{fig:Cov1w}. In case (a), $\rank\left(\widehat{\Pi}_1^t\right)$ and $\rank\left(\widehat{\Pi}_1\right)$ has been found to be always equal to 2 (the true rank), $\rank\left(\widehat{\Gamma}^t\right)$ has been found to be equal to 2 (the true rank) in 81\% of the cases and else to 1. % in 19\% of the cases. 
%We used $\rank\left(\widehat{\Pi}_1^t\right)$ for $\widehat{\beta}_{pt}$ and subsequently  $\rank\left(\widehat{\Gamma}^t\right)$, $\widetilde{\beta}^{(1)}$ and $\widetilde{\beta}^{(2)}$, even though it did not performed well for such a small sample size.  
In case (b), $\rank\left(\widehat{\Pi}_1^t\right)$, $\rank\left(\widehat{\Pi}_1\right)$, $\rank\left(\widehat{\Gamma}^t\right)$ have been found to be always equal to 2 (the true ranks). 
%We used $\rank\left(\widehat{\Pi}_1^t\right)$ for $\widehat{\beta}_{pt}$ and subsequently  $\rank\left(\widehat{\Gamma}^t\right)$, $\widetilde{\beta}^{(1)}$ and $\widetilde{\beta}^{(2)}$, even though it did not performed well for such a small sample size.  

\begin{table}[!ht]
\begin{minipage}{0.48\textwidth}
  \centering
  \caption{$N=T=50$, Within}
\label{fig:Covw}
{\small
                       \begin{tabular}{|c|c|c|c|c|c|}
  \hline
 & $\widehat{\beta}^{LS}$  & $\widehat{\beta}$& $\widehat{\beta}_{pt}$& $\widetilde{\beta}^{(1)}$ & $\widetilde{\beta}^{(2)}$ \\
  \hline
   MSE  & 0.049 & 0.016 & 5 $10^{-4}$ &  0.001 & 0.002 \\
   bias  & 0.221 & 0.124 & -4 $10^{-5}$ &  0.024 & 0.020\\
   std  &   0.025 & 0.018 & 0.023 &   0.025 &  0.044 \\
  \hline
      \end{tabular}    }
%\end{table}
%\begin{table}[!ht]
\end{minipage}
\begin{minipage}{0.48\textwidth}
\centering
\caption{$N=T=150$, Within}
  \label{fig:Cov1w}
    {\small
    \begin{tabular}{|c|c|c|c|c|c|}
  \hline
 & $\widehat{\beta}^{LS}$ & $\widehat{\beta}$& $\widehat{\beta}_{pt}$& $\widetilde{\beta}^{(1)}$ & $\widetilde{\beta}^{(2)}$ \\
  \hline
   MSE  & 0.049 & 0.007 & 5 $10^{-5}$&  5 $10^{-5}$ & 2 $10^{-5}$ \\
   bias  & 0.222 & 0.081 & -1 $10^{-4}$&  4 $10^{-4}$ & 7 $10^{-7}$\\
   std  &   0.014 & 0.007 & 0.007 &   0.007&  0.004 \\
  \hline
      \end{tabular}    
%        \caption{$N=T=200$ using within transforms}
%  \label{fig:Cov1w}
%    {\small
%    \begin{tabular}{|c|c|c|c|c|c|}
%  \hline
% & $\widehat{\beta}^{LS}$ & $\widehat{\beta}$& $\widehat{\beta}_{pt}$& $\widetilde{\beta}^{(1)}$ & $\widetilde{\beta}^{(2)}$ \\
%  \hline
%   MSE  & 0.049 & 0.005 & 3 $10^{-5}$&  3 $10^{-5}$ & 9 $10^{-6}$ \\
%   bias  & 0.222 & 0.073 & -9 $10^{-6}$&  -3 $10^{-6}$ & 5 $10^{-5}$\\
%   std  &   0.012 & 0.006 & 0.005 &   0.005&  0.003 \\
%  \hline
%      \end{tabular}     
      }
\end{minipage}
\end{table}
Table \ref{fig:Cov2w} assesses the coverage probabilities in the different cases.
\begin{table}[!ht]
        \caption{Coverage of 95\% confidence intervals based on the two-stage approaches.}
            \label{fig:Cov2w}
   \center {\small \begin{tabular}{|c|c|c|c|}
  \hline
  Within transforms & (N,T) & $\widetilde{\beta}^{(1)}$ & $\widetilde{\beta}^{(2)}$\\
  \hline
 No &  (50,50) & %0.866 
 0.87 & %0.844
 0.84\\
 Yes &  (50,50) & 
 0.81 & 
 0.76\\
  No &  (150,150) & 0.95 %0.953
   & %0.949
   0.94\\
%     No &  (200,200) & 0.95 %0.953
%   & %0.949
%   0.95\\
  Yes &  (150,150) & 0.95 &  0.94\\
%  No &  (200,200) & 0.95 %0.953
%   & %0.949
%   0.95\\
%  Yes &  (200,200) & %0.954 
%  0.95 & %0.946
%  0.95\\
   \hline
      \end{tabular}    }
\end{table}

\bibliographystyle{abbrv}
\bibliography{refIFE}

%\appendix
%\appendixpage
\section*{Appendix}
%\blue{Il faudra ne garder que ce dont on a besoin\\}
%We start by recalling a few result from linear algebra.  

\subsubsection*{Proof of Proposition \ref{equiv}}
\noindent By definition of $\widehat{\beta}$ and $\widehat{\Gamma}$, we have, for all $\beta\in \mathbb{R}^K$ and $\Gamma \in \mathcal{M}_{NT}$, 
$$\frac{1}{\sqrt{NT}}\left\vert  Y-\sum_{k=1}^K\widehat{\beta}_k X_k-\widehat{\Gamma}  \right \vert_2
+\frac{\lambda}{NT} \left|\widehat{\Gamma}\right|_*\le \frac{1}{\sqrt{NT}}\left\vert  Y-\sum_{k=1}^K\beta_k X_k-\Gamma  \right \vert_2
+\frac{\lambda}{NT} \norm{\Gamma}_*.$$
By definition of $P_X$ and of the estimator, for all $\beta\in \mathbb{R}^K$ and $\Gamma \in \mathcal{M}_{NT}$, we have 
\begin{align*}\frac{1}{\sqrt{NT}}\left\vert  M_X\left(Y-\widehat{\Gamma}\right)  \right \vert_2
+\frac{\lambda}{NT} \left|\widehat{\Gamma}\right|_*&\le \frac{1}{\sqrt{NT}}\left\vert  Y-\sum_{k=1}^K\widehat{\beta}_k X_k-\widehat{\Gamma}  \right \vert_2
+\frac{\lambda}{NT} \left|\widehat{\Gamma}\right|_*\\
%\frac{1}{\sqrt{NT}}\left\vert  M_X\left(Y-\widehat{\Gamma}\right)  \right \vert_2
%+\frac{\lambda}{NT} \left|\widehat{\Gamma}\right|_*
&\le \frac{1}{\sqrt{NT}}\left\vert  Y-\sum_{k=1}^K\beta_k X_k-\Gamma  \right \vert_2
+\frac{\lambda}{NT} \norm{\Gamma}_*.
\end{align*}
 By choosing %, for given $\Gamma$, 
 $\beta$ such that 
$\sum_{k=1}^K\beta_k X_k=P_X\left(Y-\Gamma\right)$, % in the right-hand side, 
we obtain, for all $\Gamma \in \mathcal{M}_{NT}$, 
$$\frac{1}{\sqrt{NT}}\left\vert  M_X\left(Y-\widehat{\Gamma}\right)  \right \vert_2
+\frac{\lambda}{NT} \left|\widehat{\Gamma}\right|_*\le \frac{1}{\sqrt{NT}}\left\vert  M_X\left(Y-\Gamma  \right)\right\vert_2
+\frac{\lambda}{NT} \norm{\Gamma}_*,$$ hence the result.
%Therefore, it holds that 
%\begin{equation*}
%\widehat{\Gamma}
%\in 
%\arg \min_{\Gamma\in\mathcal{M}_{NT}} \frac{1}{\sqrt{NT}}\left\vert M_X\left( Y-\Gamma\right)  \right \vert_2
%+\frac{\lambda}{NT} \norm{\Gamma}_*.
%\end{equation*}

\subsubsection*{Proof of Proposition\ref{pkappa}}
\noindent The first inequality is obtained using trace duality and \eqref{eq:rank}. The second is obtained by \eqref{orth} and the Pythagorean theorem.

\subsubsection*{Proof of Theorem \ref{orac}}
\noindent The techniques are similar to those in \cite{belloni2011square,gautier2011high}. 
Take $\widetilde{\Gamma}\in\mathcal{M}_{NT}$ and denote by $\Delta=\widehat{\Gamma}-\Gamma$. Remark that
\begin{align}
 \notag \left|\widehat{\Gamma}\right|_*
 %&=\norm{\Gamma+\Delta}_*\\
% \notag&=\left|\Gamma-\widetilde{\Gamma}+\widetilde{\Gamma}+\Delta\right|_*\\
%\notag
&=\left|\Gamma-\widetilde{\Gamma}+\widetilde{\Gamma}+\mathcal{P}_{\widetilde{\Gamma}}(\Delta)+\mathcal{P}_{\widetilde{\Gamma}}^{\perp}(\Delta)\right|_*\\
\label{finally2}&\ge \left|\widetilde{\Gamma}+\mathcal{P}_{\widetilde{\Gamma}}^{\perp}(\Delta)\right|_*-\left|\Gamma-\widetilde{\Gamma}\right|_*-\left|\mathcal{P}_{\widetilde{\Gamma}}(\Delta)\right|_*\\
\label{Trois2}%\left|\widehat{\Gamma}\right|_*
&\ge \left|\widetilde{\Gamma}\right|_*+\left|\mathcal{P}_{\widetilde{\Gamma}}^{\perp}(\Delta)\right|_*-\left|\Gamma-\widetilde{\Gamma}\right|_*-\left|\mathcal{P}_{\widetilde{\Gamma}}(\Delta)\right|_*\quad\text{(by \eqref{add})}. 
\end{align}
 Now, by \eqref{eq:apend} and the definition of $\widehat{\Gamma}$, we have
 \begin{equation}\label{minimizer2}  \frac{1}{\sqrt{NT}} \left\vert M_X\left(Y-\widehat{\Gamma}\right)\right \vert_2 +\frac{\lambda}{NT} \left|\widehat{\Gamma}\right|_*\le \frac{1}{\sqrt{NT}}\left\vert M_X\left(Y-\Gamma \right)\right \vert_2+\frac{\lambda}{NT} \left|\Gamma\right|_*. \end{equation}
By convexity, trace duality, and $\lambda\rho\norm{M_X(E)}_2/\sqrt{NT}\ge \norm{M_X(E)}_{\op}$, if $M_X(E) \ne 0$, we have
     \begin{align}
     \notag \frac{1}{\sqrt{NT}}\left\vert  M_X\left(Y-\widehat{\Gamma}\right)\right \vert_2-\frac{1}{\sqrt{NT}}\left\vert M_X\left(Y-\Gamma \right)\right \vert_2
%     \\
%     \notag 
     &\ge  -\frac{1}{\sqrt{NT}\norm{M_X(E)}_2}\left\langle M_X(E),\widehat{\Gamma}-\Gamma\right\rangle \\
%     \notag &\ge -\frac{\norm{M_X(E)}_{\op}}{\sqrt{NT}\norm{M_X(E)}_2}\norm{\Delta}_*\\
    \label{ineq2}  &\ge-\frac{\lambda\rho}{NT}  \norm{\Delta}_*.
     \end{align}
 \eqref{ineq2} also holds if $M_X(E)=0$ because $\left\vert  M_X\left(Y-\widehat{\Gamma}\right)\right \vert_2\ge 0$.
 This and \eqref{minimizer2} imply
  \begin{equation}\label{biz}\left|\widehat{\Gamma}\right|_*\le \rho\norm{\Delta}_*+\norm{\Gamma}_*.\end{equation}
Using \eqref{Trois2}, we get
  $$\left|\widetilde{\Gamma}\right|_*+\left|\mathcal{P}_{\widetilde{\Gamma}}^{\perp}(\Delta)\right|_*-\left|\Gamma-\widetilde{\Gamma}\right|_*-\left|\mathcal{P}_{\widetilde{\Gamma}}(\Delta)\right|_*\le \rho\left|\Delta\right|_*+\left|\Gamma\right|_*$$
and $\left|\Gamma\right|_*\le \left|\Gamma-\widetilde{\Gamma}\right|_*+\left|\widetilde{\Gamma}\right|_*$ yields 
   $$\left|\mathcal{P}_{\widetilde{\Gamma}}^{\perp}(\Delta)\right|_*-\left|\mathcal{P}_{\widetilde{\Gamma}}(\Delta)\right|_*\le \rho\left|\Delta\right|_*+2\left|\Gamma-\widetilde{\Gamma}\right|_*.$$
Then, because $\left|\Delta\right|_*\le\left|\mathcal{P}_{\widetilde{\Gamma}}^{\perp}(\Delta)\right|_*+\left|\mathcal{P}_{\widetilde{\Gamma}}(\Delta)\right|_*$, we have
  \begin{equation} \label{newcone2} (1-\rho)\left|\mathcal{P}_{\widetilde{\Gamma}}^{\perp}(\Delta)\right|_*\le (1+\rho) \left|\mathcal{P}_{\widetilde{\Gamma}}(\Delta)\right|_*+2 \left|\Gamma-\widetilde{\Gamma}\right|_*. \end{equation}
Also, by \eqref{minimizer2},
$$\frac{1}{\sqrt{NT}} \left\vert M_X\left(Y-\widehat{\Gamma}\right)\right \vert_2-\frac{1}{\sqrt{NT}}\left\vert M_X\left(Y-\Gamma \right)\right \vert_2 \le \frac{\lambda}{NT}\left(\left|\Gamma\right|_*-\left|\widehat{\Gamma}\right|_*\right)$$
and
\begin{align}
\left|\Gamma\right|_*-\left|\widehat{\Gamma}\right|_*&\le\left|\widetilde{\Gamma}\right|_*+\left|\Gamma-\widetilde{\Gamma}\right|_*-\left|\widehat{\Gamma}\right|_*\notag \\
&= 2\left|\Gamma-\widetilde{\Gamma}\right|_*+ \left|\widetilde{\Gamma}\right|_*-\left|\Gamma-\widetilde{\Gamma}\right|_*-\left|\widehat{\Gamma}\right|_* \notag\\
&\le 2\left|\Gamma-\widetilde{\Gamma}\right|_*+\left|\mathcal{P}_{\widetilde{\Gamma}}(\Delta)\right|_*-\left|\mathcal{P}_{\widetilde{\Gamma}}^{\perp}(\Delta)\right|_* \quad\text{(by \eqref{Trois2})} \notag,
\end{align}
hence we have
\begin{equation}\label{nref}
\frac{1}{\sqrt{NT}} \left\vert M_X\left(Y-\widehat{\Gamma}\right)\right \vert_2-\frac{1}{\sqrt{NT}}\left\vert M_X\left(Y-\Gamma \right)\right \vert_2 \le \frac{\lambda}{NT}\left(2\left|\Gamma-\widetilde{\Gamma}\right|_*+\left|\mathcal{P}_{\widetilde{\Gamma}}(\Delta)\right|_*\right).
\end{equation}
\noindent Let $\widetilde{\rho}>0$ and consider two cases.\\
\emph{Case 1.} If $\widetilde{\rho}\left|\mathcal{P}_{\widetilde{\Gamma}}(\Delta)\right|_*\le 2\left|\Gamma-\widetilde{\Gamma}\right|_* $, we have, by \eqref{newcone2}, 
   \begin{align}
   \notag \left|\Delta\right|_*&\le \left|\mathcal{P}_{\widetilde{\Gamma}}(\Delta)\right|_*+\left|\mathcal{P}_{\widetilde{\Gamma}}^{\perp}(\Delta)\right|_*\\
   \notag &\le \frac{2}{1-\rho}\left(\left|\mathcal{P}_{\widetilde{\Gamma}}(\Delta)\right|_*+\left|\Gamma-\widetilde{\Gamma}\right|_*\right)\\
   \notag &\le \frac{2}{1-\rho}\left(\frac{2}{\widetilde{\rho}}+1\right)
   \left|\Gamma-\widetilde{\Gamma}\right|_*.
   \end{align}
This yields the first part of the first inequality of Theorem \ref{orac}. The first part of the second inequality is obtained by combining \eqref{ineq2} and \eqref{nref}.\\
\emph{Case 2.} Otherwise, if $\widetilde{\rho}\left|\mathcal{P}_{\widetilde{\Gamma}}(\Delta)\right|_*> 2\left|\Gamma-\widetilde{\Gamma}\right|_* $, we obtain, by \eqref{newcone2}, that
   $$\left|\mathcal{P}_{\widetilde{\Gamma}}^{\perp}(\Delta)\right|_*\le c\left(\rho,\widetilde{\rho}\right)\left|\mathcal{P}_{\widetilde{\Gamma}}(\Delta)\right|_*, $$
which implies that $\Delta\in C_{\widetilde{\Gamma}}$ and $\norm{\Delta}_*\le\left(1+c\left(\rho,\widetilde{\rho}\right)\right)
\left|\mathcal{P}_{\widetilde{\Gamma}}(\Delta)\right|_*$.
We have
$$\frac{1}{NT}\left \vert M_X\left(Y-\widehat{\Gamma}\right) \right \vert_2^2-\frac{1}{NT}\norm{M_X\left(Y-\Gamma\right)}_2^2= \frac{1}{NT}\left \vert M_X\left(\widehat{\Gamma}-\Gamma\right) \right \vert_2^2-\frac{2}{NT}\left\langle M_X(E),\widehat{\Gamma}-\Gamma\right\rangle $$
hence, because $\lambda\rho\norm{M_X(E)}_2/\sqrt{NT}\ge \norm{M_X(E)}_{\op}$, 
%$$ \frac{1}{NT}\left\vert M_X\left(\widehat{\Gamma}-\Gamma \right)  \right \vert_2^2\le 
% \frac{1}{NT}\left \vert M_X\left(Y-\widehat{\Gamma}\right) \right \vert_2^2-\frac{1}{NT}\norm{M_X\left(Y-\Gamma\right)}_2^2+2\lambda\rho\frac{\norm{M_X(E)}_2}{(NT)^{(\rho+1)}}  \norm{\Delta}_*.$$
%As $\Delta \in C_{\widetilde{\Gamma}}$, hence
 \begin{align}
\notag &\frac{1}{NT}\left\vert  M_X\left(\widehat{\Gamma}-\Gamma \right)  \right \vert_2^2\\ 
%\notag 
%&\le 
% \frac{1}{NT}\left \vert  M_X\left(Y-\widehat{\Gamma}\right) \right \vert_2^2-\frac{1}{NT}\norm{M_X\left(Y-\Gamma\right)}_2^2+2\lambda\rho\frac{\norm{M_X(E)}_2}{(NT)^{(\rho+1)}}\left(  \left|\mathcal{P}_{\widetilde{\Gamma}}(\Delta)\right|_*+\left|\mathcal{P}_{\widetilde{\Gamma}}^{\perp}(\Delta)\right|_*\right) \\
 \label{mid2}& \le  \frac{1}{NT}\left \vert  M_X\left(Y-\widehat{\Gamma}\right)\right \vert _2^2-\frac{1}{NT}\norm{M_X\left(Y-\Gamma\right)}_2^2+2\lambda\rho\left(1+c\left(\rho,\widetilde{\rho}\right)\right)\frac{\norm{M_X(E)}_2}{(NT)^{\frac{3}{2}}}\left|\mathcal{P}_{\widetilde{\Gamma}}(\Delta)\right|_*
 \end{align}
and, by \eqref{nref},
% \blue{debut copie\\}
% Also, we have, by \eqref{minimizer2},
%$$\frac{1}{\sqrt{NT}} \left\vert M_X\left(Y-\widehat{\Gamma}\right)\right \vert_2-\frac{1}{\sqrt{NT}}\left\vert M_X\left(Y-\Gamma \right)\right \vert_2 \le \frac{\lambda}{NT}\left(\left|\Gamma\right|_*-\left|\widehat{\Gamma}\right|_*\right).$$
%Next, remark that 
%\begin{align*}
%\left|\Gamma\right|_*-\left|\widehat{\Gamma}\right|_*&\le\left|\widetilde{\Gamma}\right|_*+\left|\Gamma-\widetilde{\Gamma}\right|_*-\left|\widehat{\Gamma}\right|_* \\
%&= 2\left|\Gamma-\widetilde{\Gamma}\right|_*+ \left|\widetilde{\Gamma}\right|_*-\left|\Gamma-\widetilde{\Gamma}\right|_*-\left|\widehat{\Gamma}\right|_* \\
%&\le 2\left|\Gamma-\widetilde{\Gamma}\right|_*+\left|\mathcal{P}_{\widetilde{\Gamma}}(\Delta)\right|_*-\left|\mathcal{P}_{\widetilde{\Gamma}}^{\perp}(\Delta)\right|_* \quad\text{(by \eqref{Trois2})}\\
%&\le (1+\widetilde{\rho})\left|\mathcal{P}_{\widetilde{\Gamma}}(\Delta)\right|_*-\left|\mathcal{P}_{\widetilde{\Gamma}}^{\perp}(\Delta)\right|_* \quad\left(\text{because $2\left|\Gamma-\widetilde{\Gamma}\right|_* \le \widetilde{\rho}\left|\mathcal{P}_{\widetilde{\Gamma}}(\Delta)\right|_*$}\right)\\
%&\le (1+\widetilde{\rho})\left|\mathcal{P}_{\widetilde{\Gamma}}(\Delta)\right|_*.
%\end{align*}
%Therefore, we obtain 
$$\frac{1}{\sqrt{NT}}\left\vert  M_X\left(Y-\widehat{\Gamma}\right)\right \vert_2-\frac{1}{\sqrt{NT}}\left\vert  M_X\left(Y-\Gamma \right)\right \vert_2 \le  \frac{(1+\widetilde{\rho})\lambda}{NT}\left|\mathcal{P}_{\widetilde{\Gamma}}(\Delta)\right|_*$$
which, combined with \eqref{ineq2}, yields
% and the fact that $\norm{\Delta}_*\le\left(1+c\left(\rho,\widetilde{\rho}\right)\right)\left|\mathcal{P}_{\widetilde{\Gamma}}(\Delta)\right|_*$, yields []
%$$-\frac{\rho\lambda}{NT}\norm{\Delta}_*\le\frac{1}{\sqrt{NT}}\left\vert  M_X\left(Y-\widehat{\Gamma}\right)\right \vert_2-\frac{1}{\sqrt{NT}}\left\vert  M_X\left(Y-\Gamma \right)\right \vert_2 \le  \frac{(1+\widetilde{\rho})\lambda}{NT}\left|\mathcal{P}_{\widetilde{\Gamma}}(\Delta)\right|_*.$$
%Then, as $\Delta\in C_{\widetilde{\Gamma}}$, we have   $-\frac{\rho\lambda}{NT}\norm{\Delta}_*\ge- \frac{\rho\lambda}{NT}\left(\left|\mathcal{P}_{\widetilde{\Gamma}}(\Delta)\right|_*+\norm{\mathcal{P}_{\Gamma}^{\perp}(\Delta)}_*\right)\ge- \frac{(2+\widetilde{\rho})\rho\lambda}{(1-\rho)NT}\left|\mathcal{P}_{\widetilde{\Gamma}}(\Delta)\right|_*$, hence
\begin{align}\label{eq:varianceest}
\left \vert\frac{1}{\sqrt{NT}}\left\vert M_X\left(Y-\widehat{\Gamma}\right)\right \vert_2-\frac{1}{\sqrt{NT}}\left\vert M_X\left(Y-\Gamma \right)\right \vert_2\right \vert \le d\left(\rho,\widetilde{\rho}\right)
%\max\left(1+\widetilde{\rho},\rho\frac{2+\widetilde{\rho}}{1-\rho}\right)
\frac{\lambda}{NT}\left|\mathcal{P}_{\widetilde{\Gamma}}(\Delta)\right|_*.
\end{align}
Now, using 
\begin{align} \notag &\frac{1}{NT} \left \vert M_X\left(Y-\widehat{\Gamma}\right) \right \vert_2^2-\frac{1}{NT}\norm{M_X\left(Y-\Gamma\right)}_2^2\\ 
%\notag & =\left(\frac{1}{\sqrt{NT}} \left\vert M_X\left(Y-\widehat{\Gamma}\right)\right \vert_2-\frac{1}{\sqrt{NT}}\left\vert M_X\left(Y-\Gamma \right)\right \vert_2\right)\left(\frac{1}{\sqrt{NT}}\left\vert  M_X\left(Y-\widehat{\Gamma}\right)\right \vert_2+\frac{1}{\sqrt{NT}}\left\vert M_X\left(Y-\Gamma \right)\right \vert_2\right)\\
  \notag &  =\left(\frac{1}{\sqrt{NT}} \left\vert M_X\left(Y-\widehat{\Gamma}\right)\right \vert_2-\frac{1}{\sqrt{NT}}\left\vert M_X\left(Y-\Gamma \right)\right \vert_2\right)\\
 &\notag \quad\times\left(\frac{1}{\sqrt{NT}} \left\vert M_X\left(Y-\widehat{\Gamma}\right)\right \vert_2-\frac{1}{\sqrt{NT}}\left\vert M_X\left(Y-\Gamma \right)\right \vert_2+\frac{2}{\sqrt{NT}} \left\vert M_X\left(Y-\Gamma \right)\right \vert_2\right)
 \end{align}
and \eqref{eq:varianceest}, we obtain
 \begin{align}
  \label{finu2} &\frac{1}{NT} \left \vert M_X\left(Y-\widehat{\Gamma}\right) \right \vert_2^2-\frac{1}{NT}\norm{M_X\left(Y-\Gamma\right)}_2^2\\
% \notag &  \le \max\left(1+\widetilde{\rho},\rho\frac{2+\widetilde{\rho}}{1-\rho}\right)\frac{\lambda}{NT}\left|\mathcal{P}_{\widetilde{\Gamma}}(\Delta)\right|_*\left(\max\left(1+\widetilde{\rho},\rho\frac{2+\widetilde{\rho}}{1-\rho}\right)\frac{\lambda}{NT}\left|\mathcal{P}_{\widetilde{\Gamma}}(\Delta)\right|_*+\frac{2}{\sqrt{NT}} \left\vert M_X\left(Y-\Gamma \right)\right \vert_2\right)\\
 \notag &  \le d\left(\rho,\widetilde{\rho}\right)\frac{\lambda}{NT}\left|\mathcal{P}_{\widetilde{\Gamma}}(\Delta)\right|_*\left(d\left(\rho,\widetilde{\rho}\right)\frac{\lambda}{NT}\left|\mathcal{P}_{\widetilde{\Gamma}}(\Delta)\right|_*+\frac{2\left\vert M_X\left(E \right)\right \vert_2}{\sqrt{NT}}\right).\end{align} 
Combining \eqref{mid2} and \eqref{finu2}, we get
\begin{align*}
&\frac{1}{NT}\left \vert M_X\left(\widehat{\Gamma}-\Gamma \right) \right \vert_2^2 \le 
 \left(
 d\left(\rho,\widetilde{\rho}\right) \frac{\lambda}{NT}\left|\mathcal{P}_{\widetilde{\Gamma}}(\Delta)\right|_*\right)^2+2e\left(\rho,\widetilde{\rho}\right)\frac{\lambda\left\vert  M_X\left(E \right)\right \vert_2}{(NT)^{\frac{3}{2}}}\left|\mathcal{P}_{\widetilde{\Gamma}}(\Delta)\right|_*.
\end{align*}
By definition of  $\kappa_{\widetilde{\Gamma},c\left(\rho,\widetilde{\rho}\right)}$, this implies %that 
%$$\left\vert M_X\left(\Delta \right) \right \vert_2^2\le
%  \left(\frac{4\sqrt{2\mathrm{rank}\left(\widetilde{\Gamma}\right)}\lambda}{\sqrt{NT}\kappa_{\widetilde{\Gamma}}}\right)^2
%\left\vert M_X\left(\Delta \right)  \right \vert_2^2+\frac{16\lambda\sqrt{2\mathrm{rank}\left(\widetilde{\Gamma}\right)}\left\vert M_X\left(E \right)\right \vert_2}{\sqrt{NT}\kappa_{\widetilde{\Gamma}}}\left\vert  M_X\left(\Delta \right)  \right \vert_2 $$
%Therefore, we obtain
\begin{align}
\left\vert M_X\left(\Delta \right)  \right \vert_2&\le 2\left(1-   \left(d\left(\rho,\widetilde{\rho}\right)\frac{\sqrt{2\mathrm{rank}\left(\widetilde{\Gamma}\right)}\lambda}{\sqrt{NT}\kappa_{\widetilde{\Gamma},c\left(\rho,\widetilde{\rho}\right)}}\right)^2\right)_+^{-1}
e\left(\rho,\widetilde{\rho}\right)\frac{\lambda\sqrt{2\mathrm{rank}\left(\widetilde{\Gamma}\right)}\norm{M_X(E)}_2}{\sqrt{NT}\kappa_{\widetilde{\Gamma},c\left(\rho,\widetilde{\rho}\right)}},\notag\\
\label{eq:pnuc}
\left\vert \mathcal{P}_{\Gamma}(\Delta) \right \vert_*&\le 4\left(1-   \left(d\left(\rho,\widetilde{\rho}\right)\frac{\sqrt{2\mathrm{rank}\left(\widetilde{\Gamma}\right)}\lambda}{\sqrt{NT}\kappa_{\widetilde{\Gamma},c\left(\rho,\widetilde{\rho}\right)}}\right)^2\right)_+^{-1}e\left(\rho,\widetilde{\rho}\right)\frac{\lambda\mathrm{rank}\left(\widetilde{\Gamma}\right)\norm{M_X(E)}_2}{\sqrt{NT}\kappa_{\widetilde{\Gamma},c\left(\rho,\widetilde{\rho}\right)}^2},%\\
%\left\vert\Delta \right \vert_*&\le4\left(1+c\left(\rho,\widetilde{\rho}\right)\right) \left(1-   \left(\frac{d\left(\rho,\widetilde{\rho}\right)\sqrt{2\mathrm{rank}\left(\widetilde{\Gamma}\right)}\lambda}{\sqrt{NT}\kappa_{\widetilde{\Gamma},c\left(\rho,\widetilde{\rho}\right)}}\right)^2\right)_+^{-1}e\left(\rho,\widetilde{\rho}\right)\frac{\lambda\mathrm{rank}\left(\widetilde{\Gamma}\right)\norm{M_X(E)}_2}{\sqrt{NT}\kappa_{\widetilde{\Gamma},c\left(\rho,\widetilde{\rho}\right)}^2},\notag
\end{align}
which yields the first result.
The second result follows from \eqref{eq:varianceest} and \eqref{eq:pnuc}.

\subsubsection*{Proof of Proposition \ref{sigmacons}}
%\noindent We start by proving the following result.
\begin{lemma}\label{tech1}It holds that $\left|P_X(E)\right|_2= O_P(1)$ and $\left|P_X(E)\right|_{\op}= O_P\left(\mu_N/\sqrt{NT}\right)$.
\end{lemma}
\noindent {\bf Proof.}  Let $|\cdot|$ denote the $\ell_2$ or operator norm. We use that, due to Assumption \ref{ER} \eqref{ERii}, w.p.a. 1, $\left|P_X(E)\right|=\left|X(X^{\top}X)^{-1}X^{\top}e\right|$ and 
\begin{align*}
\left|X(X^{\top}X)^{-1}X^{\top}e\right|&= \left|\sum_{k=1}^KX_k\left((X^{\top}X)^{-1}X^{\top}e\right)_k\right|
\le \sqrt{\sum_{k=1}^K\left|X_k\right|^2}
\left|(X^{\top}X)^{-1}X^{\top}e\right|_2.
\end{align*}
Due to Assumption \ref{ER} \eqref{ERii} and \eqref{ERiii}, we have 
\begin{align}
\left|(X^{\top}X)^{-1}X^{\top}e\right|_2&\le 
\left|\left(\frac{X^{\top}X}{NT}\right)^{-1}\right|_{\op}\left|\frac{X^{\top}e}{NT}\right|_2=O_P\left(\frac{1}{\sqrt{NT}}\right)\label{ref:consbeta}
\end{align} 
and $\left|X_k\right|_2=\sqrt{(X^\top X)_{kk}}=O_P\left(\sqrt{NT}\right)$ hence the result.
 \hfill $\Box$
 
\noindent By Lemma \ref{tech1} and the inverse triangle inequality, we have
$$\left|\frac{\left|M_X(E)\right|_2}{\sqrt{NT}}-\frac{\left|E\right|_2}{\sqrt{NT}}\right|\le   \frac{\left|P_X(E)\right|_2}{\sqrt{NT}} \xrightarrow{\mathbb{P}} 0$$
and we conclude by Assumption \ref{ER} \eqref{ERi}. %this implis that $\frac{\left|M_X(E)\right|_2}{\sqrt{NT}}\xrightarrow{\mathbb{P}} \sigma$. 
For the operator norm, we use Assumption \ref{ER} \eqref{ERiv} and
$$\left|\norm{M_X(E)}_{\op}-\norm{E}_{\op}\right|\le\norm{P_X(E)}_{\op}.$$

\subsubsection*{Proof of Proposition \ref{propLB2}}
\noindent Let us consider a cone with constant $c$. We work for any draw of $X$ and $\Gamma^l$ and consider the matrices fixed. 
By the computations in the proof of Lemma \ref{tech1}, 
$$\left|P_X(\Delta)\right|_2\le
\frac{\left|X\right|_2}{NT}
\left|\left(\frac{X^{\top}X}{NT}\right)^{-1}\right|_{\op}\left|X^{\top}\delta\right|_2
.$$
Also, for $k\in\{1,\dots,K\}$, using the cone and the trace duality in the third display, we obtain
\begin{align*}
\left|\left\langle X_k,\Delta\right\rangle\right|&\le\left|\left\langle X_k,\mathcal{P}_{\Gamma^l}\left(\Delta\right)\right\rangle\right|+\left|\left\langle X_k,\mathcal{P}_{\Gamma^l}^{\perp}\left(\Delta\right)\right\rangle\right|\\
&=\left|\left\langle \mathcal{P}_{\Gamma^l}\left(X_k\right),\mathcal{P}_{\Gamma^l}\left(\Delta\right)\right\rangle\right|+\left|\left\langle \mathcal{P}_{\Gamma^l}^{\perp}\left(X_k\right),\mathcal{P}_{\Gamma^l}^{\perp}\left(\Delta\right)\right\rangle\right|\\
&\le\min\left(\left|\mathcal{P}_{\Gamma^l}\left(X_k\right)\right|_{\op},\left|X_k\right|_{\op}\right)
\left|\mathcal{P}_{\Gamma^l}\left(\Delta\right)\right|_*
+\left|\mathcal{P}_{\Gamma^l}^{\perp}\left(X_k\right)\right|_{\op}\left|\mathcal{P}_{\Gamma^l}^{\perp}\left(\Delta\right)\right|_*,
\end{align*}
hence
\begin{align*}
\left|P_X(\Delta)\right|_2^2
&\le \sum_{k=1}^K
\left(b_k\left|\mathcal{P}_{\Gamma^l}\left(\Delta\right)\right|_*+b_{\perp k}\left|\mathcal{P}_{\Gamma^l}^\perp\left(\Delta\right)\right|_*\right)^2.
\end{align*}
Also, by homogeneity, we have
$$\kappa_{\Gamma^l,c}^2=2\text{rank}\left(\Gamma^l\right)\inf_{\Delta\in C_{\Gamma^l}:\ \left|\mathcal{P}_{\Gamma^l}\left(\Delta\right)\right|_*=1}\left(\norm{\Delta}^2 - \left|P_X(\Delta)\right|_2^2\right).$$
Denote by $\{\sigma_k\}$ and $\left\{\sigma_{\perp k}\right\}$ the singular values of $\mathcal{P}_{\Gamma^l}\left(\Delta\right)$ and $\mathcal{P}_{\Gamma^l}^{\perp}\left(\Delta\right)$. The rank of the first (resp. the second) matrix  is at most 
$2\text{rank}\left(\Gamma^l\right)$ (resp. 
$p_N$) so, %. Hence, %from the above and 
by the Pythagorean theorem, %we have
\begin{align}
\kappa_{\Gamma^l,c}^2&\ge 2\rank\left(\Gamma^l\right)\inf_{\substack{\sum_k \sigma_k=1\\ 
|\sigma|_0\le2\text{rank}\left(\Gamma^l\right)\\
\sum_k \sigma_{\perp k}\le c\\
\left|\sigma_\perp\right|_0\le p_N\\
\sigma\ge0,\sigma_\perp\ge0}}
\left(\sum_k \sigma_k^2+\sum_k \sigma_{\perp k}^2- \sum_{k=1}^K
\left(b_k+b_{\perp k}\left(\sum_k \sigma_{\perp k}\right)\right)^2\right)\notag\\
&= 1+2\rank\left(\Gamma^l\right)
 \inf_{\substack{\sum_k \sigma_{\perp k}\le c\\
 \left|\sigma_\perp\right|_0\le p_N\\
\sigma_\perp\ge0}}
\left(\sum_k \sigma_{\perp k}^2- \sum_{k=1}^K
\left(b_k+b_{\perp k}\left(\sum_k \sigma_{\perp k}\right)\right)^2\right)\label{bebete}\\
&= 1+2\rank\left(\Gamma^l\right)\inf_{0\le u\le c} \inf_{\substack{\sum_k \sigma_{\perp k}=u\\
\left|\sigma_\perp\right|_0\le p_N\\
\sigma_\perp\ge0}}
\left(\sum_k \sigma_{\perp k}^2-\sum_{k=1}^K
\left(b_k+b_{\perp k}u\right)^2\right)\notag\\
&= 1+2\rank\left(\Gamma^l\right)\min_{0\le u\le c} 
\left(\frac{u^2}{p_N}- \sum_{k=1}^K
\left(b_k+b_{\perp k}u\right)^2\right).\label{em}
\end{align}
The degree 2 polynomial in the bracket has a minimum at $u_*$ given by
$u_*\left(1-p_N\left|b_{\perp}\right|_2^2\right)= p_N\langle b_{\perp},b\rangle$. 
%Because we carry restricted minimisation, 
%The minimum is either attained at 0, $u_*$, or $c$.\\
If $p_N\left|b_{\perp}\right|_2^2\ge 1$ then the minimum is at 0 in which case
$\kappa_{\Gamma^l,c}^2\ge1-2\rank\left(\Gamma^l\right)|b|_2^2$, else, if $p_N\langle b_{\perp},b\rangle< c\left(1-p_N\left|b_{\perp}\right|_2^2\right)$
%(\emph{i.e.}, $p_N\left|b_{\perp}\right|_2^2<1-p_N\langle b_{\perp},b\rangle/c$) 
the minimum is at $u_*$ and %we have 
%is 
%$$\frac{\langle b_{\perp},b\rangle^2p_N}{\left(1-p_N\left|b_{\perp}\right|_2^2\right)^2} - \sum_{k=1}^K
%\left(b_k+b_{\perp k}
%\frac{p_N\langle b_{\perp},b\rangle}{\left(1-p_N\left|b_{\perp}\right|_2^2\right)^2}\right)^2$$
%yielding 
\begin{align*}
\kappa_{\Gamma^l,c}^2&\ge1-2\rank\left(\Gamma^l\right)\left( 
\left|b+b_{\perp }
\frac{p_N\langle b_{\perp},b\rangle}{1-p_N\left|b_{\perp}\right|_2^2}\right|_2^2-\frac{p_N\langle b_{\perp},b\rangle^2}{\left(1-p_N\left|b_{\perp}\right|_2^2\right)^2} \right),
%\\
%&=1-2\rank\left(\Gamma^l\right)\left( |b|_2^2+\frac{p_N\langle b_{\perp},b\rangle^2(1+p_N|b_\perp|_2^2)}{\left(1-p_N\left|b_{\perp}\right|_2^2\right)^2} \right).
\end{align*}
else, the minimum is %attained 
at $c$ and %we have
$$\kappa_{\Gamma^l,c}^2\ge1-2\rank\left(\Gamma^l\right)
\left(
\left|b+b_{\perp}c\right|_2^2-\frac{c^2}{p_N}\right).$$

\begin{remark}\label{ra} 
%By inspection of the proof of Theorem \ref{orac}, 
Denoting by $\mathcal{P}_{A,U\times V}^{\perp}$ the operator defined like 
$\mathcal{P}_{A}$ using annihilators which project onto the orthogonal of the vector space spanned by the columns of $A$ and $U$ (resp. $A$ and $V$) for $U$ and $V$ such that the vector spaces have common dimension $r\left(A,U\times V\right)$, noting that to obtain \eqref{finally2} it is enough that
$\widetilde{\Gamma}\mathcal{P}_{\widetilde{\Gamma},U\times V}^{\perp}(\Delta)^{\top}=0$ and $\widetilde{\Gamma}^{\top}\mathcal{P}_{\widetilde{\Gamma},U\times V}^{\perp}(\Delta)=0$, the result of Theorem \ref{orac} holds replacing $\kappa_{\widetilde{\Gamma},c\left(\rho,\widetilde{\rho}\right)}$ by a compatibility constant replacing  $\mathcal{P}_{\widetilde{\Gamma}}^{\perp}$  by $\mathcal{P}_{\widetilde{\Gamma},U\times V}^{\perp}$, $\mathcal{P}_{\widetilde{\Gamma}}$  by $\mathcal{P}_{\widetilde{\Gamma},U\times V}$, everywhere $\rank\left(\widetilde{\Gamma}\right)$ by $r\left(\widetilde{\Gamma},U\times V\right)$, and with an infimum over $U$ and $V$ after the infimum over $\widetilde{\Gamma}$. The freedom over $U$ and $V$ allows to annihilate low-rank components of $X_k$ if it has a component with a factor structure and deliver constants $b_{\perp k}$ which are  $O_P\left(\sqrt{\max(N,T)}\right)$. 
\end{remark} 

\subsubsection*{Proof of Theorem \ref{const1}}\label{sec:th2}
\noindent The first inequalities follow from Theorem \ref{orac} so we only prove \eqref{ratebeta}. 
Due to Assumption \ref{ER} \eqref{ERii}, w.p.a. 1,  $\widehat{\beta}-\beta=\left(X^\top X\right)^{-1}X^\top (\gamma-\widehat{\gamma})+\left(X^\top X\right)^{-1}X^\top e$, also
%Next, we have 
 \begin{align*}\left|X^\top (\gamma-\widehat{\gamma})\right|_2^2%&=\left|\sum_{k=1}^K x_k^\top (\gamma-\widehat{\gamma})\right|_2\\
 &= \sum_{k=1}^K\left\langle X_k,\widehat{\Gamma}-\Gamma\right\rangle^2\le \sum_{k=1}^K \left|X_k\right|_{\op}^2\left|\widehat{\Gamma}-\Gamma\right|_*^2 \quad\text{(by trace duality)},\\
%hence, 
\left|\left(X^\top X\right)^{-1}X^\top (\gamma-\widehat{\gamma})\right|_2&\le \frac{1}{NT}\left|\left(\frac{X^\top X}{NT}\right)^{-1}\right|_{\op}\sqrt{\sum_{k=1}^K \left|X_k\right|_{\op}^2}\left|\widehat{\Gamma}-\Gamma\right|_*.
 \end{align*}
By Assumption \ref{ER} and \eqref{rategamma}, we obtain
$\left|\left(X^\top X\right)^{-1}X^\top (\gamma-\widehat{\gamma})\right|_2=O_P\left(\lambda_N r_{N}\mu_N/(NT)\right)$. 
% Now, using Theorem \ref{const1} and Assumption \ref{ER}, we obtain
%  $$\left|\left(X^\top X\right)^{-1}X^\top (\gamma-\widehat{\gamma})\right|_2=o_P\left(\frac{\lambda r_{N}}{\sqrt{NT}}\right).$$
  Next, by \eqref{ref:consbeta}, we have $\left|\left(X^\top X\right)^{-1}X^\top e\right|_2=O_P(1/\sqrt{NT})$. This yields the result. %because $\left\{\lambda_N\right\}$ goes to infinity.
  
\subsubsection*{Proof of Proposition \ref{l2}}
%Denote again by  $\Delta=\widehat{\Gamma}-\Gamma$. 
\noindent The proof techniques are similar to those in \cite{koltchinskii2011nuclear}. 
We make use of the fact that if $Z\in \partial |\cdot|_*\left(\widetilde{\Gamma}\right)$, \emph{i.e.}, is of the form
$$ Z=\sum_{k=1}^{\rank\left(\widetilde{\Gamma}\right)}u_k\left(\widetilde{\Gamma}\right)v_k\left(\widetilde{\Gamma}\right)^{\top}+M_{u\left(\widetilde{\Gamma}\right)}WM_{v\left(\widetilde{\Gamma}\right)},$$ for $W$ such that $|W|_{\op}\le1$, then
\begin{equation}\label{propSG}
\left\langle\widehat{ Z}- Z,\widehat{\Gamma}-\widetilde{\Gamma}\right\rangle\ge 0
\end{equation}
and, for a well chosen matrix $W$ (see \cite{koltchinskii2011nuclear} page 2308),
\begin{align*}
\left\langle M_{u\left(\widetilde{\Gamma}\right)}WM_{v\left(\widetilde{\Gamma}\right)},
\widetilde{\Gamma}-\widehat{\Gamma}\right\rangle
=-\left|M_{u\left(\widetilde{\Gamma}\right)}\widehat{\Gamma}M_{v\left(\widetilde{\Gamma}\right)}\right|_*=-\left|\mathcal{P}_{\widetilde{\Gamma}}^{\perp}\left(\widetilde{\Gamma}-\widehat{\Gamma}\right)\right|_*.
\end{align*}
Now, by \eqref{eq:apend} and \eqref{propSG}, we obtain
\begin{align}
&\left\langle M_{ X}\left(\Gamma-\widehat{\Gamma}\right),\widetilde{\Gamma}-\widehat{\Gamma}\right\rangle\notag\\
&\le \lambda\widehat{\sigma}\left\langle Z,\widetilde{\Gamma}-\widehat{\Gamma}\right\rangle-\left\langle M_{ X}\left( E\right),\widetilde{\Gamma}-\widehat{\Gamma}\right\rangle\notag\\
%&\le \lambda\widehat{\sigma}\left\langle\left(\sum_{k=1}^{\rank\left(\widetilde{\Gamma}\right)}u\left(\widetilde{\Gamma}\right)_kv\left(\widetilde{\Gamma}\right)_k^{\top}\right)^{\top}\left(\widetilde{\Gamma}-\widehat{\Gamma}\right)\right)- \lambda\widehat{\sigma}\left|M_{u\left(\widetilde{\Gamma}\right)}\widehat{\Gamma}M_{v\left(\widetilde{\Gamma}\right)}\right|_*-\left\langleM_{ X}\left( E\right)^{\top}\left(\widetilde{\Gamma}-\widehat{\Gamma}\right)\right)\\
&\le \lambda\widehat{\sigma}\left|\mathcal{P}_{\widetilde{\Gamma}}\left(\widetilde{\Gamma}-\widehat{\Gamma}\right)\right|_{*}\wedge\left|P_{u\left(\widetilde{\Gamma}\right)}\left(\widetilde{\Gamma}-\widehat{\Gamma}\right)P_{v\left(\widetilde{\Gamma}\right)}\right|_{*}
- \lambda\widehat{\sigma}\left|\mathcal{P}_{\widetilde{\Gamma}}^{\perp}\left(\widetilde{\Gamma}-\widehat{\Gamma}\right)\right|_*-\left\langle M_{ X}\left( E\right),\widetilde{\Gamma}-\widehat{\Gamma}\right\rangle.\label{eref}
\end{align}
We now use
\begin{align}\label{inrem}
2\left\langle M_{ X}\left(\Gamma-\widehat{\Gamma}\right),\widetilde{\Gamma}-\widehat{\Gamma}\right\rangle&=\left|M_{ X}\left(\Gamma-\widehat{\Gamma}\right)\right|_2^2+
\left|M_{ X}\left(\widetilde{\Gamma}-\widehat{\Gamma}\right)\right|_2^2-\left|M_{ X}\left(\Gamma-\widetilde{\Gamma}\right)\right|_2^2
\end{align}
and consider cases 
(1) $\left\langle M_{ X}\left(\Gamma-\widehat{\Gamma}\right),\widetilde{\Gamma}-\widehat{\Gamma}\right\rangle\le0$ and (2) $\left\langle M_{ X}\left(\Gamma-\widehat{\Gamma}\right),\widetilde{\Gamma}-\widehat{\Gamma}\right\rangle
>0$.\\
In case (1), due to \eqref{inrem}, we have
$\left|M_{ X}\left(\Gamma-\widehat{\Gamma}\right)\right|_2^2\le\left|M_{ X}\left(\Gamma-\widetilde{\Gamma}\right)\right|_2^2$, 
hence the result.\\
In case (2), we have 
\begin{align*}
\lambda\widehat{\sigma}\left|\mathcal{P}_{\widetilde{\Gamma}}^{\perp}\left(\widetilde{\Gamma}-\widehat{\Gamma}\right)\right|_*
&\le  \lambda\widehat{\sigma}\left|\mathcal{P}_{\widetilde{\Gamma}}\left(\widetilde{\Gamma}-\widehat{\Gamma}\right)\right|_{*}-\left\langle M_{ X}\left( E\right),\widetilde{\Gamma}-\widehat{\Gamma}\right\rangle,
\end{align*}
thus, because $\rho\lambda\widehat{\sigma}
%\min\left(\widehat{\sigma},\sigma\right)
\ge \norm{M_X(E)}_{\op}$, 
$\widetilde{\Gamma}-\widehat{\Gamma} \in C_{\widetilde{\Gamma}}$.
%also $\left\langleM_{ X}\left( E\right)^{\top}\left(\widetilde{\Gamma}-\widehat{\Gamma}\right)\right)=\left\langle E^{\top}\left(\widetilde{\Gamma}-\widehat{\Gamma}\right)\right)-\left\langleP_{ X}\left( E\right)^{\top}\left(\widetilde{\Gamma}-\widehat{\Gamma}\right)\right)$, hence
%\begin{align*}
%\left|\left\langleM_{ X}\left( E\right)^{\top}\left(\widetilde{\Gamma}-\widehat{\Gamma}\right)\right)\right|\le 
%\left(\left| E\right|_{\op}+\frac{\left|\left\langle X^{\top} E\right)\right|\left| X\right|_{\op}}{\left| X\right|_{F}^2}\right)\left(\left|\mathcal{P}_{\widetilde{\Gamma}}\left(\widetilde{\Gamma}-\widehat{\Gamma}\right)\right|_{*}+
%\left|\mathcal{P}_{\widetilde{\Gamma}}^{\perp}\left(\widetilde{\Gamma}-\widehat{\Gamma}\right)\right|_{*}
%\right)
%\end{align*}
%
%
%$$\left|\left\langleM_{ X}\left( E\right)^{\top}\left(\widetilde{\Gamma}-\widehat{\Gamma}\right)\right)\right|\le 
%\lambda\widehat{\sigma}\delta\left(\left|\mathcal{P}_{\widetilde{\Gamma}}\left(\widetilde{\Gamma}-\widehat{\Gamma}\right)\right|_{*}+
%\left|\mathcal{P}_{\widetilde{\Gamma}}^{\perp}\left(\widetilde{\Gamma}-\widehat{\Gamma}\right)\right|_{*}
%\right).$$
Moreover, by \eqref{inrem} and \eqref{eref}, we have
 \begin{align*}
&\left|M_{ X}\left(\Gamma-\widehat{\Gamma}\right)\right|_2^2+
\left|M_{ X}\left(\widetilde{\Gamma}-\widehat{\Gamma}\right)\right|_2^2+2\lambda\widehat{\sigma}\left|\mathcal{P}_{\widetilde{\Gamma}}^{\perp}\left(\widetilde{\Gamma}-\widehat{\Gamma}\right)\right|_*\\
&\le \left|M_{ X}\left(\Gamma-\widetilde{\Gamma}\right)\right|_2^2+2\lambda\widehat{\sigma}
%\sqrt{\rank\left(\widetilde{\Gamma}\right)}
\left|\mathcal{P}_{\widetilde{\Gamma}}\left(\widetilde{\Gamma}-\widehat{\Gamma}\right)\right|_*
-2\left\langle M_{ X}\left( E\right),\widetilde{\Gamma}-\widehat{\Gamma}\right\rangle\\
%%&\le \left|M_{ X}\left(\Gamma-\widetilde{\Gamma}\right)\right|_2^2+2\lambda\widehat{\sigma}\sqrt{\rank\left(\widetilde{\Gamma}\right)}\left|\mathcal{P}_{\widetilde{\Gamma}}\left(\widetilde{\Gamma}-\widehat{\Gamma}\right)\right|_2
%%-2\left\langle M_{ X}\left( E\right),\widetilde{\Gamma}-\widehat{\Gamma}\right\rangle,
%\end{align*}
%hence %, by Lemma \ref{l0},
%%$$\left|\left\langle M_{X}\left(E\right),M\right\rangle\right|\le \frac{\lambda\widehat{\sigma}}{2}|M_{X}\left(E\right)|_2\left(
%%\sqrt{2\rank\left(\Gamma\right)}\left|\mathcal{P}_{\Gamma}(M)\right|_2+
%%\left|\mathcal{P}_{\Gamma}^{\perp}\left(M\right)\right|_*\right).$$
%\begin{align*}
%&\left|M_{ X}\left(\Gamma-\widehat{\Gamma}\right)\right|_2^2+
%\left|M_{ X}\left(\widetilde{\Gamma}-\widehat{\Gamma}\right)\right|_2^2+2\lambda\widehat{\sigma}\left|\mathcal{P}_{\widetilde{\Gamma}}^{\perp}\left(\widetilde{\Gamma}-\widehat{\Gamma}\right)\right|_*\\
&\le 
 \left|M_{ X}\left(\Gamma-\widetilde{\Gamma}\right)\right|_2^2+2\lambda\widehat{\sigma}
%\sqrt{\rank\left(\widetilde{\Gamma}\right)}
\left|\mathcal{P}_{\widetilde{\Gamma}}\left(\widetilde{\Gamma}-\widehat{\Gamma}\right)\right|_*
+2\rho\lambda\widehat{\sigma}\left(\left|\mathcal{P}_{\widetilde{\Gamma}}\left(\widetilde{\Gamma}-\widehat{\Gamma}\right)\right|_*+
\left|\mathcal{P}_{\widetilde{\Gamma}}^{\perp}\left(\widetilde{\Gamma}-\widehat{\Gamma}\right)\right|_*\right)%\\
% &\quad+2
%NT\sqrt{\widehat{Q}\left(\widetilde{\Gamma}\right)}\lambda\delta\left(
%\sqrt{2\rank\left(\widetilde{\Gamma}\right)}\left|\mathcal{P}_{\widetilde{\Gamma}}\left(\widetilde{\Gamma}-\widehat{\Gamma}\right)\right|_2+
%\left|\mathcal{P}_{\widetilde{\Gamma}}^{\perp}\left(\widetilde{\Gamma}-\widehat{\Gamma}\right)\right|_*\right)
\end{align*}
and, by definition of $\kappa_{\widetilde{\Gamma},c(\rho)}$, 
\begin{align*}
\left|M_{ X}\left(\Gamma-\widehat{\Gamma}\right)\right|_2^2+
\left|M_{ X}\left(\widetilde{\Gamma}-\widehat{\Gamma}\right)\right|_2^2&\le 
 \left|M_{ X}\left(\Gamma-\widetilde{\Gamma}\right)\right|_2^2+2\lambda(1+\rho)\widehat{\sigma}\left|\mathcal{P}_{\widetilde{\Gamma}}\left(\widetilde{\Gamma}-\widehat{\Gamma}\right)\right|_*\\
&\le 
 \left|M_{ X}\left(\Gamma-\widetilde{\Gamma}\right)\right|_2^2+2\lambda(1+\rho)\widehat{\sigma}\frac{\sqrt{2\rank\left(\widetilde{\Gamma}\right)}}{\kappa_{\widetilde{\Gamma},c(\rho)}} \left|M_{ X}\left(\widetilde{\Gamma}-\widehat{\Gamma}\right)\right|_2,
\end{align*}
hence
\begin{equation*}
\frac{1}{NT}\left|M_{ X}\left(\Gamma-\widehat{\Gamma}\right)\right|_2^2\le 
\frac{1}{NT} \left|M_{ X}\left(\Gamma-\widetilde{\Gamma}\right)\right|_2^2
+\frac{2(\lambda(1+\rho)\widehat{\sigma})^2}{NT}\frac{\rank\left(\widetilde{\Gamma}\right)}{\kappa_{\widetilde{\Gamma},c(\rho)}^2} .
\end{equation*}

\subsubsection*{Proof of Proposition \ref{l1}} 
\noindent \eqref{eq:apend} yields, for all $k=1,\dots,\rank\left(\widehat{\Gamma}\right)$,
\begin{align*}
u_k\left(\widehat{\Gamma}\right)^{\top}M_{ X}
\left(\Gamma^l-\widehat{\Gamma}\right)v_k\left(\widehat{\Gamma}\right)&=\lambda\widehat{\sigma}
-u_k\left(\widehat{\Gamma}\right)^{\top}M_{ X}\left(\Gamma^d+E\right)v_k\left(\widehat{\Gamma}\right)\\
&=\lambda\widehat{\sigma}-\left\langle M_{ X}\left(\Gamma^d+ E\right),u_k\left(\widehat{\Gamma}\right)v_k\left(\widehat{\Gamma}\right)^{\top}\right\rangle,\\
&\ge \lambda(1-\rho)\widehat{\sigma}-\left|\Gamma^d\right|_{\op},
\end{align*}
and, by summing the inequalities, 
\begin{equation}\label{euseful}\left\langle \sum_{k=1}^{\rank\left(\widehat{\Gamma}\right)}u\left(\widehat{\Gamma}\right)_kv\left(\widehat{\Gamma}\right)_k^{\top}, P_{u\left(\widehat{\Gamma}\right)}M_{ X}\left(\Gamma^l-\widehat{\Gamma}\right)P_{v\left(\widehat{\Gamma}\right)}\right\rangle \ge  \left(\lambda(1-\rho)\widehat{\sigma}-\left|\Gamma^d\right|_{\op}\right)\rank\left(\widehat{\Gamma}\right),
\end{equation}
thus
 $$\left|P_{u\left(\widehat{\Gamma}\right)}M_{ X}\left(\Gamma^l-\widehat{\Gamma}\right)P_{v\left(\widehat{\Gamma}\right)}\right|_2\ge \left(\lambda(1-\rho)\widehat{\sigma}-\left|\Gamma^d\right|_{\op}\right)\sqrt{\rank\left(\widehat{\Gamma}\right)}.$$
%This yields
%$$\sqrt{\rank\left(\widehat{\Gamma}\right)}\eta\lambda\widehat{\sigma}\le \left|M_{ X}\left(\Gamma^l-\widehat{\Gamma}\right)\right|_2.$$

\subsubsection*{Proposition \ref{c2}} 
\begin{proposition}\label{c2}
Let 
$m=\left(\frac{|X|_{\op}}{NT}\left|\left(\frac{X^\top X}{NT}\right)^{-1}\right|_{\op}\right)^2\left(\sum_{k=1}^K \left|X_k\right|_{\op}^2\right)\left(\rank\left(\Gamma\right)+\rank\left(\widehat{\Gamma}\right)\right),$ we have
\begin{align*}
&\left|P_{ X}\left(\Gamma-\widehat{\Gamma}\right)\right|_2^2
\le\frac{m}
{\left(1-m\right)_+}\left|M_{ X}\left(\Gamma-\widehat{\Gamma}\right)\right|_2^2,\ \left|\Gamma-\widehat{\Gamma}\right|_2^2
\le\left(1+\frac{m}
{\left(1-m\right)_+}
\right)\left|M_{ X}\left(\Gamma-\widehat{\Gamma}\right)\right|_2^2.
\end{align*}
\end{proposition} 
\noindent {\bf Proof.}
 By Theorem C.5 in \cite{giraud2014introduction}, the definition of $P_{ X}$, and the computations in the proof of Theorem \ref{const1}, we have, w.p.a. 1, 
\begin{align*}
\left|P_{ X}\left(\Gamma-\widehat{\Gamma}\right)\right|_2^2
%&\le\left(\frac{| X|_{2,\rank\left(\Gamma-\widehat{\Gamma}\right)}}{| X|_2}\right)^2\left|\Gamma-\widehat{\Gamma}\right|_2^2\\
&\le
\left(\frac{|X|_{\op}}{NT}\left|\left(\frac{X^\top X}{NT}\right)^{-1}\right|_{\op}\right)^2\left(\sum_{k=1}^K \left|X_k\right|_{\op}^2\right)\rank\left(\Gamma-\widehat{\Gamma}\right)\left|\widehat{\Gamma}-\Gamma\right|_2^2\le m\left|\widehat{\Gamma}-\Gamma\right|_2^2.
\end{align*}
We conclude 
%by propositions \ref{l2} and \ref{l1}
%\begin{equation*}%\label{UBrank}
%\rank\left(\Gamma-\widehat{\Gamma}\right)\le 
%\left(1+\frac{9}{2\kappa_{\widetilde{\Gamma}}^2}\right)\rank(\Gamma)
%\end{equation*}
%and 
by the Pythagorean theorem.\hfill $\Box$
%\begin{align*}
%\left|P_{ X}\left(\Gamma-\widehat{\Gamma}\right)\right|_2^2
%&\le\left(\frac{| X|_{\op}}{| X|_2}\right)^2
%\left(\rank\left(\Gamma\right)+\rank\left(\widehat{\Gamma}\right)\right)\left(\left|M_{ X}\left(\Gamma-\widehat{\Gamma}\right)\right|_2^2+\left|P_{ X}\left(\Gamma-\widehat{\Gamma}\right)\right|_2^2\right).
%\end{align*}
%This yields
%\begin{align*}
%\left|P_{ X}\left(\Gamma-\widehat{\Gamma}\right)\right|_2^2
%&\le\frac{\left(\frac{|X|_{\op}}{NT}\left|\left(\frac{X^\top X}{NT}\right)^{-1}\right|_{\op}\right)^2\left(\sum_{k=1}^K \left|X_k\right|_{\op}^2\right)\left(\rank\left(\Gamma\right)+\rank\left(\widehat{\Gamma}\right)\right)}
%{\left(1-\left(\frac{|X|_{\op}}{NT}\left|\left(\frac{X^\top X}{NT}\right)^{-1}\right|_{\op}\right)^2\left(\sum_{k=1}^K \left|X_k\right|_{\op}^2\right)\left(\rank\left(\Gamma\right)+\rank\left(\widehat{\Gamma}\right)\right)\right)_+}\left|M_{ X}\left(\Gamma-\widehat{\Gamma}\right)\right|_2^2.
%\end{align*}
%We obtain the second display by the Pythagorean theorem.

\subsubsection*{Proof of Proposition \ref{p:op}}
\noindent By \eqref{eq:apend}, we have
$\Gamma^l-\widehat{\Gamma}= \sum_{k=1}^K\left(\widehat{\beta}_k-\beta_k\right) X_k -\Gamma^d- E +\lambda_N\widehat{\sigma}\widehat{Z}$, 
hence
$$\left|\Gamma-\widehat{\Gamma}\right|_{\op}\le \left|\widehat{\beta}-\beta\right|_2\sqrt{\sum_{k=1}^K\left|X_k\right|_{\op}^2}+\left|\Gamma^d\right|_{\op}+\left|E\right|_{\op} +\lambda_N\widehat{\sigma}$$
and we conclude using Theorem \ref{const1} and Assumption \ref{Pen} \eqref{Penii}. 

\subsubsection*{Proof of Proposition \ref{p:th}}
\noindent The Weyl's inequality, yields, for $k\in\{1,\dots,\min(N,T)\}$, 
$$ \left|\sigma_k\left(\Gamma^l\right)-\sigma_k\left(\widehat{\Gamma}\right)\right|\le \left|\Gamma^l-\widehat{\Gamma}\right|_{\op}.$$
This implies, for $k\le \rank\left(\Gamma^l\right)$, 
\begin{equation}\label{eq1}
\sigma_k\left(\widehat{\Gamma}\right)\ge \sigma_k\left(\Gamma^l\right)-\left|\Gamma^l-\widehat{\Gamma}\right|_{\op}
\end{equation}
and, for $k>\rank\left(\Gamma^l\right)$, 
\begin{equation}\label{eq2}
\sigma_k\left(\widehat{\Gamma}\right)\le \left|\Gamma^l-\widehat{\Gamma}\right|_{\op}.
\end{equation}
By Assumption \ref{ass:SF} \eqref{ass:SFi} and Proposition \ref{p:op}, we have
$\mathbb{P}\left(\left|\Gamma^l-\widehat{\Gamma}\right|_{\op}
\le  (\rho+1)\lambda_Nh\sigma\right)\rightarrow 1$. By Theorem \ref{const1} and %condition 
$\lambda _N^2r_{N}=o(NT)$, %(implied from $\lambda _N^3r_{N}=O(NT)$), 
we obtain $\mathbb{P}\left(  (\rho+1)\lambda_Nh\sigma < t\right)\rightarrow 1$ and, by \eqref{eq2}, 
\begin{equation}\label{eq11031} \mathbb{P}\left( \forall k>\rank\left(\Gamma^l\right),\ 
t>\sigma_k\left(\widehat{\Gamma}\right) \right)\rightarrow 1.\end{equation}
By Assumption \ref{ass:SF} \eqref{ass:SFii}, we have
$\mathbb{P}\left( \sigma_k\left(\Gamma^l\right)-\left|\Gamma^l-\widehat{\Gamma}\right|_{\op}\le  (\rho+1)\lambda_Nh^3\sigma\right)\rightarrow 1$. By Theorem \ref{const1} and %the condition 
$\lambda _N^2r_{N}=o(NT)$, we obtain $\mathbb{P}\left(  t <  (\rho+1)\lambda_Nh^3\sigma\right)\rightarrow 1$ and, by \eqref{eq1}, 
\begin{equation}\label{eq11032}\mathbb{P}\left(\forall k\le \rank\left(\Gamma^l\right),\ 
t< \sigma_k\left(\widehat{\Gamma}\right)
 \right)\rightarrow 1.\end{equation}
 Combining \eqref{eq11031} and \eqref{eq11032}, we obtain the first result. The other results are obtained similarly.   

\subsubsection*{Proof of Proposition \ref{p:proj}}
\noindent Because
\begin{align*}
\left|M_{v\left(\widehat{\Gamma}^t\right)}-M_{v\left(\Gamma^l\right)}\right|_2^2&=\left|P_{v\left(\widehat{\Gamma}^t\right)}-P_{v\left(\Gamma^l\right)}\right|_2^2=\rank\left(\widehat{\Gamma}^t\right)+\rank\left(\Gamma^l\right)-2\sum_{k=1}^{\rank\left(\Gamma^l\right)}v_k\left(\Gamma^l\right)^{\top}P_{v\left(\widehat{\Gamma}^t\right)}v_k\left(\Gamma^l\right),\\
%\left|M_{v\left(\widehat{\Gamma}^t\right)}-M_{v\left(\Gamma^l\right)}\right|_2^2
&=\rank\left(\widehat{\Gamma}^t\right)-\rank\left(\Gamma^l\right)+2\sum_{k=1}^{\rank\left(\Gamma^l\right)}v_k\left(\Gamma^l\right)^{\top}M_{v\left(\widehat{\Gamma}^t\right)}v_k\left(\Gamma^l\right)\\
\left|\Gamma^l M_{v\left(\widehat{\Gamma}^{t}\right)}\right|_2^2&=\sum_{k=1}^{\rank\left(\Gamma^l\right)}\sigma_k\left(\Gamma^l\right)^2v_k\left(\Gamma^l\right)^{\top} M_{v\left(\widehat{\Gamma}^{t}\right)}v_k\left(\Gamma^l\right),
\end{align*}
the result follows from
\begin{align*}
\left|M_{v\left(\widehat{\Pi}_v^t\right)}-M_{v\left(\Pi_v^l\right)}\right|_2^2&\le\left|\rank\left(\widehat{\Gamma}^t\right)-\rank\left(\Gamma^l\right)\right|+
\frac{2}{\sigma_{\rank\left(\Gamma^l\right)}\left(\Gamma^l\right)^2}\left|\Gamma^l M_{v\left(\widehat{\Gamma}^t\right)}\right|_2^2\\
&\le\left|\rank\left(\widehat{\Gamma}^t\right)-\rank\left(\Gamma^l\right)\right|+
\frac{2}{\sigma_{\rank\left(\Gamma^l\right)}\left(\Gamma^l\right)^2}\left|\Gamma^l -\widehat{\Gamma}^t\right|_{\op}^2\left|M_{v\left(\widehat{\Gamma}^t\right)}\right|_2^2\\
&\le o_P(1)+2r_N\left(\frac{(\rho+1)\lambda_N(h^2+1)\left(\sigma+o_P\left(1\right)\right)}{\sigma_{\rank\left(\Gamma^l\right)}\left(\Gamma^l\right)} \right)^2.
\end{align*}

\subsubsection*{Proof of Theorem \ref{twostep}}
\noindent Using that $M_{u\left(\widehat{\Pi}_u^t\right)}$ and $M_{v\left(\widehat{\Pi}_v^t\right)}$ are self-adjoint, a solution to \eqref{ewan} satisfies, for $l=1,\dots,K$, 
$\left\langle M_{u\left(\widehat{\Pi}_u^t\right)}X_lM_{v\left(\widehat{\Pi}_v^t\right)} ,Y-\sum_{k=1}^K\widetilde{\beta}^{(1)}_kX_k\right\rangle=0$, hence
\begin{align*}
&\left\langle M_{u\left(\Pi_u^l\right)}X_lM_{v\left(\Pi_v^l\right)} ,\Gamma^d+E+\sum_{k=1}^K\left(\beta_k-\widetilde{\beta}^{(1)}_k\right)X_k \right\rangle\\
&=\left\langle\left(M_{u\left(\Pi_u^l\right)}-M_{u\left(\widehat{\Pi}_u^t\right)}\right)X_lM_{v\left(\Pi_v^l\right)} ,\Gamma^d+E+\sum_{k=1}^K\left(\beta_k-\widetilde{\beta}^{(1)}_k\right)X_k\right\rangle\\
&\quad+\left\langle M_{u\left(\Pi_u^l\right)}X_l\left(M_{v\left(\Pi_v^l\right)}-M_{v\left(\widehat{\Pi}_v^t\right)}\right) ,\Gamma^d+E+\sum_{k=1}^K\left(\beta_k-\widetilde{\beta}^{(1)}_k\right) X_k\right\rangle\\
&\quad-\left\langle \left(M_{u\left(\Pi_u^l\right)}-M_{u\left(\widehat{\Pi}_u^t\right)}\right) X_l\left(M_{v\left(\Pi_v^l\right)}-M_{v\left(\widehat{\Pi}_v^t\right)}\right) ,\Gamma+E+\sum_{k=1}^K\left(\beta_k-\widetilde{\beta}^{(1)}_k\right) X_k \right\rangle, 
\end{align*}
so
\begin{align}
&\sum_{k=1}^K\left(\beta_k-\widetilde{\beta}^{(1)}_k\right)\Bigg(
\left\langle M_{u\left(\Pi_u^l\right)}X_lM_{v\left(\Pi_v^l\right)} ,X_k\right\rangle
-\left\langle\left(M_{u\left(\Pi_u^l\right)}-M_{u\left(\widehat{\Pi}_u^t\right)}\right)X_lM_{v\left(\Pi_v^l\right)} ,X_k \right\rangle\notag\\
& \quad 
-\left\langle M_{u\left(\Pi_u^l\right)}X_l\left(M_{v\left(\Pi_v^l\right)}-M_{v\left(\widehat{\Pi}_v^t\right)}\right) ,X_k \right\rangle\notag\\
& \quad +\left\langle \left(M_{u\left(\Pi_u^l\right)}-M_{u\left(\widehat{\Pi}_u^t\right)}\right) X_l\left(M_{v\left(\Pi_v^l\right)}-M_{v\left(\widehat{\Pi}_v^t\right)}\right) ,X_k \right\rangle\Bigg)\notag\\
&=-\left\langle M_{u\left(\Pi_u^l\right)}X_lM_{v\left(\Pi_v^l\right)} ,\Gamma^d+E\right\rangle+\left\langle\left(M_{u\left(\Pi_u^l\right)}-M_{u\left(\widehat{\Pi}_u^t\right)}\right)X_lM_{v\left(\Pi_v^l\right)} ,\Gamma^d+E\right\rangle\notag\\
&\quad+\left\langle M_{u\left(\Pi_u^l\right)}X_l\left(M_{v\left(\Pi_v^l\right)}-M_{v\left(\widehat{\Pi}_v^t\right)}\right) ,\Gamma^d+E\right\rangle\notag\\
&\quad -\left\langle \left(M_{u\left(\Pi_u^l\right)}-M_{u\left(\widehat{\Pi}_u^t\right)}\right) X_l\left(M_{v\left(\Pi_v^l\right)}-M_{v\left(\widehat{\Pi}_v^t\right)}\right),\Gamma+E\right\rangle\label{emain}
%+\left\langle \left(M_{u\left(\Pi_u^l\right)}-M_{u\left(\widehat{\Pi}_u^t\right)}\right) X_l\left(M_{v\left(\Pi_v^l\right)}-M_{v\left(\widehat{\Pi}_v^t\right)}\right) ,\Gamma \right\rangle
.
\end{align}
Let us show that $\left\langle M_{u\left(\Pi_u^l\right)}X_lM_{v\left(\Pi_v^l\right)} ,X_k\right\rangle$, which by Assumption \ref{ass:asnorm} \eqref{ass:asnormv} diverges like $NT$, is the high-order term %in the bracket 
multiplying $\left(\beta_k-\widetilde{\beta}^{(1)}_k\right)$ in \eqref{emain}. This also yields the consistency of the estimator of the %variance-
covariance matrix.  
By self-adjointness of the projections, Theorem C.5 in \cite{giraud2014introduction}, and Proposition \ref{p:th} with the modifications of Section \ref{sec:modif} which imply $\rank\left(M_{u\left(\Pi_u^l\right)}-M_{u\left(\widehat{\Pi}_u^t\right)}\right)\le 2\overline{r}_N$ w.p.a. 1,  denoting, for a matrix $M$ and $r\in\N$ by $|M|_{2,r}^2=\sum_{k=1}^{r}\sigma_k(M)^2$, we have, 
\begin{align*}
&\left|\left\langle\left(M_{u\left(\Pi_u^l\right)}-M_{u\left(\widehat{\Pi}_u^t\right)}\right)X_lM_{v\left(\Pi_v^d\right)} ,X_k \right\rangle\right|\\
&\le (1+o_P(1))
\left|M_{u\left(\Pi_u^l\right)}-M_{u\left(\widehat{\Pi}_u^t\right)}\right|_2\left|X_l M_{v\left(\Pi_u^l\right)}X_k^\top\right|_{2,2r_N}\\
%\left|\left(M_{u\left(\Pi_u^l\right)}-M_{u\left(\widehat{\Pi}_u^t\right)}\right)X_lM_{v\left(\Pi_v^d\right)}\right|_*\left|X_kM_{v\left(\Pi_v^d\right)}\right|_{\op}\\ 
%&\le \left(\sqrt{2\overline{r}_N}+o_P(1)\right)\left|\left(M_{u\left(\Pi_u^l\right)}-M_{u\left(\widehat{\Pi}_u^t\right)}\right)X_lM_{v\left(\Pi_v^l\right)}\right|_2\left|X_kM_{v\left(\Pi_v^l\right)}\right|_{\op}\\
&\le \left(\sqrt{2\overline{r}_N}+o_P(1)\right)\left|M_{u\left(\Pi_u^l\right)}-M_{u\left(\widehat{\Pi}_u^t\right)}\right|_2\left|X_lM_{v\left(\Pi_v^l\right)}\right|_{\op}\left|X_kM_{v\left(\Pi_v^l\right)}\right|_{\op},
\end{align*}
%\blue{(on peut peut etre faire un peu plus fin avec $P_{u,\widehat{u}}$)}
hence, by Proposition \ref{p:proj} with the modifications of Section \ref{sec:modif}, 
\begin{align*}
&\left|\left\langle\left(M_{u\left(\Pi_u^l\right)}-M_{u\left(\widehat{\Pi}_u^t\right)}\right)X_lM_{v\left(\widehat{\Pi}_v^t\right)} ,X_k \right\rangle\right|
\\
&
\le 
\frac{2 (\rho+1)\overline{r}_N\lambda_N}{\overline{v}_N}\left((h^2+1)\widetilde{\sigma}+o_P\left(1\right)\right)\left|\left(\Pi_l^d+U_l\right)M_{v\left(\Pi_v^l\right)}\right|_{\op}\left|\left(\Pi_k^d+U_k\right)M_{v\left(\Pi_v^l\right)}\right|_{\op}.
\end{align*}
We treat similarly $\left|\left\langle M_{u\left(\Pi_u^l\right)}X_l\left(M_{v\left(\Pi_v^l\right)}-M_{v\left(\widehat{\Pi}_v^t\right)}\right) ,X_k\right\rangle\right|$,
and, for the fourth term, use that
%, using that $\rank\left(M_{u\left(\Pi_u^l\right)}-M_{u\left(\widehat{\Pi}_u^t\right)}\right)\le 2r_N$ w.p.a. 1 and Lemma C2 item 1 in \cite{giraud2014introduction},
\begin{align*}
&\left|\left\langle \left(M_{u\left(\Pi_u^l\right)}-M_{u\left(\widehat{\Pi}_u^t\right)}\right) X_l\left(M_{v\left(\Pi_v^l\right)}-M_{v\left(\widehat{\Pi}_v^t\right)}\right) ,X_k \right\rangle
\right|\\
&\le \left| \left(M_{u\left(\Pi_u^l\right)}-M_{u\left(\widehat{\Pi}_u^t\right)}\right) X_l\left(M_{v\left(\Pi_v^l\right)}-M_{v\left(\widehat{\Pi}_v^t\right)}\right)\right|_* \left|X_k \right|_{\op}\\
&\le  \left(\sqrt{2\overline{r}_N}+o_P(1)\right)\left| \left(M_{u\left(\Pi_u^l\right)}-M_{u\left(\widehat{\Pi}_u^t\right)}\right) X_l\left(M_{v\left(\Pi_v^l\right)}-M_{v\left(\widehat{\Pi}_v^t\right)}\right)\right|_2 \left|X_k \right|_{\op}\\
&\le  \left(\sqrt{2\overline{r}_N}+o_P(1)\right)\left|M_{u\left(\Pi_u^l\right)}-M_{u\left(\widehat{\Pi}_u^t\right)}\right|_{\op}\left|X_l\left(M_{v\left(\Pi_v^l\right)}-M_{v\left(\widehat{\Pi}_v^t\right)}\right)\right|_2 \left|X_k \right|_{\op}\\
&\le  \left(\sqrt{2\overline{r}_N}+o_P(1)\right)\left|M_{u\left(\Pi_u^l\right)}-M_{u\left(\widehat{\Pi}_u^t\right)}\right|_{2}\left|X_l\right|_{\op}\left|M_{v\left(\Pi_v^l\right)}-M_{v\left(\widehat{\Pi}_v^t\right)}\right|_2 \left|X_k \right|_{\op}\\
&\le\frac{ (\rho+1)^2(2\overline{r}_N)^{3/2}\lambda_N^2}{\overline{v}_N^2}\left((h^2+1)^2\widetilde{\sigma}^2+o_P\left(1\right)\right)|X_l|_{\op}|X_k|_{\op},
\end{align*}
where we use Proposition \ref{p:th}  in the third display and Proposition \ref{p:proj} (with the modifications of Section \ref{sec:modif}) in the last display. 
%\begin{align*}
%&\left|
%\left\langle \left(M_{u\left(\Pi_u^l\right)}-M_{u\left(\widehat{\Pi}_u^t\right)}\right) X_l\left(M_{v\left(\Pi_v^l\right)}-M_{v\left(\widehat{\Pi}_v^t\right)}\right) ,X_k \right\rangle
%\right|\le 
%\frac{9\sqrt{2}r_N^{3/2}\lambda_N^2}{2v_N^2}\left((c^2+1)^2\sigma^2+o_P\left(1\right)\right)|X_l|_{\op}|X_k|_{\op}.
%\end{align*}
%\red{
%Je dirais plutôt 
%$$O_P(2r_N)\frac{9r_N\lambda_N^2}{2v_N^2}\left((c^2+1)^2\sigma^2+o_P\left(1\right)\right)|X_l|_{\op}|X_k|_{\op}.$$
%}
% is of even smaller order \red{(Pas forcément car il y a les projecteurs qui donnent $\left|M_{u(\Gamma)}X_l\right|_{\op}$ et $\left|X_lM_{v\left(\Pi_v^l\right)}\right|_{\op}$ dans les deux autres termes, si jamais $\left|M_{u(\Gamma)}X_l\right|_{\op}$ et $\left|X_lM_{v\left(\Pi_v^l\right)}\right|_{\op}$ sont négligeables par rapport à $\left|X_l\right|_{\op}$, alors le dernier terme n'est pas forcément plus petit je pense)}. \blue{(le faire alors, j'ai commente dans un email)} \red{(On trouve $(9r_N^2\lambda_N^2\left|X_l\right|_{\op}/v_N^2)\left((c^2+1)^2\sigma^2+o_P\left(1\right)\right)\left|X_k\right|_{\op}$, no ?)}\\
Let us consider now the quantities on the right-hand side in \eqref{emain}. Proceeding like above,
% and using \eqref{rateopt}, 
we have
\begin{align*}
&\left|\left\langle\left(M_{u\left(\Pi_u^l\right)}-M_{u\left(\widehat{\Pi}_u^t\right)}\right)X_lM_{v\left(\Pi_v^l\right)} ,\Gamma^d+E \right\rangle\right|\notag\\
&\le (1+o_P(1))
\left|M_{u\left(\Pi_u^l\right)}-M_{u\left(\widehat{\Pi}_u^t\right)}\right|_2\left|X_l M_{v\left(\Pi_v^l\right)}\left(\Gamma^d+E\right)^\top\right|_{2,2r_N}\notag\\
&\le 
\frac{2(\rho+1)\overline{r}_N\lambda_N\left((h^2+1)\widetilde{\sigma}+o_P\left(1\right)\right)
}{\overline{v}_N}\left(\rho\lambda_N\sigma+\left|\Gamma^dM_{v\left(\Pi_v^l\right)}\right|_{\op}\right)
\left(\rho\lambda_N\sigma_l+\left|\Pi_l^dM_{v\left(\Pi_v^l\right)}\right|_{\op}\right)%\label{eder}
\end{align*}
and treat similarly $\left\langle M_{u\left(\Pi_u^l\right)}X_l\left(M_{v\left(\Pi_v^l\right)}-M_{v\left(\widehat{\Pi}_v^t\right)}\right) ,\Gamma^d+E\right\rangle$.
%\blue{(Du fait de ce terme (et de notre facon (trop) simple de le traiter, un seul annihilateur aurait marche aussi bien n'est-ce pas?)}
%\red{(Je trouve plutôt $$ \frac{2\sqrt{2r_N}\lambda_N}{v_N}\max\left(\left|X_lM_{v\left(\Pi_v^l\right)}\right|_{\op},\left|M_{u(\Gamma)}X_l\right|_{\op}\right)\left(\frac{5}{2}\lambda_N(\sigma+o_P(1))\right).$$)}
%\blue{(merci ca m'a permis de corriger a plein d'endroits, par contre tu dois obtenir quelque chose homogene a du $\sigma^2$ car il y a un produit de 2 quantites stochastiques)}  \red{Effectivement, }
With the same arguments, the absolute value of the last term of \eqref{emain} is smaller than\\
$$\frac{(\rho+1)\sqrt{2}(2\overline{r}_N)^{3/2}\lambda_N^2\left((h^2+1)^2\widetilde{\sigma}^2+o_P\left(1\right)\right)}{\overline{v}_N^2}\left|X_l\right|_{\op}
\left(\left|\Gamma\right|_{\op}+ \rho\lambda_N(h^2+1)\widetilde{\sigma}+o_P(1)\right).$$
%\red{
%Je dirais plutôt 
%\begin{align*}
%&\left|\left\langle\left(M_{u\left(\Pi_u^l\right)}-M_{u\left(\widehat{\Pi}_u^t\right)}\right)X_l\left(M_{v\left(\Pi_v^l\right)}-M_{v\left(\widehat{\Pi}_v^t\right)}\right) ,\Gamma+E \right\rangle\right|\\
%&\le \left|\left(M_{u\left(\Pi_u^l\right)}-M_{u\left(\widehat{\Pi}_u^t\right)}\right)X_l\left(M_{v\left(\Pi_v^l\right)}-M_{v\left(\widehat{\Pi}_v^t\right)}\right)\right|_*\left(\left|\Gamma\right|_{\op}+\left|E \right|_{\op}\right)\\
%&\le \text{rank}\left(M_{u\left(\Pi_u^l\right)}-M_{u\left(\widehat{\Pi}_u^t\right)}\right)\left|\left(M_{u\left(\Pi_u^l\right)}-M_{u\left(\widehat{\Pi}_u^t\right)}\right)X_l\left(M_{v\left(\Pi_v^l\right)}-M_{v\left(\widehat{\Pi}_v^t\right)}\right)\right|_{\op}\left(\left|\Gamma\right|_{\op}+\left|E \right|_{\op}\right)\\
%&=O_P(2r_N)\left|M_{u\left(\Pi_u^l\right)}-M_{u\left(\widehat{\Pi}_u^t\right)}\right|_{\text{F}} \left|X_l\right|_{\op}\left|M_{v\left(\Pi_v^l\right)}-M_{v\left(\widehat{\Pi}_v^t\right)}\right|_{\text{F}}\left(\left|\Gamma\right|_{\op}+\left|E \right|_{\op}\right)\\
%&\le O_P(2r_N)\frac{3\sqrt{2r_N}\lambda_N}{2v_N}\left((c^2+1)\sigma+o_P\left(1\right)\right)\mu_N\frac{3\sqrt{2r_N}\lambda_N}{2v_N}\left((c^2+1)\sigma+o_P\left(1\right)\right)\left(\left|\Gamma\right|_{\op}+\frac12\lambda_N\left(\sigma+o_P\left(1\right)\right)\right)\\
%&=O_P(r_N^2)\frac{9\lambda_N^2}{v_N^2}\left((c^2+1)^2\sigma^2+o_P\left(1\right)\right)\mu_N\left(\left|\Gamma\right|_{\op}+\frac12\lambda_N\left(\sigma+o_P\left(1\right)\right)\right).
%\end{align*}
%}
Let us now look at the first terms on the left-hand side and on the right-hand side of \eqref{emain}. By \eqref{ass:asnormiv}, %$\lambda_N^2\max(r_{kN},r_{lN})=o(NT/\sqrt{\overline{r}}_N)=o(NT)$ so, 
for all $k,l\in\{1,\dots,K\}$, 
$$\left\langle M_{u\left(\Pi_u^l\right)}X_lM_{v\left(\Pi_v^l\right)} ,X_k\right\rangle=\left\langle M_{u\left(\Pi_u^l\right)}U_lM_{v\left(\Pi_v^l\right)} ,U_k\right\rangle+o_P(NT)$$
so,  by \eqref{ass:asnormv}, $\left\langle M_{u\left(\Pi_u^l\right)}X_lM_{v\left(\Pi_v^l\right)} ,X_k\right\rangle$ are the high-order terms on the left-hand side of \eqref{emain}. Similarly, by \eqref{ass:asnormiv}, the high-order terms on the right-hand side of \eqref{emain} are $\left\langle M_{u\left(\Pi_u^l\right)}U_lM_{v\left(\Pi_v^l\right)} ,E\right\rangle$. 
%\red{(Rather, $(9r_N^2\lambda_N^2\left|X_l\right|_{\op}/v_N^2)\left((c^2+1)^2\sigma^2+o_P\left(1\right)\right)(\left|\Gamma\right|_{\op}+\left|E\right|_{\op})$, no ?)}
%\blue{(bien vu, e propose encore autre chose qu vu de ta remarque precedente, ca te vas?)} \red{Je ne comprends pas pourquoi tu as $r_N^{3/2}$ seuleument, pour moi il y a 2 $\sqrt{r_N}$ qui viennent des deux différences et un $r_N$ qui vient de la norme nucléaire. Sinon, je suis d'accord, j'avais oublié des choses car j'avais mal appliqué la proposition 11} 
%\blue{(revoir la suite et mettre les hypotheses qu'il faut dans l'ennonce)}  This is $o_P(\sqrt{NT})$ by Assumption \ref{ass:asnorm} (ii). 
%Moreover, by Assumption \ref{ass:asnorm} (v), we have $\left(\left\langle M_{u\left(\Pi_u^l\right)}X_lM_{v\left(\Pi_v^l\right)} ,E\right\rangle\right)/\sqrt{NT}\xrightarrow{\mathbb{d}}\mathcal{N}\left(0,\sigma\Sigma_{\perp}\right)$, hence this is the high-order term.\\ 
As a result, $\widetilde{\beta}^{(1)}$ is asymptotically equivalent to the ideal estimator $\overline{\beta}$
\begin{equation}\label{ewanideal}
\overline{\beta}\in\argmin{\beta\in\mathbb{R}^K}\left|\mathcal{P}_{\Pi^l}^\perp\left(Y-\sum_{k=1}^K\beta_kU_k\right)\right|_2^2.
\end{equation}
Hence, w.p.a. 1, $\overline{\beta}-\beta=\left(\mathcal{P}_{\Pi^l}^\perp(U)^\top \mathcal{P}_{\Pi^l}^\perp(U)\right)^{-1}\mathcal{P}_{\Pi^l}^\perp(U)^\top e$ and we conclude by usual arguments.\\
To obtain the first part of the second statement we use that $U^\top U-\mathcal{P}_{\Pi^l}^\perp(U)^\top\mathcal{P}_{\Pi^l}^\perp(U)$ is symmetric positive definite. It is clearly symmetric. The positive definiteness follows from the following computations. Let $b\in\mathbb{R}^K$, we have
%Moreover, because $M_{u(\Pi^l)}$ is an orthogonal projector, for all $b\in\mathbb{R}^K$, $\sum_{k=1}^Kb_k^2\left|M_{u(\Pi^l)}X_kM_{v\left(\Pi_v^l\right)}\right|_2^2\le \sum_{k=1}^Kb_k^2\left|X_kM_{v\left(\Pi_v^l\right)}\right|_2^2$. Also, using that $M_{v\left(\Pi_v^l\right)}$ is self-adjoint, we obtain  
%$\sum_{k,l}b_kb_l\tr\left(M_{v\left(\Pi_v^l\right)}X_k^\top X_l\right)=\sum_{k,l}b_kb_l\tr\left(X_k^\top X_l\right)- 
%\sum_{k,l}b_kb_l\tr\left(P_{v\left(\Pi^l\right)}X_k^\top X_l\right)$, hence
\begin{align*}
\sum_{k,l}b_kb_l\tr\left(U_k^\top U_l\right)&\ge
\sum_{k,l}b_kb_l\tr\left(M_{v\left(\Pi_v^l\right)}U_k^\top U_l\right)\\
&=
\sum_{k,l}b_kb_l\tr\left(M_{v\left(\Pi_v^l\right)}U_k^\top M_{u\left(\Pi_v^l\right)}U_lM_{v\left(\Pi_v^l\right)}\right)+\sum_{k,l}b_kb_l\tr\left(M_{v\left(\Pi_v^l\right)}U_k^\top P_{u\left(\Pi_v^l\right)}U_lM_{v\left(\Pi_v^l\right)}\right)\\
&\ge\sum_{k,l}b_kb_l\tr\left(\mathcal{P}_{\Pi^l}^\perp(U_k)^\top
\mathcal{P}_{\Pi^l}^\perp(U_l)\right).\end{align*} Because $U^\top U$ has a fixed dimension, all norms are equivalent and 
$\left|U^\top U-\mathcal{P}_{\Pi^l}^\perp(U)^\top\mathcal{P}_{\Pi^l}^\perp(U)\right|_{\op}\le \tr\left(U^\top U-\mathcal{P}_{\Pi^l}^\perp(U)^\top\mathcal{P}_{\Pi^l}^\perp(U)\right)=\left|\mathcal{P}_{\Pi^l}(U)\right|_2^2=o_P(\left|U\right|_2^2)$. 
We conclude using that $\left|U\right|_2^2\le K\left|U^\top U\right|_{\op}$. Also, from the above, $\mathcal{P}_{\Pi^l}^\perp(U)^\top\mathcal{P}_{\Pi^l}^\perp(U)=\mathcal{P}_{\Pi^l}^\perp(U)^\top\mathcal{P}_{\Pi^l}^\perp(U)+M$ where $M$ is a smaller order term by condition (iv). We obtain the last part of the second statement using the next lemma. 

\begin{lemma}\label{lU} %Let $k\in\{1,\dots,K\}$ and 
Assume 
%$U_k\in\mathcal{M}_{NT}$ has iid mean zero entries with variance $\sigma_k^2$, that 
$U$ and $(\Pi_u^l,\Pi_v^l)$ are independent, and $\mathbb{E}\left[\max\left(\rank\left(\Pi_u^l\right),\rank\left(\Pi_v^l\right)\right)\right]=o\left(\sqrt{\min(N,T)}\right)$, then 
$ \left|\mathcal{P}_{\Pi^l}(U)\right|_{2}^2/(NT)=o_P(1)$, hence
$%\left|
\mathcal{P}^\perp_{\Gamma^r}(U)^\top
\mathcal{P}^\perp_{\Gamma^r}(U)
%\right|_{2}^2
/(NT)\xrightarrow{\mathbb{P}}\mathbb{E}\left[U^\top U\right]$. 
\end{lemma}
\noindent {\bf Proof.} We prove that, for $k\in\left\{1,\dots,K\right\}$, 
% By the Pythagoren theorem, it is enough to show that 
$ \left|\mathcal{P}_{\Pi^l}(U_k)\right|_{2}^2/(NT)$ converges to 0 in $L^1$. This relies on \eqref{eP} and the facts that $M_{u(\Pi^l)}$ is a contraction for the Euclidian norm and 
\begin{align*}\mathbb{E}\left[\left|U_kP_{v(\Pi^l)}\right|_2^2\right]&= \mathbb{E}\left[\mathbb{E}\left[\left|U_kP_{v(\Pi^l)}\right|_2^2\left| \Pi_u^l,\Pi_v^l\right.\right]\right]\\
%&\le \mathbb{E}\left[\mathbb{E}\left[\left|UP_{v(\Gamma^r)}\right|_2^2\vert \Gamma^r\right]\right]\\
&=\mathbb{E}\left[\mathbb{E}\left[\sum_{i=1}^N\left|U_{i\cdot}P_{v(\Pi^l)}\right|_2^2\left| \Pi_u^l,\Pi_v^l\right.\right]\right]
=N\mathbb{E}\left[\rank\left(\Pi_v^l\right)\right]u^2=o(NT)
\end{align*}
and similarly for $\mathbb{E}\left[\left|P_{u(\Pi^l)}U_k\right|_2^2\right]$. By the arguments in the previous proof 
$U^\top U/(NT)$ and 
$%\left|
\mathcal{P}^\perp_{\Gamma^r}(U)^\top
\mathcal{P}^\perp_{\Gamma^r}(U)
%\right|_{2}^2
/(NT)$ have same limit, hence the result by the law of large numbers.
\hfill $\Box$

%\subsection*{Proof of Theorem \ref{twostep3}}
%\noindent The main difference is in the treatment of \eqref{eder}. We have
%\begin{align*}
%&\left|\left\langle\left(M_{u\left(\Gamma^l\right)}-M_{u\left(\widehat{\Gamma}^t\right)}\right)X_lM_{v\left(\Pi_v^l\right)} ,E \right\rangle\right|\\
%&\le (1+o_P(1))\left|P_{u\left(\Gamma^l\right)}-P_{u\left(\widehat{\Gamma}^l\right)}\right|_2\left|
%X_lM_{v\left(\Gamma^l\right)}E^\top\right|_{2,2r_N}
%\end{align*}
%and, given  $\left(M_{u\left(\Gamma^l\right)},M_{u\left(\Gamma^l\right)},X\right)$, the columns of $X_lM_{v\left(\Gamma^l\right)}E^\top$ are independent mean zero Gaussian vectors in $\mathbb{R}^N$ with covariance matrix $X_lM_{v\left(\Gamma^l\right)}X_l^\top$. []
% 
%and the second random matrix has independent columns, on pourrait utiliser une hypothese du type $(X,\Gamma)$ independant de $E$. Ce terme ne semble pas leur poser de probleme. On pourrait travailler avec les Grassmaniennes aussi).

\end{document}